# INTRODUCTION TO BIMATRICES

W. B. Vasantha Kandasamy
Florentin Smarandache
K. Ilanthenral

# INTRODUCTION TO BIMATRICES

**W. B. Vasantha Kandasamy**
Department of Mathematics
Indian Institute of Technology, Madras
Chennai – 600036, India
e-mail: vasantha@iitm.ac.in
web: http://mat.iitm.ac.in/~wbv

**Florentin Smarandache**
Department of Mathematics
University of New Mexico
Gallup, NM 87301, USA
e-mail: smarand@gallup.unm.edu

**K. Ilanthenral**
Editor, Maths Tiger, Quarterly Journal
Flat No.11, Mayura Park,
16, Kazhikundram Main Road, Tharamani,
Chennai – 600 113, India
e-mail: ilanthenral@gmail.com

# CONTENTS

Chapter Three

## INTRODUCTION TO NEUTROSOPHIC BIMATRICES

Chapter Four

## SUGGESTED PROBLEMS 

## BIBLIOGRAPHY 

## INDEX 

## ABOUT THE AUTHORS

## Preface

Matrix theory has been one of the most utilised concepts in fuzzy models and neutrosophic models. From solving equations to characterising linear transformations or linear operators, matrices are used. Matrices find their applications in several real models. In fact it is not an exaggeration if one says that matrix theory and linear algebra (i.e. vector spaces) form an inseparable component of each other.

The study of bialgebraic structures led to the invention of new notions like birings, Smarandache birings, bivector spaces, linear bialgebra, bigroupoids, bisemigroups, etc. But most of these are abstract algebraic concepts except, the bisemigroup being used in the construction of biautomatons. So we felt it is important to construct non-abstract bistructures which can give itself for more and more lucid applications.

So, as a first venture we have defined the notion of bimatrices. These bimatrices in general cannot be given an algebraic operation like bimatrix addition and bimatrix multiplication, so that the collection of bimatrices become closed with respect to these operations. In fact this property will be nice in a way, for, in all our analysis we would not in general get a solution from a set we have started with. Only this incompleteness led to the invention of complex number or the imaginary number "i". (Complex number i = $\sqrt{-1}$ ).

But, after the advent of neutrosophic logic, it is more appropriate to accept the notion of indeterminacy. It is common for any one to say that a particular thing is indeterminate than saying it is imaginary; more so in mathematics. Because when one says something is imaginary, to a common mind it means that the thing does not exist, but indeterminacy is acceptable, for we can say to any layman, "We cannot determine what you ask for", but we cannot say, "what you ask for is imaginary". Thus when we are in a meek position as we cannot perceive -1 all the more it is difficult to perceive the square root of –1, so we felt it is appropriate under these circumstances to introduce the notion of indeterminacy I where $I^2 = I$. Using this indeterminacy we construct neutrosophic bimatrices, which can be used in neutrosophic models.

By no means we make any claim that the concept of indeterminacy would be an alternative of the imaginary. We have already defined the concept of neutrosophic matrices and have used them in Neutrosophic Cognitive Maps model and in the Neutrosophic Relational Equations models, which are analogous to Fuzzy Cognitive Map and Fuzzy Relational Equations models respectively. These bimatrices, fuzzy bimatrices, or neutrosophic bimatrices will certainly find their application when the study of the model is time-dependent and where relative comparison is needed.

Thus this book, which is an introduction to bimatrices and neutrosophic bimatrices contains 4 chapters. First chapter introduces the notion of bimatrices and analyses its properties. Chapter 2 introduces bivector spaces and defines notions like bieigen vectors, bipolynomials etc. In Chapter 3 neutrosophic bimatrices and fuzzy bimatrices are introduced. The final chapter gives some problems, which will make the reader to understand these new notions. We have given a possible list of books for the reader to understand background.

We thank Dr. Kandasamy for proof reading and Meena Kandasamy for cover designing and Kama Kandasamy for formatting this book.

W.B.VASANTHA KANDASAMY
FLORENTIN SMARANDACHE
K. ILANTHENRAL

**Chapter One**

# BIMATRICES AND THEIR PROPERTIES

Matrices provide a very powerful tool for dealing with linear models. Bimatrices which we are going to define in this chapter are still a powerful and an advanced tool which can handle over one linear model at a time. Bimatrices will be useful when time bound comparisons are needed in the analysis of the model.

Thus for the first time the notion of bimatrices are introduced. Unlike matrices bimatrices can be of several types.

This chapter introduces bimatrices and illustrates them with examples. This chapter has seven sections. In section one we just define the notion of bimatrix and give some basic operations on them. In section two we give some basic properties of bimatrices.

Symmetric and skew symmetric bimatrices are introduced in section three. Section four introduces the concept of subbimatrix and some related results. The basic notion of a determinant in case of bimatrices is introduced in section five. Laplace equation for a bimatrix is dealt in section six and the final section defines the notion of overlap bimatrices.

Throughout this book Z will denote the set of integers both positive, negative with zero, R the reals, Q the set of rationals and C the set of complex numbers. $Z_n$ will denote the set of integers modulo n.

## 1.1 Definition and basic operations on bimatrices

In this section we for the first time define the notion of bimatrix and illustrate them with examples and define some of basic operation on it.

**DEFINITION 1.1.1:** *A bimatrix $A_B$ is defined as the union of two rectangular away of numbers $A_1$ and $A_2$ arranged into rows and columns. It is written as follows $A_B = A_1 \cup A_2$ where $A_1 \neq A_2$ with*

$$A_1 = \begin{bmatrix} a_{11}^1 & a_{12}^1 & \cdots & a_{1n}^1 \\ a_{21}^1 & a_{22}^1 & \cdots & a_{2n}^1 \\ \vdots & & & \\ a_{m1}^1 & a_{m2}^1 & \cdots & a_{mn}^1 \end{bmatrix}$$

*and*

$$A_2 = \begin{bmatrix} a_{11}^2 & a_{12}^2 & \cdots & a_{1n}^2 \\ a_{21}^2 & a_{22}^2 & \cdots & a_{2n}^2 \\ \vdots & & & \\ a_{m1}^2 & a_{m2}^2 & \cdots & a_{mn}^2 \end{bmatrix}$$

*'$\cup$' is just the notational convenience (symbol) only.*

The above array is called a m by n bimatrix (written as $B(m \times n)$ since each of $A_i$ $(i = 1, 2)$ has m rows and n columns. It is to be noted a bimatrix has no numerical value associated with it. It is only a convenient way of representing a pair of arrays of numbers.

*Note:* If $A_1 = A_2$ then $A_B = A_1 \cup A_2$ is not a bimatrix. A bimatrix $A_B$ is denoted by $\left(a_{ij}^1\right) \cup \left(a_{ij}^2\right)$. If both $A_1$ and $A_2$ are $m \times n$ matrices then the bimatrix $A_B$ is called the $m \times n$ rectangular bimatrix. But we make an assumption the zero bimatrix is a union of two zero matrices even if $A_1$ and $A_2$ are one and the same; i.e., $A_1 = A_2 = (0)$.

***Example 1.1.1:*** The following are bimatrices

i. $A_B = \begin{bmatrix} 3 & 0 & 1 \\ -1 & 2 & 1 \end{bmatrix} \cup \begin{bmatrix} 0 & 2 & -1 \\ 1 & 1 & 0 \end{bmatrix}$

is a rectangular $2 \times 3$ bimatrix.

ii. $A'_B = \begin{bmatrix} 3 \\ 1 \\ 2 \end{bmatrix} \cup \begin{bmatrix} 0 \\ -1 \\ 0 \end{bmatrix}$

is a column bimatrix.

iii. $A''_B = (3, -2, 0, 1, 1) \cup (1, 1, -1, 1, 2)$

is a row bimatrix.

In a bimatrix $A_B = A_1 \cup A_2$ if both $A_1$ and $A_2$ are $m \times n$ rectangular matrices then the bimatrix $A_B$ is called the rectangular $m \times n$ matrix.

**DEFINITION 1.1.2:** *Let $A_B = A_1 \cup A_2$ be a bimatrix. If both $A_1$ and $A_2$ are square matrices then $A_B$ is called the square bimatrix.*

*If one of the matrices in the bimatrix $A_B = A_1 \cup A_2$ is square and other is rectangular or if both $A_1$ and $A_2$ are rectangular matrices say $m_1 \times n_1$ and $m_2 \times n_2$ with $m_1 \neq m_2$ or $n_1 \neq n_2$ then we say $A_B$ is a mixed bimatrix.*

The following are examples of a square bimatrix and the mixed bimatrix.

***Example 1.1.2:*** Given

$$A_B = \begin{bmatrix} 3 & 0 & 1 \\ 2 & 1 & 1 \\ -1 & 1 & 0 \end{bmatrix} \cup \begin{bmatrix} 4 & 1 & 1 \\ 2 & 1 & 0 \\ 0 & 0 & 1 \end{bmatrix}$$

is a 3 × 3 square bimatrix.

$$A'_B = \begin{bmatrix} 1 & 1 & 0 & 0 \\ 2 & 0 & 0 & 1 \\ 0 & 0 & 0 & 3 \\ 1 & 0 & 1 & 2 \end{bmatrix} \cup \begin{bmatrix} 2 & 0 & 0 & -1 \\ -1 & 0 & 1 & 0 \\ 0 & -1 & 0 & 3 \\ -3 & -2 & 0 & 0 \end{bmatrix}$$

is a 4 × 4 square bimatrix.

***Example 1.1.3:*** Let

$$A_B = \begin{bmatrix} 3 & 0 & 1 & 2 \\ 0 & 0 & 1 & 1 \\ 2 & 1 & 0 & 0 \\ 1 & 0 & 1 & 0 \end{bmatrix} \cup \begin{bmatrix} 1 & 1 & 2 \\ 0 & 2 & 1 \\ 0 & 0 & 4 \end{bmatrix}$$

then $A_B$ is a mixed square bimatrix.
Let

$$A'_B = \begin{bmatrix} 2 & 0 & 1 & 1 \\ 0 & 1 & 0 & 1 \\ -1 & 0 & 2 & 1 \end{bmatrix} \cup \begin{bmatrix} 2 & 0 \\ 4 & -3 \end{bmatrix}$$

$A'_B$ is a mixed bimatrix.

Now we proceed on to give the bimatrix operations.

Let $A_B = A_1 \cup A_2$ and $C_B = C_1 \cup C_2$ be two bimatrices we say $A_B$ and $C_B$ are equal written as $A_B = C_B$ if and only if $A_1$ and $C_1$ are identical and $A_2$ and $C_2$ are identical i.e., $A_1 = C_1$ and $A_2 = C_2$.

If $A_B = A_1 \cup A_2$ and $C_B = C_1 \cup C_2$, we say $A_B$ is not equal to $C_B$ we write $A_B \neq C_B$ if and only if $A_1 \neq C_1$ or $A_2 \neq C_2$.

***Example 1.1.4:*** Let

$$A_B = \begin{bmatrix} 3 & 2 & 0 \\ 2 & 1 & 1 \end{bmatrix} \cup \begin{bmatrix} 0 & -1 & 2 \\ 0 & 0 & 1 \end{bmatrix}$$

and

$$C_B = \begin{bmatrix} 1 & 1 & 1 \\ 0 & 0 & 0 \end{bmatrix} \cup \begin{bmatrix} 2 & 0 & 1 \\ 1 & 0 & 2 \end{bmatrix}$$

clearly $A_B \neq C_B$. Let

$$A_B = \begin{bmatrix} 0 & 0 & 1 \\ 1 & 1 & 2 \end{bmatrix} \cup \begin{bmatrix} 0 & 4 & -2 \\ -3 & 0 & 0 \end{bmatrix}$$

$$C_B = \begin{bmatrix} 0 & 0 & 1 \\ 1 & 1 & 2 \end{bmatrix} \cup \begin{bmatrix} 0 & 0 & 0 \\ 1 & 0 & 1 \end{bmatrix}$$

clearly $A_B \neq C_B$.

If $A_B = C_B$ then we have $C_B = A_B$.

We now proceed on to define multiplication by a scalar. Given a bimatrix $A_B = A_1 \cup B_1$ and a scalar $\lambda$, the product of $\lambda$ and $A_B$ written $\lambda\, A_B$ is defined to be

$$\lambda A_B = \begin{bmatrix} \lambda a_{11} & \cdots & \lambda a_{1n} \\ \vdots & & \vdots \\ \lambda a_{m1} & \cdots & \lambda a_{mn} \end{bmatrix} \cup \begin{bmatrix} \lambda b_{11} & \cdots & \lambda b_{1n} \\ \vdots & & \vdots \\ \lambda b_{m1} & \cdots & \lambda b_{mn} \end{bmatrix}$$

each element of $A_1$ and $B_1$ are multiplied by $\lambda$. The product $\lambda\, A_B$ is then another bimatrix having m rows and n columns if $A_B$ has m rows and n columns.

We write

$$\begin{aligned} \lambda\, A_B &= [\lambda a_{ij}] \cup [\lambda b_{ij}] \\ &= [a_{ij}\lambda] \cup [b_{ij}\lambda] \\ &= A_B\, \lambda. \end{aligned}$$

***Example 1.1.5:*** Let

$$A_B = \begin{bmatrix} 2 & 0 & 1 \\ 3 & 3 & -1 \end{bmatrix} \cup \begin{bmatrix} 0 & 1 & -1 \\ 2 & 1 & 0 \end{bmatrix}$$

and $\lambda = 3$ then

$$3A_B = \begin{bmatrix} 6 & 0 & 3 \\ 9 & 9 & -3 \end{bmatrix} \cup \begin{bmatrix} 0 & 3 & -3 \\ 6 & 3 & 0 \end{bmatrix}.$$

If $\lambda = -2$ for

$$A_B = [3\ 1\ 2\ -4] \cup [0\ 1\ -1\ 0],$$
$$\lambda A_B = [-6\ -2\ -4\ 8] \cup [0\ -2\ 2\ 0].$$

Let $A_B = A_1 \cup B_1$ and $C_B = A_2 \cup B_2$ be any two $m \times n$ bimatrices. The sum $D_B$ of the bimatrices $A_B$ and $C_B$ is defined as $D_B = A_B + C_B = [A_1 \cup B_1] + [A_2 \cup B_2] = (A_1 + A_2) \cup [B_2 + B_2]$; where $A_1 + A_2$ and $B_1 + B_1$ are the usual addition of matrices i.e., if

$$A_B = \left(a_{ij}^1\right) \cup \left(b_{ij}^1\right)$$

and

$$C_B = \left(a_{ij}^2\right) \cup \left(b_{ij}^2\right)$$

then

$$A_B + C_B = D_B = \left(a_{ij}^1 + b_{ij}^2\right) \cup \left(b_{ij}^1 + b_{ij}^2\right) \ \left(\forall\, ij\right).$$

If we write in detail

$$A_B = \begin{bmatrix} a_{11}^1 & \cdots & a_{1n}^1 \\ \vdots & & \\ a_{m1}^1 & \cdots & a_{mn}^1 \end{bmatrix} \cup \begin{bmatrix} b_{11}^1 & \cdots & b_{1n}^1 \\ \vdots & & \\ b_{m1}^1 & \cdots & b_{mn}^1 \end{bmatrix}$$

$$C_B = \begin{bmatrix} a_{11}^2 & \cdots & a_{1n}^2 \\ \vdots & & \\ a_{m1}^2 & \cdots & a_{mn}^2 \end{bmatrix} \cup \begin{bmatrix} b_{11}^2 & \cdots & b_{1n}^2 \\ \vdots & & \\ b_{m1}^2 & \cdots & b_{mn}^2 \end{bmatrix}$$

$A_B + C_B =$

$$\begin{bmatrix} a^1_{11}+a^2_{11} & \dots & a^1_{1n}a^2_{1n} \\ \vdots & & \vdots \\ a^1_{m1}+a^2_{m1} & \dots & a^1_{mn}+a^2_{mn} \end{bmatrix} \cup \begin{bmatrix} b^1_{11}+b^2_{11} & \dots & b^1_{1n}+b^2_{1n} \\ \vdots & & \vdots \\ b^1_{m1}+b^2_{m1} & \dots & b^1_{mn}+b^2_{mn} \end{bmatrix}.$$

The expression is abbreviated to

$$\begin{aligned} D_B &= A_B + C_B \\ &= (A_1 \cup B_1) + (A_2 \cup B_2) \\ &= (A_1 + A_2) \cup (B_1 + B_2). \end{aligned}$$

Thus two bimatrices are added by adding the corresponding elements only when compatibility of usual matrix addition exists.

*Note*: If $A_B = A^1 \cup A^2$ be a bimatrix we call $A^1$ and $A^2$ as the components matrices of the bimatrix $A_B$.

***Example 1.1.6:***

(i) Let

$$A_B = \begin{bmatrix} 3 & 1 & 1 \\ -1 & 0 & 2 \end{bmatrix} \cup \begin{bmatrix} 4 & 0 & -1 \\ 0 & 1 & -2 \end{bmatrix}$$

and

$$C_B = \begin{bmatrix} -1 & 0 & 1 \\ 2 & 2 & -1 \end{bmatrix} \cup \begin{bmatrix} 3 & 3 & 1 \\ 0 & 2 & -1 \end{bmatrix},$$

then, $D_B = A_B + C_B$

$$= \begin{bmatrix} 3 & 1 & 1 \\ -1 & 0 & 2 \end{bmatrix} + \begin{bmatrix} -1 & 0 & 1 \\ 2 & 2 & -1 \end{bmatrix} \cup \begin{bmatrix} 4 & 0 & -1 \\ 0 & 1 & -2 \end{bmatrix} + \begin{bmatrix} 3 & 3 & 1 \\ 0 & 2 & -1 \end{bmatrix}$$

$$= \begin{bmatrix} 2 & 1 & 2 \\ 1 & 2 & 1 \end{bmatrix} \cup \begin{bmatrix} 7 & 3 & 0 \\ 0 & 3 & -3 \end{bmatrix}.$$

(ii) Let $A_B = (3\ 2\ –1\ 0\ 1) \cup (0\ 1\ 1\ 0\ –1)$ and
$C_B = (1\ 1\ 1\ 1\ 1) \cup (5\ –1\ 2\ 0\ 3)$,
$A_B + C_B = (4\ 3\ 0\ 1\ 2) \cup (5\ 0\ 3\ 0\ 2)$.

***Example 1.1.7:*** Let

$$A_B = \begin{bmatrix} 6 & -1 \\ 2 & 2 \\ 1 & -1 \end{bmatrix} \cup \begin{bmatrix} 3 & 1 \\ 0 & 2 \\ -1 & 3 \end{bmatrix}$$

and

$$C_B = \begin{bmatrix} 2 & -4 \\ 4 & -1 \\ 3 & 0 \end{bmatrix} \cup \begin{bmatrix} 1 & 4 \\ 2 & 1 \\ 3 & 1 \end{bmatrix}.$$

$$A_B + A_B = \begin{bmatrix} 12 & -2 \\ 4 & 4 \\ 2 & -2 \end{bmatrix} \cup \begin{bmatrix} 6 & 2 \\ 0 & 4 \\ -2 & 6 \end{bmatrix} = 2A_B$$

$$C_B + C_B = \begin{bmatrix} 4 & -8 \\ 8 & -2 \\ 6 & 0 \end{bmatrix} \cup \begin{bmatrix} 2 & 8 \\ 4 & 2 \\ 6 & 2 \end{bmatrix} = 2C_B.$$

Similarly we can add

$$A_B + A_B + A_B = 3A_B = \begin{bmatrix} 18 & -3 \\ 6 & 6 \\ 3 & -3 \end{bmatrix} \cup \begin{bmatrix} 9 & 3 \\ 0 & 6 \\ -3 & 9 \end{bmatrix}.$$

*Note:* Addition of bimatrices are defined if and only if both the bimatrices are $m \times n$ bimatrices.
Let

$$A_B = \begin{bmatrix} 3 & 0 & 1 \\ 1 & 2 & 0 \end{bmatrix} \cup \begin{bmatrix} 1 & 1 & 1 \\ 0 & 2 & -1 \end{bmatrix}$$

$$C_B = \begin{bmatrix} 3 & 1 \\ 2 & 1 \\ 0 & 0 \end{bmatrix} \cup \begin{bmatrix} 1 & 1 \\ 2 & -1 \\ 3 & 0 \end{bmatrix}.$$

The addition of $A_B$ with $C_B$ is not defined for $A_B$ is a $2 \times 3$ bimatrix where as $C_B$ is a $3 \times 2$ bimatrix.

Clearly $A_B + C_B = C_B + A_B$ when both $A_B$ and $C_B$ are $m \times n$ matrices.

Also if $A_B$, $C_B$, $D_B$ be any three $m \times n$ bimatrices then $A_B + (C_B + D_B) = (A_B + C_B) + D_B$.

Subtraction is defined in terms of operations already considered for if

$$A_B = A_1 \cup A_2$$

and

$$B_B = B_1 \cup B_2$$

then

$$\begin{aligned} A_B - B_B &= A_B + (-B_B) \\ &= (A_1 \cup A_2) + (-B_1 \cup -B_2) \\ &= (A_1 - B_1) \cup (A_2 - B_2) \\ &= [A_1 + (-B_1)] \cup [A_2 + (-B_2)]. \end{aligned}$$

***Example 1.1.8:***

i. Let

$$A_B = \begin{bmatrix} 3 & 1 \\ -1 & 2 \\ 0 & 3 \end{bmatrix} \cup \begin{bmatrix} 5 & -2 \\ 1 & 1 \\ 3 & -2 \end{bmatrix}$$

and

$$B_B = \begin{bmatrix} 8 & -1 \\ 4 & 2 \\ -1 & 3 \end{bmatrix} \cup \begin{bmatrix} 9 & 2 \\ 2 & 9 \\ -1 & 1 \end{bmatrix}$$

$A_B - B_B = A_B + (-B_B)$.

$$= \left\{\begin{bmatrix} 3 & 1 \\ -1 & 2 \\ 0 & 3 \end{bmatrix} \cup \begin{bmatrix} 5 & -2 \\ 1 & 1 \\ 3 & -2 \end{bmatrix}\right\} + - \left\{\begin{bmatrix} 8 & -1 \\ 4 & 2 \\ -1 & 3 \end{bmatrix} \cup \begin{bmatrix} 9 & 2 \\ 2 & 9 \\ -1 & 1 \end{bmatrix}\right\}$$

$$= \left\{\begin{bmatrix} 3 & 1 \\ -1 & 2 \\ 0 & 3 \end{bmatrix} \cup \begin{bmatrix} 8 & -1 \\ 4 & 2 \\ -1 & 3 \end{bmatrix}\right\} + \left( - \left\{\begin{bmatrix} 5 & -2 \\ 1 & 1 \\ 3 & -2 \end{bmatrix} \cup \begin{bmatrix} 9 & 2 \\ 2 & 9 \\ -1 & 1 \end{bmatrix}\right\}\right)$$

$$= \begin{bmatrix} -5 & 2 \\ -5 & 0 \\ 1 & 0 \end{bmatrix} \cup \begin{bmatrix} 4 & -4 \\ -1 & -8 \\ 4 & -3 \end{bmatrix}.$$

ii. Let

$$A_B = (1, 2, 3, -1, 2, 1) \cup (3, -1, 2, 0, 3, 1)$$

and

$$B_B = (-1, 1, 1, 1, 1, 0) \cup (2, 0, -2, 0, 3, 0)$$

then

$$A_B + (-B_B) = (2, 1, 2, -2, 1, 1) \cup (1, -1, 4, 0, 0, 1).$$

Now we have defined addition and subtraction of bimatrices. Unlike in matrices we cannot say if we add two bimatrices the sum will be a bimatrix.

## 1.2 Some Basic properties of bimatrices

First we wish to mention that when we have a collection of m × n bimatrices say $M_B$ then $M_B$ need not be even closed with respect to addition. Further we make a definition, a m × m zero bimatrix to be $0_B = (0)_{m\times m} \cup (0)_{m\times m}$. Similarly the unit square bimatrix denoted by $I_B = I^1_{m\times m} \cup I^2_{m\times m}$. Thus we make the following special type of concession in case of zero and unit m × m bimatrices. Appropriate changes are made in case of zero m × n bimatrix.

We first illustrate this by the following example; that in general sum of two bimatrices is not a bimatrix.

***Example 1.2.1:*** Let

$$A_B = \begin{bmatrix} 1 & 1 & 1 \\ -1 & 0 & 3 \end{bmatrix} \cup \begin{bmatrix} 1 & 2 & 0 \\ 2 & 1 & 2 \end{bmatrix}$$

and

$$B_B = \begin{bmatrix} 2 & 1 & 0 \\ 1 & 1 & 1 \end{bmatrix} \cup \begin{bmatrix} 2 & 0 & 1 \\ -2 & 0 & 2 \end{bmatrix}.$$

$A_B + B_B$

$$= \left\{\begin{bmatrix} 1 & 1 & 1 \\ -1 & 0 & 3 \end{bmatrix} + \begin{bmatrix} 2 & 1 & 0 \\ 1 & 1 & 1 \end{bmatrix}\right\} \cup \left\{\begin{bmatrix} 1 & 2 & 0 \\ 2 & 1 & 2 \end{bmatrix} + \begin{bmatrix} 2 & 0 & 1 \\ -2 & 0 & 2 \end{bmatrix}\right\}$$

$$= \begin{bmatrix} 3 & 2 & 1 \\ 0 & 1 & 4 \end{bmatrix} \cup \begin{bmatrix} 3 & 2 & 1 \\ 0 & 1 & 4 \end{bmatrix}.$$

Clearly $A_B + B_B$ is not a bimatrix as $A_1 \cup B_1 = A_2 \cup B_2$ where $A_B = A_1 \cup A_2$ and $B_B = B_1 \cup B_2$.

Thus it is highly erratic, for we see the sum of two bimatrices is not in general a bimatrix.

Now in the following theorem we give a necessary and sufficient condition for the sum of two bimatrices to be a bimatrix.

**THEOREM 1.2.1:** *Let $A_B = A_1 \cup A_2$ and $B_B = B_1 \cup B_2$ be two $m \times n$ bimatrices.*

$$A_B + B_B = (A_1 + B_1) \cup (A_2 + B_2)$$

*is a bimatrix if and only if $A_1 + B_1 \neq A_2 + B_2$.*

*Proof:* Let $A_B = A_1 \cup A_2$ and $B_B = B_1 \cup B_2$ be two $m \times n$ bimatrices. If $A_B + B_B = (A_1 + B_1) \cup (A_2 + B_2)$ is a bimatrix then $A_1 + B_1 \neq A_2 + B_2$.

On the other hand in $A_B + B_B$ if $A_1 + B_1 \neq A_2 + B_2$ then clearly $A_B + B_B$ is a bimatrix. Hence the theorem.

**COROLLARY:** *If $A_B = A_1 \cup A_2$ and $B_B = B_1 \cup B_2$ be two m × n bimatrices then $A_B - B_B = A_B + (-B_B)$ is a bimatrix if and only if $A_1 + (-B_1) \neq A_2 + (-B_2)$.*

Now all the above results continue to be true when the bimatrix under consideration are square bimatrices.
In case of mixed bimatrices we have to make some more additional observations.

***Example 1.2.2:*** Let

$$A_B = \begin{bmatrix} 3 & 1 \\ 2 & 4 \end{bmatrix} \cup \begin{bmatrix} 1 \\ 0 \end{bmatrix}$$

and

$$B_B = \begin{bmatrix} 4 & -1 \\ 2 & 0 \end{bmatrix} \cup \begin{bmatrix} 7 \\ -2 \end{bmatrix}$$

be two mixed bimatrices. The sum can be obtained in case of mixed bimatrices only when the addition is compatible i.e., when the addition is defined between the component matrices of $A_B$ and $B_B$.

$$A_B + B_B = \begin{bmatrix} 7 & 0 \\ 4 & 4 \end{bmatrix} \cup \begin{bmatrix} 8 \\ -2 \end{bmatrix}.$$

Now when we say $A_B = A_1 \cup A_2$ is a square bimatrix we mean both $A_1$ and $A_2$ are m × m matrices, and from now on wards we denote it by $A_B^{m\times m}$ when we say $A_B = A_1 \cup A_2$ is a m × n bimatrix if and only if both $A_1$ and $A_2$ are both m × n matrices. We shall denote a m × n bimatrix $A_B$ by $A_B^{m\times n}$.

But when $A_B = A_1 \cup A_2$ is a mixed square bimatrix we may have $A_1$ to be a m × m square matrix and $A_2$ to be a t × t square matrix. t ≠ m. Then we denote $A_B = A_1^{m\times m} \cup A_2^{t\times t}$.

Suppose $A_1$ is a m × m square matrix and $A_2$ to be p × q rectangular matrix then $A_B = A_1^{m\times m} \cup A_2^{p\times q}$ (p ≠ m or q ≠ m) or $A_1$ to a be m × n rectangular matrix $A_2$ to be a s × t rectangular matrix (m ≠ s or n ≠ t).

$$A_B = A_1^{m\times n} \cup A_2^{s\times t}$$

In all cases we call $A_B$ to be a mixed bimatrix. So if $A_B$ and $C_B$ are mixed bimatrices with compatible addition of their components then their sum is always a bimatrix for $A_1 \cup B_1 \neq A_2 \cup B_2$ as both the components in $A_B + C_B$ are of different orders.

**THEOREM 1.2.2:** *Let $A_B$ and $C_B$ be mixed bimatrices with compatible addition then $A_B + C_B$ is always a mixed bimatrix.*

*Proof:* Follows from the very definition.

Now we proceed on to define bimatrix multiplication. First it is clear that bimatrix multiplication is not defined when the bimatrices are not square bimatrix. Secondly in case of mixed bimatrix, multiplication is defined when both $A_1$ and $A_2$ are square matrices or when $A_B = A_1^{m\times n} \cup A_2^{p\times s}$ and $C_B = B_1^{n\times u} \cup B_2^{s\times t}$.

We first illustrate these operations by examples.

***Example 1.2.3:*** If $A_B = A_1 \cup A_2$ and $B_B = B_1 \cup B_1$ are both m × n bimatrices m ≠ n, the bimatrix multiplication is not defined.

***Example 1.2.4:*** Let

$$A_B^{m\times m} = A_1 \cup A_2$$

and

$$B_B^{m\times m} = B_1 \cup B_2$$

then

$$A_B^{m\times m} \times B_B^{m\times m} = (A_1B_1) \cup (A_2B_2).$$

$A_B^{m\times m} \times B_B^{m\times m}$ is bimatrix if and only if $A_1 B_1 \neq A_2 B_2$ i.e., components of the bimatrix are compatible with matrix multiplication.

**THEOREM 1.2.3:** *Let $A_B = A_1 \cup A_2$ and $B_B = B_1 \cup B_1$ be two $m \times m$ square bimatrices $A_B$ . $B_B = A_1 B_1 \cup A_2 B_2$ is a bimatrix if and only if $A_1 B_1 \neq A_2 B_2$.*

*Proof:* Given $A_B = A_1 \cup A_2$ and $B_B = B_1 \cup B_2$ are two $m \times m$ square bimatrices and $A_B B_B$ is a bimatrix i.e.,

$$A_B B_B = (A_1 \cup A_2) (B_1 \cup B_2) = A_1 B_1 \cup A_2 B_2$$

is a bimatrix. Hence $A_1B_1 \neq A_2B_2$.

Clearly if $A_1 B_1 \neq A_2 B_2$ then $A_B B_B = A_1B_1 \cup A_2B_2$ is a $m \times m$ bimatrix.

***Example 1.2.5:*** Let

$$A_B = \begin{bmatrix} 3 & 0 \\ 1 & 2 \end{bmatrix} \cup \begin{bmatrix} 1 & 1 \\ 0 & 2 \end{bmatrix}$$

and

$$B_B = \begin{bmatrix} -2 & 0 \\ 1 & 0 \end{bmatrix} \cup \begin{bmatrix} 2 & -1 \\ 0 & 0 \end{bmatrix}$$

be two $2 \times 2$ square bimatrices.

$$A_BB_B = \begin{bmatrix} 3 & 0 \\ 1 & 2 \end{bmatrix}\begin{bmatrix} -2 & 0 \\ 1 & 0 \end{bmatrix} \cup \begin{bmatrix} 1 & 1 \\ 0 & 2 \end{bmatrix}\begin{bmatrix} 2 & -1 \\ 0 & 0 \end{bmatrix}$$

$$= \begin{bmatrix} -6 & 0 \\ 0 & 0 \end{bmatrix} \cup \begin{bmatrix} 2 & -1 \\ 0 & 0 \end{bmatrix}$$

is a bimatrix.

Let $A_B = A_1^{m\times m} \cup A_2^{t\times t}$ and $B_B = B_1^{m\times m} \cup B_2^{t\times t}$. Clearly $A_B B_B$ is a mixed square bimatrix.

***Example 1.2.6:*** Let

$$A_B = \begin{bmatrix} 3 & 0 & 1 \\ 1 & 0 & 1 \\ 0 & 1 & 0 \end{bmatrix} \cup \begin{bmatrix} 2 & 0 \\ 1 & 1 \end{bmatrix}$$

and

$$B_B = \begin{bmatrix} 1 & 0 & 1 \\ 0 & 0 & 1 \\ 1 & 0 & 0 \end{bmatrix} \cup \begin{bmatrix} 0 & 3 \\ 1 & 0 \end{bmatrix}$$

$$A_B\, B_B \;=\; \begin{bmatrix} 3 & 0 & 1 \\ 1 & 0 & 1 \\ 0 & 1 & 0 \end{bmatrix} \begin{bmatrix} 1 & 0 & 1 \\ 0 & 0 & 1 \\ 1 & 0 & 0 \end{bmatrix} \cup \begin{bmatrix} 2 & 0 \\ 1 & 1 \end{bmatrix} \begin{bmatrix} 0 & 3 \\ 1 & 0 \end{bmatrix}$$

$$= \begin{bmatrix} 4 & 0 & 3 \\ 2 & 0 & 1 \\ 0 & 0 & 1 \end{bmatrix} \cup \begin{bmatrix} 0 & 6 \\ 1 & 3 \end{bmatrix}$$

is a mixed square bimatrix.

***Example 1.2.7:*** Let

$$A_B = \begin{bmatrix} 3 & 2 \\ -1 & 4 \\ 0 & 3 \end{bmatrix} \cup \begin{bmatrix} 2 & 0 & 1 \\ 0 & 2 & 1 \end{bmatrix}$$

and

$$B_B = \begin{bmatrix} 0 & 1 & 1 \\ 2 & 0 & -1 \end{bmatrix} \cup \begin{bmatrix} 2 & 0 \\ 1 & 0 \\ 0 & -1 \end{bmatrix}.$$

The product $A_B\, B_B$ is defined as

$$A_B\, B_B \;=\; \begin{bmatrix} 3 & 2 \\ -1 & 4 \\ 0 & 3 \end{bmatrix} \begin{bmatrix} 0 & 1 & 1 \\ 2 & 0 & -1 \end{bmatrix} \cup \begin{bmatrix} 2 & 0 & 1 \\ 0 & 2 & 1 \end{bmatrix} \begin{bmatrix} 2 & 0 \\ 1 & 0 \\ 0 & -1 \end{bmatrix}$$

$$= \begin{bmatrix} 4 & 3 & 1 \\ 8 & -1 & -5 \\ 6 & 0 & -3 \end{bmatrix} \cup \begin{bmatrix} 4 & -1 \\ 2 & -1 \end{bmatrix}$$

is a square mixed matrix.

Thus we have the following theorem the proof of which is left for the reader as an exercise.

**THEOREM 1.2.4:** *Let* $A_B = A_1^{m\times n} \cup A_2^{p\times q}$ *be a mixed bimatrix and* $B_B = B_1^{n\times m} \cup B_2^{q\times p}$ *be another mixed bimatrix. Then the product is defined and is a mixed square bimatrix as*

$$A_B\, B_B = C_1^{m\times m} \cup C_2^{p\times p}$$

*with*

$$C_1 = C_1^{m\times m} \cup C_2^{p\times p}$$

*where*

$$C_1^{m\times m} = A_1^{m\times n} \times B_1^{n\times m}$$

*and*

$$C_2^{p\times p} = A_2^{p\times q} \times B_2^{q\times p}.$$

It is important to note that in general the product of two mixed bimatrices need not be a square mixed bimatrix.

The following example is evident to show that the product of two mixed bimatrices in general is not a square mixed bimatrix.

***Example 1.2.8:*** Let

$$A_B = \begin{bmatrix} 3 & 2 & 1 \\ 0 & 0 & 1 \\ 1 & 0 & 1 \\ 2 & 3 & -4 \end{bmatrix} \cup \begin{bmatrix} 3 & 2 & 1 & 3 \\ 0 & 1 & 0 & -1 \end{bmatrix}$$

and

$$B_B = \begin{bmatrix} 0 & 1 \\ 3 & 0 \\ -1 & 2 \end{bmatrix} \cup \begin{bmatrix} 3 & 0 & 0 \\ 0 & -1 & 0 \\ 1 & 0 & 0 \\ -1 & 2 & 1 \end{bmatrix}$$

be two mixed bimatrices. Clearly the product $A_B$ $B_B$ is defined and it is a mixed bimatrix but not a square mixed bimatrix.

It can also happen that the product of two mixed bimatrices can be a square bimatrix. This is illustrated by the following example.

***Example 1.2.9:*** Let

$$A_B = \begin{bmatrix} 3 & 2 \\ 1 & 0 \\ 3 & 1 \end{bmatrix} \cup \begin{bmatrix} 0 & 1 & 0 & 0 \\ 1 & 0 & 0 & 0 \\ 0 & 1 & 0 & 1 \end{bmatrix}$$

and

$$B_B = \begin{bmatrix} 0 & 1 & 2 \\ -1 & 0 & -1 \end{bmatrix} \cup \begin{bmatrix} 0 & 1 & 1 \\ -1 & 0 & 0 \\ 0 & 1 & 1 \\ 0 & 1 & 0 \end{bmatrix}.$$

Now both $A_B$ and $B_B$ are mixed bimatrices. Clearly then product $A_B$ $B_B$ is defined, for

$$A_BB_B = \begin{bmatrix} 3 & 2 \\ 1 & 0 \\ 3 & 1 \end{bmatrix} \begin{bmatrix} 0 & 1 & 2 \\ -1 & 0 & -1 \end{bmatrix} \cup \begin{bmatrix} 0 & 1 & 0 & 0 \\ 1 & 0 & 0 & 0 \\ 0 & 1 & 0 & 1 \end{bmatrix} \begin{bmatrix} 0 & 1 & 1 \\ -1 & 0 & 0 \\ 0 & 1 & 1 \\ 0 & 1 & 0 \end{bmatrix}$$

$$= \begin{bmatrix} -2 & 3 & 4 \\ 0 & 1 & 2 \\ -1 & 3 & 5 \end{bmatrix} \cup \begin{bmatrix} -1 & 0 & 0 \\ 0 & 1 & 1 \\ -1 & 1 & 0 \end{bmatrix}.$$

It is clear that $A_B$ $B_B$ is only a square bimatrix, it is a $3 \times 3$ square bimatrix. Thus the product of a mixed bimatrix can be just a square bimatrix i.e., need not in general be a mixed bimatrix.

It is important to make note of the following:

*1. If $A_B$ and $B_B$ be two $m \times m$ square bimatrices in general $A_B\, B_B \neq B_B\, A_B$.*

The proof is shown by the following example:

***Example 1.2.10:*** Let

$$A_B = \begin{bmatrix} 3 & 0 & 1 \\ 0 & 1 & -1 \\ 0 & 1 & 0 \end{bmatrix} \cup \begin{bmatrix} 2 & 0 & 0 \\ 1 & 0 & 1 \\ 0 & 2 & -1 \end{bmatrix}$$

and

$$B_B = \begin{bmatrix} 5 & 2 & -1 \\ 0 & 0 & -3 \\ 1 & 2 & 0 \end{bmatrix} \cup \begin{bmatrix} 1 & 0 & 1 \\ 0 & 3 & 1 \\ -1 & 0 & 2 \end{bmatrix}$$

be $3 \times 3$ bimatrices.

$A_B B_B$

$$= \begin{bmatrix} 3 & 0 & 1 \\ 0 & 1 & -1 \\ 0 & 1 & 0 \end{bmatrix} \begin{bmatrix} 5 & 2 & -1 \\ 0 & 0 & -3 \\ 1 & 2 & 0 \end{bmatrix} \cup \begin{bmatrix} 2 & 0 & 0 \\ 1 & 0 & 1 \\ 0 & 2 & -1 \end{bmatrix} \begin{bmatrix} 1 & 0 & 1 \\ 0 & 3 & 1 \\ -1 & 0 & 2 \end{bmatrix}$$

$$= \begin{bmatrix} 16 & 8 & -3 \\ -1 & -2 & -3 \\ 0 & 0 & -3 \end{bmatrix} \cup \begin{bmatrix} 2 & 0 & 2 \\ 0 & 0 & 3 \\ 1 & 6 & 0 \end{bmatrix}.$$

Consider $B_B A_B$

$$= \begin{bmatrix} 5 & 2 & -1 \\ 0 & 0 & -3 \\ 1 & 2 & 0 \end{bmatrix} \begin{bmatrix} 3 & 0 & 1 \\ 0 & 1 & -1 \\ 0 & 1 & 0 \end{bmatrix} \cup \begin{bmatrix} 1 & 0 & 1 \\ 0 & 2 & 1 \\ -1 & 0 & 2 \end{bmatrix} \begin{bmatrix} 2 & 0 & 0 \\ 1 & 0 & 1 \\ 0 & 2 & -1 \end{bmatrix}$$

$$= \begin{bmatrix} 15 & 1 & -2 \\ 0 & -3 & 0 \\ 3 & 2 & -1 \end{bmatrix} \cup \begin{bmatrix} 2 & 2 & -1 \\ 3 & 2 & 2 \\ -2 & 4 & -2 \end{bmatrix}.$$

Clearly $A_B\ B_B \neq B_B\ A_B$ but both of them are $3 \times 3$ square bimatrices.

*2. In some cases for the bimatrices $A_B$ and $B_B$ only one type of product $A_B\ B_B$ may be defined and $B_B\ A_B$ may not be even defined.*

This is shown by the following example:

***Example 1.2.11:*** Let

$$A_B = \begin{bmatrix} 0 & 1 & 1 & 1 \\ 1 & 0 & 0 & 2 \\ -1 & 0 & 1 & -1 \end{bmatrix} \cup \begin{bmatrix} 2 & 1 \\ 3 & 0 \\ 1 & 2 \end{bmatrix}$$

and

$$B_B = \begin{bmatrix} 1 & 2 \\ 0 & 1 \\ 3 & 1 \\ - & -1 \end{bmatrix} \cup \begin{bmatrix} 3 & 0 & 1 \\ 2 & -1 & 4 \end{bmatrix}.$$

Clearly $A_B\ B_B$

$$= \begin{bmatrix} 0 & 1 & 1 & 1 \\ 1 & 0 & 0 & 2 \\ -1 & 0 & 1 & -1 \end{bmatrix} \begin{bmatrix} 1 & 2 \\ 0 & 1 \\ 3 & 0 \\ 1 & -1 \end{bmatrix} \cup \begin{bmatrix} 2 & 1 \\ 3 & 0 \\ 1 & 2 \end{bmatrix} \begin{bmatrix} 3 & 0 & 1 \\ 2 & -1 & 4 \end{bmatrix}$$

$$= \begin{bmatrix} 4 & 0 \\ 3 & 0 \\ -2 & -1 \end{bmatrix} \cup \begin{bmatrix} 8 & -1 & 6 \\ 9 & 0 & 3 \\ 7 & -2 & 9 \end{bmatrix}.$$

is a mixed bimatrix. But $B_B A_B$ is not even defined.

It is left for the reader to verify the following identities in case of square bimatrices with examples.

1. $A_B (B_B C_B) = (A_B B_B) C_B = A_B B_B C_B$ i.e., if $A_B = A_1 \cup A_2$, $B_B = B_1 \cup B_2$ and $C_B = C_1 \cup C_2$.
   Clearly
   $$A_B (B_B C_B) = A_B [B_1 C_1 \cup B_2 C_2] = A_1 (B_1 C_1) \cup A_2 (B_2 C_2)$$
   as matrix multiplication is associative we have
   $A_B (B_B C_B) = (A_B B_B) C_B = A_B B_B C_B$.

2. Further bimatrix multiplication satisfies the distributive law for if $A_B = A_1 \cup A_2$, $B_B = B_1 \cup B_2$ and $C_B = C_1 \cup C_2$; then

$$\begin{aligned} A_B (B_B + C_B) &= A_B [(B_1 + C_1) \cup (B_2 \cup C_2) \\ &= (A_1 (B_1 + C_1) \cup A_2 (B_2 + C_2) \\ &= (A_1B_1 + A_1 C_1) \cup (A_2 B_2 + A_2 C_2) \\ &= (A_1 B_1 \cup A_2 B_2) + A_1 C_1 U A_2 C_2 \\ &= A_B B_B + A_B C_B. \end{aligned}$$

We use the fact that matrix multiplication is distributive and we obtain the bimatrix multiplication is also distributive.

**DEFINITION 1.2.1:** *Let* $A_B^{m\times m} = A_1 \cup A_2$ *be a* $m \times m$ *square bimatrix. We define* $I_B^{m\times m} = I^{m\times m} \cup I^{m\times m} = I_1^{m\times m} \cup I_2^{m\times m}$ *to be the identity bimatrix.*

Suppose $A_B = A_1^{m\times m} \cup A_2^{n\times n}$ be a mixed square bimatrix then the identity bimatrix $I_B$ is defined to be the bimatrix $I^{m\times m} \cup I^{n\times n} = I_1^{m\times m} \cup I_2^{n\times n}$.

**DEFINITION 1.2.2:** *Let* $A_B = A_1 \cup B_1$ *be a* $m \times m$ *square bimatrix. We say* $A_B$ *is a diagonal bimatrix if each of* $A_1$ *and* $B_1$ *are diagonal* $m \times m$ *matrices. Clearly the identity bimatrix is a diagonal bimatrix. If* $A_B = A_1^{m\times m} \cup A_2^{n\times n}$ *be a* *mixed square matrix we say $A_B$ is a diagonal mixed bimatrix if both $A_1$ and $A_2$ are diagonal matrices.*

It is clear that as in case of matrices diagonal bimatrix cannot be defined in case of rectangular bimatrix or mixed bimatrix which is not a square mixed bimatrix.

It is easily verified if $I^{m\times m} = I_1^{m\times m} \cup I_2^{m\times m}$ is the identity square bimatrix then for any m × m square bimatrix $A_B^{m\times m} = A_1 \cup A_2$ we have $A_B^{m\times m}\ I^{m\times m} = I^{m\times m} A_B^{m\times m} = A_B^{m\times m}$.

Further it can be verified that if $I^{m\times m} = I_1^{m\times m} \cup I_2^{m\times m}$ be the identity square bimatrix then $\left(I^{m\times m}\right)^2 = I^{m\times m}$.

A bimatrix $A_B$ all of whose elements are zero is called a null or zero bimatrix and it is denoted by $O_B = O_1 \cup O_2 = (0)_1 \cup (0)_2$.

For

$$O_B = (0\ 0\ 0\ 0) \cup \begin{bmatrix} 0 & 0 & 0 \\ 0 & 0 & 0 \end{bmatrix} \qquad (1)$$

$$O_B = \begin{bmatrix} 0 & 0 & 0 & 0 \\ 0 & 0 & 0 & 0 \\ 0 & 0 & 0 & 0 \end{bmatrix} \cup \begin{bmatrix} 0 & 0 & 0 & 0 \\ 0 & 0 & 0 & 0 \\ 0 & 0 & 0 & 0 \end{bmatrix} \qquad (2)$$

$$O_B = \begin{bmatrix} 0 & 0 & 0 \\ 0 & 0 & 0 \\ 0 & 0 & 0 \end{bmatrix} \cup \begin{bmatrix} 0 & 0 & 0 \\ 0 & 0 & 0 \\ 0 & 0 & 0 \end{bmatrix} \qquad (3)$$

$$O_B = \begin{bmatrix} 0 & 0 \\ 0 & 0 \end{bmatrix} \cup \begin{bmatrix} 0 & 0 & 0 & 0 \\ 0 & 0 & 0 & 0 \\ 0 & 0 & 0 & 0 \\ 0 & 0 & 0 & 0 \end{bmatrix} \qquad (4)$$

$$O_B = \begin{bmatrix} 0 & 0 & 0 & 0 \\ 0 & 0 & 0 & 0 \\ 0 & 0 & 0 & 0 \end{bmatrix} \cup \begin{bmatrix} 0 & 0 & 0 \\ 0 & 0 & 0 \\ 0 & 0 & 0 \end{bmatrix}. \qquad (5)$$

The $O_B$ given by (1) is a mixed zero bimatrix. The $O_B$ given by (2) is a rectangular zero bimatrix. The $O_B$ given by (3) is a square zero bimatrix. The $O_B$ given by (4) is a mixed square zero bimatrix. The $O_B$ of identity (5) is a mixed zero bimatrix.

Thus we see in case of zero bimatrix it can be rectangular bimatrix or a square bimatrix or a mixed square bimatrix or mixed bimatrix.

For every bimatrix $A_B$ there exist a zero bimatrix $O_B$ such that $A_B + O_B = O_B + A_B = A_B$.

But we cannot say $O_B \times A_B = O_B$, this is true only in case of mixed square bimatrix and square bimatrix only. Now if $A_B$ is a square bimatrix of a mixed square bimatrix we see $A_B\, A_B = A^2_B$ , $A_B.\, A_B.\, A_B = A^3_B$ and so on.

This type of product does not exist in case of rectangular bimatrix or mixed rectangular bimatrix.

Also for any scalar $\lambda$ the square bimatrix $A_B^{m\times m}$ is called a scalar bimatrix if $A_B^{m\times m} = \lambda I_B^{m\times m}$.

We also define the notion of mixed scalar bimatrix provided the bimatrix is a mixed square bimatrix. So if $A_B = A_1^{m\times m} \cup A_2^{n\times n}$ and $I_B = I_1^{m\times m} \cup I_2^{n\times n}$ then if

$$A_B = \lambda I_B = \lambda I_1^{m\times m} \cup \lambda I_2^{n\times n}.$$

Now null the bimatrix can be got for any form of matrices $A_B$ and $C_B$ provided the product $A_B\, C_B$ is defined and $A_B\, C_B = (0)$.

***Example 1.2.12:*** Let

$$A_B = \begin{bmatrix} 0 & 0 \\ 1 & 0 \end{bmatrix} \cup (0\,0\,0\,1\,0)$$

and

$$C_B = \begin{bmatrix} 0 & 0 & 0 \\ 4 & 1 & 3 \end{bmatrix} \cup \begin{bmatrix} 5 \\ 0 \\ 2 \\ 0 \\ 6 \end{bmatrix}.$$

Clearly $A_BC_B$ is a null bi matrix as $A_B\, C_B =$

$$= \begin{bmatrix} 0 & 0 \\ 1 & 0 \end{bmatrix} \begin{bmatrix} 0 & 0 & 0 \\ 4 & 1 & 3 \end{bmatrix} \cup (0\,0\,0\,1\,0) \begin{bmatrix} 5 \\ 0 \\ 2 \\ 0 \\ 6 \end{bmatrix}$$

$$= \begin{bmatrix} 0 & 0 & 0 \\ 0 & 0 & 0 \end{bmatrix} \cup (0).$$

Now we proceed on to define the transpose of a bimatrix $A_B$.

Let

$$A_B = A_1 \cup A_2 = \left(a^1_{ij}\right) \cup \left(a^2_{ij}\right)$$

be the bimatrix formed from $A_1$ and $A_2$. By interchanging rows and columns such that row i of $A_1$ and $A_2$ become column i of the transposed matrix $A_1$ and $A_2$ respectively.

The transpose is denoted by

$$A'_B = \left(a^1_{ji}\right) \cup \left(a^2_{ji}\right)$$

where

$$A_1 = \left(a^1_{ij}\right) \text{ and } A_2 = \left(a^2_{ij}\right).$$

It can be easily verified that if $A_B$ and $B_B$ are two bimatrices and if $C_B = A_B + B_B$ then $C'_B = A'_B + B'_B$.

It is left as a simple problem for the reader to prove that if $A_B\ B_B$ is defined then $(A_B\ B_B)' = B'_B\ A'_B$ that is the transpose of the product is the product of the transpose in the reverse order.

We illustrate this by the following example:

***Example 1.2.13:*** Let

$$A_B = \begin{bmatrix} 3 & 2 \\ 0 & 1 \\ 1 & 2 \end{bmatrix} \cup \begin{bmatrix} 2 & 3 & 1 & 4 \\ 0 & 1 & 0 & 2 \end{bmatrix}$$

and

$$B_B = \begin{bmatrix} 3 & 0 & 1 & 2 \\ 0 & 1 & 1 & -2 \end{bmatrix} \cup \begin{bmatrix} 3 & 3 \\ 0 & 0 \\ 1 & -1 \\ -2 & 0 \end{bmatrix}.$$

$A_B B_B$

$$= \begin{bmatrix} 3 & 2 \\ 0 & 1 \\ 1 & 2 \end{bmatrix} \begin{bmatrix} 3 & 0 & 1 & 2 \\ 0 & 1 & 1 & -2 \end{bmatrix} \cup \begin{bmatrix} 2 & 3 & 1 & 4 \\ 0 & 1 & 0 & 2 \end{bmatrix} \begin{bmatrix} 3 & 3 \\ 0 & 0 \\ 1 & -1 \\ 2 & 0 \end{bmatrix}$$

$$= \begin{bmatrix} 9 & 2 & 5 & 2 \\ 0 & 1 & 1 & -2 \\ 3 & 2 & 3 & -2 \end{bmatrix} \cup \begin{bmatrix} 15 & 5 \\ 4 & 0 \end{bmatrix}$$

$$(A_B\, B_B)' = \begin{bmatrix} 9 & 0 & 3 \\ 2 & 1 & 2 \\ 5 & 1 & 3 \\ 2 & -2 & -2 \end{bmatrix} \cup \begin{bmatrix} 15 & 4 \\ 5 & 0 \end{bmatrix}.$$

Consider

$$A'_B = \begin{bmatrix} 3 & 0 & 1 \\ 2 & 1 & 2 \end{bmatrix} \cup \begin{bmatrix} 2 & 0 \\ 3 & 1 \\ 1 & 0 \\ 4 & 2 \end{bmatrix}$$

$$B'_B = \begin{bmatrix} 3 & 0 \\ 0 & 1 \\ 1 & 1 \\ 2 & -2 \end{bmatrix} \cup \begin{bmatrix} 3 & 0 & 1 & -2 \\ 3 & 0 & -1 & 0 \end{bmatrix}$$

$$B'_BA'_B = \begin{bmatrix} 3 & 0 \\ 0 & 1 \\ 1 & 1 \\ 2 & -2 \end{bmatrix} \begin{bmatrix} 3 & 0 & 1 \\ 2 & 1 & 2 \end{bmatrix} \cup \begin{bmatrix} 3 & 0 & 1 & -2 \\ 3 & 0 & -1 & 0 \end{bmatrix} \begin{bmatrix} 2 & 0 \\ 3 & 1 \\ 1 & 0 \\ 4 & 2 \end{bmatrix}$$

$$= \begin{bmatrix} 9 & 0 & 3 \\ 2 & 1 & 2 \\ 5 & 1 & 3 \\ 2 & -2 & -2 \end{bmatrix} \cup \begin{bmatrix} -1 & -4 \\ 5 & 0 \end{bmatrix}.$$

Thus $(A_B\ B_B)' = B'_B\ A'_B$.

Now we leave it for the reader to prove that if $A_B$, $B_B$, $C_B$, …, $N_B$ be bimatrices such that their product $A_B\ B_B\ C_B$ … $N_B$ is well defined then we have $(A_B\ B_B\ C_B \ldots N_B)' = N'_B \ldots C'_B\ B'_B\ A'_B$.

## 1.3 Symmetric and skew symmetric bimatrices

In this section we just introduce the notion of symmetric and skew symmetric bimatrices and illustrate them with examples.

A symmetric bimatrix is a matrix $A_B$ for which $A_B = A'_B$ i.e., the component matrices of $A_B$ are also symmetric matrices.

***Example 1.3.1:*** Let

$$A_B = \begin{bmatrix} 3 & 0 & 2 \\ 0 & 1 & -1 \\ 2 & -1 & -5 \end{bmatrix} \cup \begin{bmatrix} 0 & 1 & 2 \\ 1 & -5 & 3 \\ 2 & 3 & 0 \end{bmatrix}$$

$$A'_B = \begin{bmatrix} 3 & 0 & 2 \\ 0 & 1 & -1 \\ 2 & -1 & -5 \end{bmatrix} \cup \begin{bmatrix} 0 & 1 & 2 \\ 1 & -5 & 3 \\ 2 & 3 & 0 \end{bmatrix}.$$

Thus $A'_B = A_B$ hence the bimatrix $A_B$ is a symmetric bimatrix.

***Example 1.3.2:*** Let

$$A_B = \begin{bmatrix} 2 & 0 \\ 0 & 1 \end{bmatrix} \cup \begin{bmatrix} 3 & 1 & 2 & 4 \\ 1 & 0 & -1 & 2 \\ 2 & -1 & 1 & -4 \\ 4 & 2 & -4 & 8 \end{bmatrix}$$

$$A'_B = \begin{bmatrix} 2 & 0 \\ 0 & 1 \end{bmatrix} \cup \begin{bmatrix} 3 & 1 & 2 & 4 \\ 1 & 0 & -1 & 2 \\ 2 & -1 & 1 & -4 \\ 4 & 2 & -4 & 8 \end{bmatrix}.$$

Thus $A_B = A'_B$. This mixed square bimatrix is symmetric.

Clearly a symmetric bimatrix must be either a square bimatrix or a mixed square bimatrix. Each of the component matrices in $A_B$ is symmetric about the main diagonal that is a reflection in the main diagonal leaves the component matrices of the bimatrix unchanged. (Recall $A_B = A_1 \cup A_2$ be a bimatrix $A_1$ and $A_2$ are called as the component matrices of the bimatrix $A_B$).

Let $A_B = A_1 \cup A_2$ be any $m \times m$ square bimatrix. This is a $m^{th}$ order square bimatrix.

This will not have $2m^2$ arbitrary elements since $a^1_{ij} = a^1_{ji}$ and $a^2_{ij} = a^2_{ji}$ (where $A_1 = (a^1_{ij})$ and $A_2 = (a^2_{ij})$) both below and above the main diagonal. The number of elements above the main diagonal of $A_B = A_1 \cup A_2$ is $(m^2 – m)$. The diagonal elements are also arbitrary. Thus the total number of arbitrary elements in an $m^{th}$ order symmetric bimatrix is $m^2 – m + 2m = m (m +1)$.

***Example 1.3.3:*** Let $A_B = A_1 \cup A_2$ be a $4 \times 4$ square symmetric bimatrix.

$$A_B = \begin{bmatrix} 3 & 0 & 1 & 1 \\ 0 & 2 & 0 & 1 \\ 1 & 0 & 5 & 0 \\ 1 & 1 & 0 & -3 \end{bmatrix} \cup \begin{bmatrix} -1 & 0 & 1 & 4 \\ 0 & 3 & 0 & 0 \\ 1 & 0 & 0 & 2 \\ 4 & 0 & 2 & 5 \end{bmatrix}.$$

Now how to define symmetry in case of square mixed bimatrices.

**DEFINITION 1.3.1:** *Let $A_B = A_1^{m \times m} \cup A_2^{n \times n}$ be a square mixed bimatrix, $A_B$ is a symmetric bimatrix if the component matrices $A_1$ and $A_2$ are symmetric matrices. i.e., $A_1 = A'_1$ and $A_2 = A'_2$. Clearly $A_B = A_1^{m \times m} \cup A_2^{n \times n}$ has a total number $\frac{m(m+1)}{2} + \frac{n(n+1)}{2}$ arbitrary elements.*

***Example 1.3.4:*** Let

$$A_B = \begin{bmatrix} 0 & 1 & 2 & 5 \\ 1 & 3 & -1 & 0 \\ 2 & -1 & 2 & -1 \\ 5 & 0 & -1 & -4 \end{bmatrix} \cup \begin{bmatrix} 3 & 0 & 1 \\ 0 & 2 & -1 \\ 1 & -1 & 5 \end{bmatrix}$$

$$A'_B = \begin{bmatrix} 0 & 1 & 2 & 5 \\ 1 & 3 & -1 & 0 \\ 2 & -1 & 2 & -1 \\ 5 & 0 & -1 & -4 \end{bmatrix} \cup \begin{bmatrix} 3 & 0 & 1 \\ 0 & 2 & -1 \\ 1 & -1 & 5 \end{bmatrix}.$$

So $A_B$ is a symmetric mixed square bimatrix. This has $\frac{4(4+1)}{2} + \frac{3(3+1)}{2} = 10 + 6 = 16$ arbitrary elements.

Now we proceed on to define the notion of skew symmetric bimatrix.

**DEFINITION 1.3.2:** *Let $A_B^{m \times m} = A_1 \cup A_2$ be a m × m square bimatrix i.e., $A_1$ and $A_2$ are m × m square matrices. A skew-symmetric bimatrix is a bimatrix $A_B$ for which $A_B = -A^1_B$. where $-A'_B = -A'_1 \cup -A'_2$ i.e., the component matrices $A_1$ and $A_2$ of $A_B$ are also skew symmetric.*

So in a skew symmetric bimatrix $a^1_{ij} = -a^1_{ji}$ and $a^2_{ij} = -a^2_{ji}$ where $A_1 = \left(a^1_{ij}\right)$ and $A_2 = \left(a^2_{ij}\right)$.

Thus in a skew symmetric bimatrix we have the diagonal elements of $A_1$ and $A_2$ are zero i.e., $a^1_{ii} = 0$ and $a^2_{ii} = 0$ and the number of arbitrary elements in an m × m (m[th] order) skew symmetric bimatrix is 2m(m – 1).

***Example 1.3.5:*** Let

$$A_B = \begin{bmatrix} 0 & -1 & 2 \\ 1 & 0 & 3 \\ -2 & -3 & 0 \end{bmatrix} \cup \begin{bmatrix} 0 & 3 & 13 \\ -3 & 0 & -2 \\ -13 & 2 & 0 \end{bmatrix}.$$

This is a 3[rd] order or 3 × 3 skew symmetric bimatrix. This has 12 arbitrary elements.

Now we can also have skew symmetric bimatrices when the bimatrix is a mixed one. Let $A_B = A_1^{m \times m} \cup A_2^{n \times n}$ be a square mixed bimatrix, $A_B$ is said to be a skew symmetric bimatrix if $A_1^{m\times m} = -\left(A_1^{m\times m}\right)'$ and $\left(A_2^{n\times n}\right) = -\left(A_2^{n\times n}\right)'$ we denote this by $A_B = -A'_B$.

Clearly $A_B$ is a skew symmetric bimatrix. The number of elements in a skew symmetric bimatrix is m(m –1) + n(n – 1). The number elements in AB is 12 + 2 = 14.

***Example 1.3.5:*** Let

$$A_B = \begin{bmatrix} 0 & 6 & 1 & 2 \\ -6 & 0 & -1 & 4 \\ -1 & 1 & 0 & 3 \\ -2 & -4 & -3 & 0 \end{bmatrix} \cup \begin{bmatrix} 0 & -1 \\ 1 & 0 \end{bmatrix}.$$

Now we prove the following results:

*RESULT 1:* Let $A_B^{m\times m} = A_1 \cup A_2$ be a square bimatrix. Then $A_B^{m\times m}$ can be written as the sum of a symmetric and a skew symmetric bimatrix. $A_B^{m\times m} = A_1 \cup A_2$. Clearly $A_1$ and $A_2$ are just square matrices, we know if A is any square matrix then

$$A = A + \frac{A'}{2} - \frac{A'}{2}$$

$$A = \frac{A + A'}{2} + \frac{A - A'}{2}$$

so each $A_1$ and $A_2$ can be represented by the above equation i.e.,

$$A_1 = \frac{A_1 + A'_1}{2} + \frac{A_1 - A'_1}{2}$$

and

$$A_2 = \frac{A_2 + A'_2}{2} + \frac{A_2 - A'_2}{2}.$$

Now $A_B^{m\times m} = A_1 \cup A_2$

$$= \frac{A_1 + A'_1}{2} + \frac{A_1 - A'_1}{2} \cup \frac{A_2 + A'_2}{2} + \frac{A_2 - A'_2}{2}$$

$$= \left[\left(\frac{A_1 + A'_1}{2}\right) \cup \left(\frac{A_2 + A'_2}{2}\right)\right] + \left[\left(\frac{A_1 - A'_1}{2}\right) \cup \left(\frac{A_2 - A'_2}{2}\right)\right]$$

is the sum of symmetric bimatrix and skew symmetric bimatrix.

Now we proceed on to define the same type of result in case of square mixed symmetric bimatrix.

*RESULT 2:* Let $A_B$ = $A_1^{n\times n} \cup A_2^{m\times m}$ be a mixed square bimatrix. $A_B$ can be written as a sum of mixed symmetric bimatrix and mixed skew symmetric bimatrix.

*Proof:* Let $A_B$ = $A_1^{n\times n} \cup A_2^{m\times m}$ be a mixed square bimatrix. Now $A_1^{n\times n}$ is a square matrix so let

$$A_1^{n\times n} = A = \frac{A + A'}{2} + \frac{A - A'}{2}$$

where $\frac{A + A'}{2}$ and $\frac{A - A'}{2}$ are symmetric and skew symmetric n × n square matrices.

Let $A_2^{m\times m}$ = B clearly B is a square matrix; So B can be written as a sum of the symmetric and skew symmetric square matrices each of order m i.e.,

$$B = \frac{B + B'}{2} + \frac{B - B'}{2}.$$

So now

$$A_B = \left(\frac{A + A'}{2} + \frac{A - A'}{2}\right) \cup \left(\frac{B + B'}{2} + \frac{B - B'}{2}\right)$$

$$= \left(\frac{A + A'}{2} \cup \frac{B + B'}{2}\right) + \left(\frac{A - A'}{2} \cup \frac{B - B'}{2}\right).$$

Thus $A_B$ is a sum of a mixed square symmetric bimatrix and a mixed square skew symmetric bimatrix.

Now we illustrate this by the following example:

***Example 1.3.7:*** Let

$$A_B = \begin{bmatrix} 3 & 2 & 1 & -1 \\ 0 & 1 & 2 & 0 \\ 3 & 4 & -1 & 0 \\ 0 & 1 & 2 & 1 \end{bmatrix} \cup \begin{bmatrix} 5 & -1 & 0 & 2 \\ 0 & 1 & 2 & 1 \\ 1 & 2 & 0 & 0 \\ 2 & 0 & 1 & 0 \end{bmatrix}$$

be a 4 × 4 square bimatrix. We show that $A_B$ can be represented as the sum of a symmetric and a skew symmetric bimatrix.

Since

$$\begin{bmatrix} 3 & 2 & 1 & -1 \\ 0 & 1 & 2 & 0 \\ 3 & 4 & -1 & 0 \\ 0 & 1 & 2 & 1 \end{bmatrix}$$

$$= \left( \frac{\left( \begin{bmatrix} 3 & 2 & 1 & -1 \\ 0 & 1 & 2 & 0 \\ 3 & 4 & -1 & 0 \\ 0 & 1 & 2 & 1 \end{bmatrix} + \begin{bmatrix} 3 & 0 & 3 & 0 \\ 2 & 1 & 4 & 1 \\ 1 & 2 & -1 & 2 \\ -1 & 0 & 0 & 1 \end{bmatrix} \right)}{2} \right.$$

$$\left. + \left( \frac{\begin{bmatrix} 3 & 2 & 1 & -1 \\ 0 & 1 & 2 & 0 \\ 3 & 4 & -1 & 0 \\ 0 & 1 & 2 & 1 \end{bmatrix} - \begin{bmatrix} 3 & 0 & 3 & 0 \\ 2 & 1 & 4 & 1 \\ 1 & 2 & -1 & 2 \\ -1 & 0 & 0 & 1 \end{bmatrix}}{2} \right) \right)$$

$$= \left( \frac{\begin{bmatrix} 6 & 2 & 4 & -1 \\ 2 & 2 & 6 & 1 \\ 4 & 6 & -2 & 2 \\ -1 & 1 & 2 & 2 \end{bmatrix} + \begin{bmatrix} 0 & 2 & -2 & -1 \\ -2 & 0 & -2 & -1 \\ 2 & 2 & 0 & -2 \\ 1 & 1 & 2 & 0 \end{bmatrix}}{2} \right).$$

Now

$$\begin{bmatrix} 5 & -1 & 0 & 2 \\ 0 & 1 & 2 & 1 \\ 1 & 2 & 0 & 0 \\ 2 & 0 & 1 & 0 \end{bmatrix}$$

$$= \left\{ \left( \frac{\begin{bmatrix} 5 & -1 & 0 & 2 \\ 0 & 1 & 2 & 1 \\ 1 & 2 & 0 & 0 \\ 2 & 0 & 1 & 0 \end{bmatrix} + \begin{bmatrix} 5 & -1 & 0 & 2 \\ 0 & 1 & 2 & 1 \\ 1 & 2 & 0 & 0 \\ 2 & 0 & 1 & 0 \end{bmatrix}}{2} \right) \right.$$

$$\left. + \left( \frac{\begin{bmatrix} 5 & -1 & 0 & 2 \\ 0 & 1 & 2 & 1 \\ 1 & 2 & 0 & 0 \\ 2 & 0 & 1 & 0 \end{bmatrix} - \begin{bmatrix} 5 & -1 & 0 & 2 \\ 0 & 1 & 2 & 1 \\ 1 & 2 & 0 & 0 \\ 2 & 0 & 1 & 0 \end{bmatrix}}{2} \right) \right\}$$

$$= \left( \frac{\begin{bmatrix} 10 & -1 & 1 & 4 \\ -1 & 2 & 4 & 1 \\ 1 & 4 & 0 & 1 \\ 4 & 1 & 1 & 0 \end{bmatrix}}{2} + \frac{\begin{bmatrix} 0 & -1 & -1 & 0 \\ 1 & 0 & 0 & 1 \\ 1 & 0 & 0 & -1 \\ 0 & -1 & 1 & 0 \end{bmatrix}}{2} \right).$$

Now

$$A_B = \left\{ \left( \frac{\begin{bmatrix} 6 & 2 & 4 & -1 \\ 2 & 2 & 6 & 1 \\ 4 & 6 & -2 & 2 \\ -1 & 1 & 2 & 2 \end{bmatrix}}{2} \cup \frac{\begin{bmatrix} 10 & -1 & 1 & 4 \\ -1 & 2 & 4 & 1 \\ 1 & 4 & 0 & 1 \\ 4 & 1 & 1 & 0 \end{bmatrix}}{2} \right) \right.$$

$$+ \left( \frac{\begin{bmatrix} 0 & 2 & -2 & -1 \\ -2 & 0 & -2 & -1 \\ 2 & 2 & 0 & -2 \\ 1 & 1 & 2 & 0 \end{bmatrix}}{2} \cup \frac{\begin{bmatrix} 0 & -1 & -1 & 0 \\ 1 & 0 & 0 & 1 \\ 1 & 0 & 0 & -1 \\ 0 & -1 & 1 & 0 \end{bmatrix}}{2} \right) \Bigg\}.$$

Thus $A_B$ is the sum of the symmetric bimatrix and the skew symmetric bimatrix.

Now we give an example for the mixed square bimatrix i.e., to show that the mixed square bimatrix can also be represented as the sum of the mixed symmetric bimatrix and the mixed skew symmetric bimatrix.

***Example 1.3.8:*** Let

$$A_B = \begin{bmatrix} 6 & 2 & 5 & -1 & 7 \\ 0 & 1 & 0 & 0 & 2 \\ 2 & 0 & 0 & 6 & 1 \\ -3 & 0 & 4 & 5 & 0 \\ 2 & -1 & 0 & -1 & 6 \end{bmatrix} \cup \begin{bmatrix} 3 & 0 & 2 \\ 1 & 1 & 1 \\ 4 & -2 & 6 \end{bmatrix}$$

be the mixed square bimatrix. Let

$$A = \begin{bmatrix} 6 & 2 & 5 & -1 & 7 \\ 0 & 1 & 0 & 0 & 2 \\ 2 & 0 & 0 & 6 & 1 \\ -3 & 0 & 4 & 5 & 0 \\ 2 & -1 & 0 & -1 & 6 \end{bmatrix}$$

and

$$B = \begin{bmatrix} 3 & 0 & 2 \\ 1 & 1 & 1 \\ 4 & -2 & 6 \end{bmatrix}$$

both of the component matrices are square matrices so can be written as the sum of a symmetric and a skew symmetric matrix.

$$A = \left\{ \frac{\begin{bmatrix} 6 & 2 & 5 & -1 & 7 \\ 0 & 1 & 0 & 0 & 2 \\ 2 & 0 & 0 & 6 & 1 \\ -3 & 0 & 4 & 5 & 0 \\ 2 & -1 & 0 & -1 & 6 \end{bmatrix} + \begin{bmatrix} 6 & 0 & 2 & -3 & 2 \\ 2 & 1 & 0 & 0 & -1 \\ 5 & 0 & 0 & 4 & 0 \\ -1 & 0 & 6 & 5 & -1 \\ 7 & 2 & 1 & 0 & 6 \end{bmatrix}}{2} \right.$$

$$\left. + \left( \frac{\begin{bmatrix} 6 & 2 & 5 & -1 & 7 \\ 0 & 1 & 0 & 0 & 2 \\ 2 & 0 & 0 & 6 & 1 \\ -3 & 0 & 4 & 5 & 0 \\ 2 & -1 & 0 & -1 & 6 \end{bmatrix} - \begin{bmatrix} 6 & 0 & 2 & -3 & 2 \\ 2 & 1 & 0 & 0 & -1 \\ 5 & 0 & 0 & 4 & 0 \\ -1 & 0 & 6 & 5 & -1 \\ 7 & 2 & 1 & 0 & 6 \end{bmatrix}}{2} \right) \right\}$$

$$= \left\{ \frac{\begin{bmatrix} 12 & 2 & 7 & -4 & 9 \\ 2 & 2 & 0 & 0 & 1 \\ 7 & 0 & 0 & 10 & 1 \\ -4 & 0 & 10 & 10 & -1 \\ 9 & 1 & 1 & -1 & 12 \end{bmatrix}}{2} + \frac{\begin{bmatrix} 0 & 2 & 3 & 2 & 5 \\ -2 & 0 & 0 & 0 & 3 \\ -3 & 0 & 0 & 2 & 1 \\ -2 & 0 & -2 & 0 & 1 \\ -5 & -3 & -1 & -1 & 0 \end{bmatrix}}{2} \right\}.$$

$$B = \left\{ \left( \frac{\begin{bmatrix} 3 & 0 & 2 \\ 1 & 1 & 1 \\ 4 & -2 & 6 \end{bmatrix} + \begin{bmatrix} 3 & 1 & 4 \\ 0 & 1 & -2 \\ 2 & 1 & 6 \end{bmatrix}}{2} \right) \right.$$

$$\left. + \left( \frac{\begin{bmatrix} 3 & 0 & 2 \\ 1 & 1 & 1 \\ 4 & -2 & 6 \end{bmatrix} - \begin{bmatrix} 3 & 1 & 4 \\ 0 & 1 & -2 \\ 2 & 1 & 6 \end{bmatrix}}{2} \right) \right\}$$

$$= \left\{ \frac{\begin{bmatrix} 6 & 1 & 6 \\ 1 & 2 & -1 \\ 6 & -1 & 12 \end{bmatrix}}{2} + \frac{\begin{bmatrix} 0 & -1 & -2 \\ 1 & 0 & 3 \\ 2 & -3 & 0 \end{bmatrix}}{2} \right\}$$

So

$$A_B = \left\{ \frac{\begin{bmatrix} 12 & 2 & 7 & -4 & 9 \\ 2 & 2 & 0 & 0 & 1 \\ 7 & 0 & 0 & 10 & 1 \\ -4 & 0 & 10 & 10 & -1 \\ 9 & 1 & 1 & -1 & 12 \end{bmatrix}}{2} + \frac{\begin{bmatrix} 0 & 2 & 3 & 2 & 5 \\ -2 & 0 & 0 & 0 & 3 \\ -3 & 0 & 0 & 2 & 1 \\ -2 & 0 & -2 & 0 & 1 \\ -5 & -3 & -1 & -1 & 0 \end{bmatrix}}{2} \right.$$

$$\left. \cup \ \frac{\begin{bmatrix} 6 & 1 & 6 \\ 1 & 2 & -1 \\ 6 & -1 & 12 \end{bmatrix}}{2} + \frac{\begin{bmatrix} 0 & -1 & -2 \\ 1 & 0 & 3 \\ 2 & -3 & 0 \end{bmatrix}}{2} \right\}$$

$$= \left\{ \frac{\begin{bmatrix} 12 & 2 & 7 & -4 & 9 \\ 2 & 2 & 0 & 0 & 1 \\ 7 & 0 & 0 & 10 & 1 \\ -4 & 0 & 10 & 10 & -1 \\ 9 & 1 & 1 & -1 & 12 \end{bmatrix}}{2} \cup \frac{\begin{bmatrix} 6 & 1 & 6 \\ 1 & 2 & -1 \\ 6 & -1 & 12 \end{bmatrix}}{2} \right.$$

$$\left. + \frac{\begin{bmatrix} 0 & 2 & 3 & 2 & 5 \\ -2 & 0 & 0 & 0 & 3 \\ -3 & 0 & 0 & 2 & 1 \\ -2 & 0 & -2 & 0 & 1 \\ -5 & -3 & -1 & -1 & 0 \end{bmatrix}}{2} \cup \frac{\begin{bmatrix} 0 & -1 & -2 \\ 1 & 0 & 3 \\ 2 & -3 & 0 \end{bmatrix}}{2} \right\}$$

= mixed symmetric bimatrix + mixed skew symmetric bimatrix i.e. a mixed square matrix $A_B$ is written as a sum of the mixed symmetric bimatrix and mixed skew symmetric bimatrix. Hence the result.

Now we proceed on to define the notion of partitioning of bimatrices. It is interesting to study some subset of elements of a bimatrix which form the subbimatrix which is given in the following section.

## 1.4 Subbimatrix

Here we introduce the notion of subbimatrix and give some properties about them.

**DEFINITION 1.4.1:** *Let $A_B$ be any bimatrix i.e., $A_B = A_1^{m\times n} \cup A_2^{p\times q}$. If we cross out all but $k_1$ rows and $s_1$ columns of the $m \times n$ matrix $A_1$ and cross out all but $k_2$ rows and $s_2$ columns of the $p \times q$ matrix $A_2$ the resulting $k_1 \times s_1$ and $k_2 \times s_2$ bimatrix is called a subbimatrix of $A_B$.*

***Example 1.4.1:*** Let

$$A_B = \begin{bmatrix} 3 & 2 & 1 & 4 \\ 6 & 0 & 1 & 2 \\ -1 & 6 & -1 & 0 \end{bmatrix} \cup \begin{bmatrix} 3 & 8 & 3 & 6 & -2 \\ 0 & 0 & 1 & 0 & 2 \\ 1 & 1 & 0 & 0 & 1 \\ 0 & 0 & 1 & 2 & 3 \\ 2 & 1 & 0 & -1 & 3 \\ -1 & 4 & 0 & 0 & 2 \end{bmatrix}$$

$= A_1 \cup A_2$.

Then a submatrix of $A_B$ is given by

$$\begin{bmatrix} 3 & 2 & 1 \\ -1 & 6 & 1 \end{bmatrix} \cup \begin{bmatrix} 3 & 8 \\ 1 & 1 \\ 2 & 1 \\ -1 & 4 \end{bmatrix}.$$

where one row and one column is crossed out in $A_1$ and 2 rows and 3 columns is crossed out in $A_2$.

We for several reasons introduce the notion of partition bimatrices in subbimatrices. Some of the main reasons are

(1) The partitioning may simplify the writing or printing of $A_B$.
(2) It exhibits some particular structure of $A_B$ which is of interest.
(3) It simplifies computation.

***Example 1.4.2:*** Let $A_B = A_1 \cup A_2$ be a bimatrix i.e., $A_B$

$$= \left[\begin{array}{ccc|c} a_{11} & a_{12} & \cdots & a_{16} \\ a_{21} & a_{22} & & a_{26} \\ \hline a_{31} & a_{32} & \cdots & a_{36} \\ \vdots & & & \\ a_{91} & a_{92} & & a_{96} \end{array}\right] \cup \left[\begin{array}{cc|ccc} b_{11} & b_{12} & b_{13} & b_{14} & b_{15} \\ b_{21} & b_{22} & b_{23} & b_{24} & b_{25} \\ \hline b_{31} & b_{32} & b_{33} & b_{34} & b_{35} \end{array}\right].$$

Now unlike usual matrices we have the following partitions if we imagine $A_B$ to be divided up by lines as shown.

Now the bimatrix is partitioned into

$$A^1_B = \begin{bmatrix} a_{11} & a_{12} & \ldots & a_{15} \\ a_{21} & a_{22} & \ldots & a_{25} \end{bmatrix} \cup \begin{bmatrix} b_{11} & b_{12} \\ b_{21} & b_{22} \end{bmatrix}$$

$$A_B^{11} = \begin{bmatrix} a_{11} & a_{12} & \ldots & a_{15} \\ a_{21} & a_{22} & \ldots & a_{25} \end{bmatrix} \cup \begin{bmatrix} b_{13} & b_{14} & b_{15} \\ b_{21} & b_{24} & b_{25} \end{bmatrix}$$

$$A_B^3 = \begin{bmatrix} a_{11} & a_{12} & \ldots & a_{15} \\ a_{21} & a_{22} & \ldots & a_{25} \end{bmatrix} \cup \begin{bmatrix} b_{31} & b_{32} \end{bmatrix}$$

$$A_B^4 = \begin{bmatrix} a_{11} & a_{12} & \ldots & a_{15} \\ a_{21} & a_{22} & \ldots & a_{25} \end{bmatrix} \cup [b_{33}\ b_{34}\ b_{35}]$$

$$A_B^5 = \begin{bmatrix} a_{16} \\ a_{26} \end{bmatrix} \cup \begin{bmatrix} b_{11} & b_{12} \\ b_{21} & b_{22} \end{bmatrix}$$

$$A_B^6 = \begin{bmatrix} a_{16} \\ a_{26} \end{bmatrix} \cup \begin{bmatrix} b_{13} & b_{14} & b_{15} \\ b_{23} & b_{24} & b_{25} \end{bmatrix}$$

$$A_B^7 = \begin{bmatrix} a_{16} \\ a_{26} \end{bmatrix} \cup [b_{31}\ \ b_{32}]$$

$$A_B^8 = \begin{bmatrix} a_{16} \\ a_{26} \end{bmatrix} \cup [b_{33}\ \ b_{34}\ \ b_{35}]$$

$$A_B^9 = \begin{bmatrix} a_{31} & a_{32} & \cdots & a_{35} \\ a_{41} & a_{42} & \cdots & a_{45} \\ & & & \\ a_{91} & a_{92} & & a_{95} \end{bmatrix} \cup \begin{bmatrix} b_{11} & b_{12} \\ b_{21} & b_{22} \end{bmatrix}$$

$$A_B^{10} = \begin{bmatrix} a_{31} & a_{32} & \cdots & a_{35} \\ a_{41} & a_{42} & \cdots & a_{45} \\ \vdots & \vdots & & \vdots \\ a_{91} & a_{92} & \cdots & a_{95} \end{bmatrix} \cup \begin{bmatrix} b_{13} & b_{14} & b_{15} \\ b_{21} & b_{24} & b_{25} \end{bmatrix}$$

$$A_B^{11} = \begin{bmatrix} a_{31} & a_{32} & \cdots & a_{35} \\ a_{41} & a_{42} & \cdots & a_{45} \\ \vdots & \vdots & & \vdots \\ a_{91} & a_{92} & \cdots & a_{95} \end{bmatrix} \cup \begin{bmatrix} b_{31} & b_{32} \end{bmatrix}$$

$$A_B^{12} = \begin{bmatrix} a_{31} & a_{32} & \cdots & a_{35} \\ a_{41} & a_{42} & \cdots & a_{45} \\ \vdots & \vdots & & \vdots \\ a_{91} & a_{92} & \cdots & a_{95} \end{bmatrix} \cup \begin{bmatrix} b_{33} & b_{34} & b_{35} \end{bmatrix}$$

$$A_B^{13} = \begin{bmatrix} a_{36} \\ \vdots \\ a_{96} \end{bmatrix} \cup \begin{bmatrix} b_{11} & b_{12} \\ b_{21} & b_{22} \end{bmatrix}$$

$$A_B^{14} = \begin{bmatrix} a_{36} \\ \vdots \\ a_{96} \end{bmatrix} \cup \begin{bmatrix} b_{13} & b_{14} & b_{15} \\ b_{23} & b_{24} & b_{25} \end{bmatrix}$$

$$A_B^{15} = \begin{bmatrix} a_{36} \\ \vdots \\ a_{96} \end{bmatrix} \cup \begin{bmatrix} b_{31} & b_{32} \end{bmatrix}$$

and

$$A_B^{14} = \begin{bmatrix} a_{36} \\ \vdots \\ a_{96} \end{bmatrix} \cup \begin{bmatrix} b_{33} & b_{34} & b_{35} \end{bmatrix}.$$

With the lines shown for $A_B$ there are 16 partition some mixed bimatrices some mixed square matrices and so on.

It is pertinent to mention here that the rule for addition of partitioned bimatrices is the same as the rule for addition of ordinary bimatrices if the subbimatrices are conformable for addition. In other words the matrices can be added by blocks. Addition of bimatrices $A_B$ and $B_B$ when partitioned is possible only when $A_B$ and $B_B$ are of the same type i.e., both are m × n rectangular bimatrix or both $A_B$ and $B_B$ are both square bimatrix or both $A_B$ and $B_B$ are both mixed square matrices or both mixed rectangular matrices with computable addition.

***Example 1.4.3:*** Let

$$A_B = \left[\begin{array}{cc|cc} 3 & 0 & 1 & 2 \\ 1 & 1 & 0 & 3 \\ \hline 0 & 1 & 0 & 2 \end{array}\right] \cup \left[\begin{array}{cc|c} 3 & 0 & 1 \\ 0 & 1 & 1 \\ 2 & 1 & 0 \\ \hline 0 & 0 & 5 \end{array}\right]$$

$$B_B = \left[\begin{array}{cc|cc} 5 & -1 & 3 & 2 \\ -1 & 0 & 0 & 5 \\ \hline 1 & 1 & 5 & -2 \end{array}\right] \cup \left[\begin{array}{cc|c} 6 & 2 & -3 \\ 0 & 1 & 2 \\ 1 & 1 & 6 \\ \hline 0 & -1 & 6 \end{array}\right]$$

$A_B + B_B$ as block sum of bimatrices is given by

$$A_B + B_B = \left[\begin{array}{cc|cc} 8 & -1 & 4 & 4 \\ 0 & 1 & 0 & 8 \\ \hline 1 & 2 & 5 & 0 \end{array}\right] \cup \left[\begin{array}{cc|c} 9 & 2 & -2 \\ 0 & 2 & 3 \\ 3 & 2 & 6 \\ \hline 0 & -1 & 11 \end{array}\right].$$

***Example 1.4.4:*** Let

$$A_B = \left[\begin{array}{ccc|cc} 3 & 1 & 1 & 2 & 5 \\ 0 & 1 & 0 & 0 & -1 \\ 1 & 2 & 3 & 0 & 0 \\ 1 & 1 & 1 & 0 & 1 \\ \hline 0 & -1 & 2 & 0 & -2 \end{array}\right] \cup \left[\begin{array}{c|c} 3 & 1 \\ 2 & 1 \\ \hline 2 & 2 \end{array}\right]$$

$$B_B = \left[\begin{array}{ccc|cc} -1 & 0 & 0 & 0 & 1 \\ 2 & 1 & 1 & 0 & 0 \\ -4 & 0 & 1 & 0 & 1 \\ 2 & 1 & 2 & 0 & 1 \\ \hline 1 & 1 & 1 & 1 & 1 \end{array}\right] \cup \left[\begin{array}{c|c} 0 & -1 \\ 2 & 1 \\ \hline 3 & 4 \end{array}\right].$$

Now if the bimatrices are divided into blocks by these lines then the addition of these can be carried out for both mixed bimatrices are compatible with respect to addition and the block division also a compatible or a similar division.

$$A_B + B_B = \left[\begin{array}{ccc|cc} 2 & 1 & 1 & 2 & 6 \\ 2 & 2 & 1 & 0 & -1 \\ -3 & 2 & 4 & 0 & 1 \\ 3 & 2 & 3 & 0 & 2 \\ \hline 1 & 0 & 3 & 1 & -1 \end{array}\right] \cup \left[\begin{array}{c|c} 3 & 0 \\ 4 & 2 \\ \hline 5 & 6 \end{array}\right].$$

So the sum of the bimatrices are also divided into blocks.

***Example 1.4.5:*** Let

$$A_B = \left[\begin{array}{cc|c} 2 & 1 & 3 \\ 0 & 1 & 1 \\ 9 & 2 & -1 \\ \hline 3 & 1 & 2 \end{array}\right] \cup \left[\begin{array}{cc|c} 3 & 5 & 1 \\ \hline 5 & 4 & 0 \end{array}\right]$$

and

$$B_B = \left[\begin{array}{cc|c} -1 & 0 & 1 \\ 1 & 4 & 2 \\ \hline 2 & 2 & 2 \\ 3 & -1 & 0 \end{array}\right] \cup \left[\begin{array}{c|cc} 0 & 1 & 7 \\ \hline 0 & 3 & 8 \end{array}\right].$$

Clearly $A_B + B_B$ cannot be added as blocks for the sum of these two bimatrices exist but as block the addition is not compatible.

Now we proceed on to define bimatrix multiplication as block bimatrices.

***Example 1.4.6:*** Let

$$A_B = \left[\begin{array}{cc|cc} 3 & 2 & -1 & 0 \\ 0 & 1 & 0 & 1 \\ \hline 1 & 0 & 0 & 0 \\ 5 & -1 & 0 & 0 \end{array}\right] \cup \left[\begin{array}{cc|cc} 4 & 0 & 0 & 0 \\ 0 & 2 & 1 & 0 \\ \hline 1 & 0 & 1 & 0 \\ 2 & 1 & 1 & 0 \end{array}\right]$$

$$B_B = \left[\begin{array}{cc|cc} 3 & 0 & 0 & 6 \\ -1 & 0 & 1 & 0 \\ \hline 0 & -1 & 0 & -1 \\ 2 & 1 & 0 & 1 \end{array}\right] \cup \left[\begin{array}{cc|cc} 1 & 0 & -1 & 0 \\ 0 & 1 & 0 & -1 \\ \hline 1 & 0 & 1 & 1 \\ 0 & 2 & 0 & 0 \end{array}\right].$$

Be two $4 \times 4$ square bimatrices and be divided into blocks as shown by the lines i.e.

$$A_B = \begin{bmatrix} A_1 & A_2 \\ A_3 & A_4 \end{bmatrix} \cup \begin{bmatrix} A_1^1 & A_2^1 \\ A_3^1 & A_4^1 \end{bmatrix}$$

and

$$B_B = \begin{bmatrix} B_1 & B_2 \\ B_3 & B_4 \end{bmatrix} \cup \begin{bmatrix} B_1^1 & B_2^1 \\ B_3^1 & B_4^1 \end{bmatrix}$$

$A_B\ B_B =$

$$\begin{bmatrix} A_1B_1 + A_2B_2 & A_1B_2 + A_2B_4 \\ A_3B_1 + A_4B_2 & A_3B_2 + A_4B_4 \end{bmatrix} \cup \begin{bmatrix} A_1^1B_1^1 + A_2^1B_2^1 & A_1^1B_2^1 + A_2^1B_4^1 \\ A_3^1B_1^1 + A_4^1B_2^1 & A_3^1B_2^1 + A_4^1B_4^1 \end{bmatrix}$$

is again the bimatrix which is partitioned. It is to be noted all bimatrices in general need be compatible with block multiplication. This sort of block multiplication can be defined even in the case of mixed square bimatrices.

Now we proceed on to interrupt the study of bimatrices for a while in order to develop some aspects of the theory of determinants which will be needed in our further investigations of the properties of bimatrices.

## 1.5 Basic concept of bideterminant in case of bimatrices

We in this section proceed onto define the concept of determinant of a bimatrix and derive some of properties analogous to matrices which we bideterminant of a bimatrix.

Let $A_B = A_1 \cup A_2$ be a square bimatrix. The bideterminant of a square bimatrix is an ordered pair $(d_1, d_2)$ where $d_1 = |A_1|$ and $d_2 = |A_2|$. $|A_B| = (d_1, d_2)$ where $d_1$ and $d_2$ are reals may be positive or negative or even zero. ( $|A|$ denotes determinant of A).

***Example 1.5.1:*** Let

$$A_B = \begin{bmatrix} 3 & 0 & 0 \\ 2 & 1 & 1 \\ 0 & 1 & 1 \end{bmatrix} \cup \begin{bmatrix} 4 & 5 \\ -2 & 0 \end{bmatrix}.$$

The bideterminant of this bimatrix is the pair (0, 10) and denoted by $|A_B|$.

***Example 1.5.2:*** Let

$$A_B = \begin{bmatrix} 3 & 0 & 1 & 2 \\ 0 & 1 & 0 & 0 \\ 0 & 2 & 1 & 5 \\ 0 & 0 & 1 & 2 \end{bmatrix} \cup \begin{bmatrix} 0 & 5 & -1 & 0 \\ 1 & 0 & 0 & 0 \\ 0 & 1 & 0 & 4 \\ 1 & 1 & 1 & 1 \end{bmatrix}$$

be the 4 × 4 square bimatrix. Then the bideterminant of $|A_B|$ = (–9, –23).

Thus we give the key to the general definition of a bideterminant in case of a n × n square bimatrix $A_B^{n \times n} = A_1 \cup A_2$.

The determinant of an n[th] order bimatrix $A_B = (a_{ij}) \cup (b_{ij})$ written as $|A_B|$, is defined to be the pair of number computed from the following sum involving $n^2 + n^2$ elements in $A_B$.

$$|A_B| = \Sigma\ (\pm)\ (a_{1i}\ a_{2j}\ \ldots\ a_{nr}),\ +\ \Sigma\ (\pm)\ (b_{1i}\ b_{2j} \ldots\ b_{nr})$$

the sum being taken over all permutations of the second subscripts. A term is assigned as plus sign if (i, j, … , r) is an even permutation of (1, 2, …, n) and a minus sign if it is an odd permutation. We shall find it convenient to refer to the bideterminant of the n[th] order bimatrix or n × n bideterminant. In case of mixed square bimatrix it will be little modified.

Now we proceed on to give some properties of bideterminants. If $A_B = A_1 \cup A_2$ where $A_B$ is a square n × n bimatrix the bideterminant of $A_B$ is denoted by

$$|A_B| = |A_1| \cup |A_2| = (d_1, d_2) = (\ |A_1|, |A_2|\ ).$$

Now we proceed on to define the notion of expansion of cofactors and the notion of bicofactors in square bimatrices of order n. Let $A_B = A_1 \cup A_2$ be a square bimatrix, of order n the bicofactor of $A_B$ denoted $(A_B)_{ij}$ can be considered as a bideterminant of order n – 1. All elements of $A_B$ appear in $(A_B)_{ij}$ except those in row i and column j. In fact except for a possible difference in sign $(A_B)_{ij}$ is the determinant of the subbimatrix formed from $A_B$ by crossing out row i and column j, if the reader examines the $(A_B)_{ij}$ .

Let us now determine the sign that should be assigned to the determinant of the subbimatrix obtained by crossing out row i and column j in matrices $A_1$ and $A_2$ of $A_B$, in order to convert it to $(A_B)_{ij}$. We do this most simply by moving row i to row 1 by i-1 inter changes of rows in both $A_1$ and $A_2$. The other rows while retaining their original order will be moved down one. Then by j – 1 interchanges column j from $A_1$ and $A_2$ will be moved to the position of column 1. We shall call this new matrix $B_B$.

Then

$$|B_B| = (-1)^{i+j-2} |A_B| = (-1)^{i+j} |A_B|.$$

However the product of the elements on the main diagonal of $B_B$ are the same as those appearing on the main diagonal of the subbimatrix whose determinant is $(A_B)_{ij}$. From the above equation, $(-1)^{i+j}$ times this term in the expansion of $|A_B|$ must be positive. Therefore $(A_B)_{ij}$ is $(-1)^{i+j}$ times the determinant of the subbimatrix formed from $A_B$ by deleting row i and column j of both $A_1$ and $A_2$. The $(A_B)_{ij}$ is called the bicofactor of $(a^1_{ij}) \cup (a^2_{ij})$ .

Now the bicofactor $(A_B)_{ij}$ can be defined for $A_B = A_1 \cup A_2$ as follows.

**DEFINITION 1.5.1:** *Let $A_B = A_1 \cup A_2$ be a $n \times n$ square bimatrix. The bicofactor $(A_B)_{ij} = (A_1)_{ij} \cup (A_2)_{ij}$ of the element $a^1_{ij} \cup a^2_{ij}$ of the bimatrix $A_B$ is $(-1)^{i+j}$ times the bideterminant subbimatrix obtained from $A_B$ by deleting row i and column j in both $A_1$ and $A_2$.*

We now provide a way to evaluate the $(A_B)_{ij}$ of $|A_B|$. This important method of evaluating bideterminants is called expansion by bicofactors.

Now we proceed on to define biminor of order k. For any n × n bimatrix $A_B = A_1 \cup A_2$ consider the $k^{th}$ order subbimatrix $R = R_1 \cup R_2$ obtained by deleting all but some k-rows and k-columns of $A_1$ and $A_2$ ($A_B = A_1 \cup A_2$) then |R| is called a $k^{th}$ order biminor of $A_B$. We could write $|A_B|$ as

$$\sum_{i=1}^{n} a_{ij}^1 A_{ij}^1 \cup \sum_{i=1}^{n} a_{ij}^2 A_{ij}^2$$

where $A_{ij}^1$ and $A_{ij}^2$ are again the cofactors of $a_{ij}^1$ and $a_{ij}^2$ respectively. This is called an expansion of $|A_B|$ by column j of $A_B = A_1 \cup A_2$.

Now we proceed onto define the notion of additional properties of bideterminants. Expansion by bicofactors can be used to prove some additional properties of bideterminants.

Consider

$$A_B = \begin{pmatrix} \lambda_1^1 a_{11}^1 + \lambda_2^1 b_{11}^1 + \lambda_3^1 c_{11}^1 & a_{12}^1 & \cdots & a_{1n}^1 \\ \lambda_1^1 a_{21}^1 + \lambda_2^1 b_{21}^1 + \lambda_3^1 c_{21}^1 & a_{22}^1 & \cdots & a_{2n}^1 \\ \vdots & \vdots & & \vdots \\ \lambda_1^1 a_{n1}^1 + \lambda_2^1 b_{n1}^1 + \lambda_3^1 c_{n1}^1 & a_{n2}^1 & \cdots & a_{nn}^1 \end{pmatrix}$$

$$\cup \begin{pmatrix} \lambda_1^2 a_{11}^2 + \lambda_2^2 b_{11}^2 + \lambda_3^2 c_{11}^2 & a_{12}^2 & \cdots & a_{1n}^2 \\ \lambda_1^2 a_{21}^2 + \lambda_2^2 b_{21}^2 + \lambda_3^2 c_{21}^2 & a_{22}^2 & \cdots & a_{2n}^2 \\ \vdots & \vdots & & \\ \lambda_1^2 a_{n1}^2 + \lambda_2^2 b_{n1}^2 + \lambda_3^2 c_{n1}^2 & a_{n2}^2 & \dots & a_{nn}^2 \end{pmatrix}$$

$$= \left( \lambda_1^1 \begin{vmatrix} a_{11}^1 & a_{12}^1 & \cdots & a_{1n}^1 \\ a_{21}^1 & a_{22}^1 & \cdots & a_{2n}^1 \\ \vdots & & & \\ a_{n1}^1 & a_{n2}^1 & \cdots & a_{nn}^1 \end{vmatrix} \cup \lambda_2^1 \begin{vmatrix} a_{11}^2 & a_{12}^2 & \cdots & a_{1n}^2 \\ a_{21}^2 & a_{22}^2 & \cdots & a_{2n}^2 \\ \vdots & \vdots & & \\ a_{n1}^2 & a_{n2}^2 & \cdots & a_{nn}^2 \end{vmatrix} \right)$$

$$+ \left( \lambda_2^1 \begin{vmatrix} b_{11}^1 & a_{12}^1 & \cdots & a_{1n}^1 \\ b_{21}^1 & a_{22}^1 & \cdots & a_{2n}^1 \\ \vdots & & & \\ b_{n1}^1 & a_{n2}^1 & \cdots & a_{nn}^1 \end{vmatrix} \cup \lambda_2^2 \begin{vmatrix} b_{11}^2 & a_{12}^2 & \cdots & a_{1n}^2 \\ a_{21}^2 & a_{22}^2 & \cdots & a_{2n}^2 \\ \vdots & \vdots & & \\ a_{n1}^2 & a_{n2}^2 & \cdots & a_{nn}^2 \end{vmatrix} \right)$$

$$+ \left( \lambda_3^1 \begin{vmatrix} c_{11}^1 & a_{12}^1 & \cdots & a_{1n}^1 \\ c_{21}^1 & a_{22}^1 & \cdots & a_{2n}^1 \\ \vdots & & & \\ c_{n1}^1 & a_{n2}^1 & \cdots & a_{nn}^1 \end{vmatrix} \cup \lambda_3^2 \begin{vmatrix} c_{11}^2 & a_{12}^2 & \cdots & a_{1n}^2 \\ c_{21}^2 & a_{22}^2 & \cdots & a_{2n}^2 \\ \vdots & \vdots & & \\ a_{n1}^2 & a_{n2}^2 & \cdots & a_{nn}^2 \end{vmatrix} \right)$$

since

$$\sum \left( \lambda_1^1 a_{1i}^1 + \lambda_2^1 b_{1i}^1 + \lambda_3^1 c_{1i}^1 \right) A_{i1}^1 \cup \sum \left( \lambda_1^2 a_{1i}^2 + \lambda_2^2 b_{1i}^2 + \lambda_3^2 c_{1i}^2 \right) A_{i1}^2$$

$$= \lambda_1^1 \sum a_{i1}^1 A_{i1}^1 \cup \lambda_1^2 \sum a_{i1}^2 A_{i1}^2 + \lambda_2^1 \sum b_{i1}^1 A_{i1}^1 \cup \lambda_2^2 \sum b_{i1}^2 A_{i1}^2 +$$
$$\lambda_3^1 \sum c_{i1}^1 A_{i1}^1 \cup \lambda_3^2 \sum c_{i1}^2 A_{i1}^2 .$$

The same sort of result holds when the i$^{th}$ row or j$^{th}$ column of $A_B$ is written as a sum of terms.

***Example 1.5.3:*** Let

$$A_B = \begin{bmatrix} 8 & 2 & 3 \\ 5 & 1 & 3 \\ 0 & 1 & 0 \end{bmatrix} \cup \begin{bmatrix} 9 & 2 & 0 & 1 \\ 8 & 1 & 0 & 1 \\ 1 & 2 & 0 & -1 \\ 0 & 1 & 1 & 2 \end{bmatrix}$$

be the

$$\left| A_B \right| = \begin{vmatrix} 2 & + & 2 & + & 2 & 2 & 3 \\ 2 & + & 2 & + & 1 & 1 & 3 \\ 0 & + & 0 & + & 0 & 1 & 0 \end{vmatrix}$$

$$\cup \begin{vmatrix} 3 & + & 3 & + & 3 & 2 & 0 & 1 \\ 2 & + & 2 & + & 2 & 1 & 0 & 1 \\ 1 & + & 0 & + & 0 & 2 & 0 & -1 \\ 0 & + & 0 & + & 0 & 1 & 1 & 2 \end{vmatrix}$$

$$= \left\{ \begin{vmatrix} 2 & 2 & 3 \\ 2 & 1 & 3 \\ 0 & 1 & 0 \end{vmatrix} + \begin{vmatrix} 2 & 2 & 3 \\ 2 & 1 & 3 \\ 0 & 1 & 0 \end{vmatrix} + \begin{vmatrix} 2 & 2 & 3 \\ 1 & 1 & 3 \\ 0 & 1 & 0 \end{vmatrix} \right\}$$

$$\cup \left\{ \begin{vmatrix} 3 & 2 & 0 & 1 \\ 2 & 1 & 0 & 1 \\ 1 & 2 & 0 & -1 \\ 0 & 1 & 1 & 2 \end{vmatrix} + \begin{vmatrix} 3 & 2 & 0 & 1 \\ 2 & 1 & 0 & 1 \\ 0 & 2 & 0 & -1 \\ 0 & 1 & 1 & 2 \end{vmatrix} + \begin{vmatrix} 3 & 2 & 0 & 1 \\ 2 & 1 & 0 & 1 \\ 0 & 2 & 0 & -1 \\ 0 & 1 & 1 & 2 \end{vmatrix} \right\}$$

$$= \left\{ \begin{vmatrix} 2 & 2 & 3 \\ 2 & 1 & 3 \\ 0 & 1 & 0 \end{vmatrix} \cup \begin{vmatrix} 3 & 2 & 0 & 1 \\ 2 & 1 & 0 & 1 \\ 1 & 2 & 0 & -1 \\ 0 & 1 & 1 & 2 \end{vmatrix} \right\}$$

$$+ \left\{ \begin{vmatrix} 2 & 2 & 3 \\ 2 & 1 & 3 \\ 0 & 1 & 0 \end{vmatrix} \cup \begin{vmatrix} 3 & 2 & 0 & 1 \\ 2 & 1 & 0 & 1 \\ 1 & 2 & 0 & -1 \\ 0 & 1 & 1 & 2 \end{vmatrix} \right\}$$

$$+ \left\{ \begin{vmatrix} 2 & 2 & 3 \\ 1 & 1 & 3 \\ 0 & 1 & 0 \end{vmatrix} \cup \begin{vmatrix} 3 & 2 & 0 & 1 \\ 2 & 1 & 0 & 1 \\ 1 & 2 & 0 & -1 \\ 0 & 1 & 1 & 2 \end{vmatrix} \right\}.$$

We can sum the first columns of the bimatrix differently also still the value would remain the same.

***Example 1.5.4:*** Let

$$A_B = \begin{bmatrix} 9 & 2 & 3 \\ 18 & 0 & 1 \\ 21 & 0 & 0 \end{bmatrix} \cup \begin{bmatrix} 90 & 3 \\ 45 & 2 \end{bmatrix}$$

$$= \begin{vmatrix} 3 & + & 3 & + & 3 & 2 & 3 \\ 6 & + & 6 & + & 6 & 0 & 1 \\ 7 & + & 7 & + & 7 & 0 & 0 \end{vmatrix} \cup \begin{bmatrix} 30 & + & 30 & + & 30 & 3 \\ 15 & + & 15 & + & 15 & 2 \end{bmatrix}$$

$$= \left\{ \begin{vmatrix} 3 & 2 & 3 \\ 6 & 0 & 1 \\ 7 & 0 & 0 \end{vmatrix} + \begin{vmatrix} 3 & 2 & 3 \\ 6 & 0 & 1 \\ 7 & 0 & 0 \end{vmatrix} + \begin{vmatrix} 3 & 2 & 3 \\ 6 & 0 & 1 \\ 7 & 0 & 0 \end{vmatrix} \right\}$$

$$\cup \left\{ \begin{vmatrix} 30 & 3 \\ 15 & 2 \end{vmatrix} + \begin{vmatrix} 30 & 3 \\ 15 & 2 \end{vmatrix} + \begin{vmatrix} 30 & 3 \\ 15 & 2 \end{vmatrix} \right\}$$

$$= \left\{ \begin{vmatrix} 3 & 2 & 3 \\ 6 & 0 & 1 \\ 7 & 0 & 0 \end{vmatrix} \cup \begin{vmatrix} 30 & 3 \\ 15 & 2 \end{vmatrix} \right\} + \left\{ \begin{vmatrix} 3 & 2 & 3 \\ 6 & 0 & 1 \\ 7 & 0 & 0 \end{vmatrix} \cup \begin{vmatrix} 30 & 3 \\ 15 & 2 \end{vmatrix} \right\}$$

$$+ \left\{ \begin{vmatrix} 3 & 2 & 3 \\ 6 & 0 & 1 \\ 7 & 0 & 0 \end{vmatrix} \cup \begin{vmatrix} 30 & 3 \\ 15 & 2 \end{vmatrix} \right\}$$

$$= (14, 15) + (14, 15) + (14, 15)$$

$$= (42, 45)$$

$$|A_B| = (42, 45).$$

Suppose one wishes to divide the first columns of the bimatrix differently:

$$\begin{vmatrix} 4 & + & 4 & + & 1 & 2 & 3 \\ 9 & + & 9 & + & 0 & 0 & 1 \\ 10 & + & 10 & + & 1 & 0 & 0 \end{vmatrix} \cup \begin{vmatrix} 40 & + & 40 & + & 10 & 3 \\ 20 & + & 20 & + & 5 & 2 \end{vmatrix} =$$

$$\left\{ \begin{vmatrix} 4 & 2 & 3 \\ 9 & 0 & 1 \\ 10 & 0 & 0 \end{vmatrix} + \begin{vmatrix} 4 & 2 & 3 \\ 9 & 0 & 1 \\ 10 & 0 & 0 \end{vmatrix} + \begin{vmatrix} 1 & 2 & 3 \\ 0 & 0 & 1 \\ 1 & 0 & 0 \end{vmatrix} \right\} \cup$$

$$\left\{ \begin{vmatrix} 40 & 3 \\ 20 & 2 \end{vmatrix} + \begin{vmatrix} 40 & 3 \\ 20 & 2 \end{vmatrix} + \begin{vmatrix} 10 & 3 \\ 5 & 2 \end{vmatrix} \right\}$$

$$= \left\{ \begin{vmatrix} 4 & 2 & 3 \\ 9 & 0 & 1 \\ 10 & 0 & 0 \end{vmatrix} \cup \begin{vmatrix} 40 & 3 \\ 20 & 2 \end{vmatrix} \right\} + \left\{ \begin{vmatrix} 4 & 2 & 3 \\ 9 & 0 & 1 \\ 10 & 0 & 0 \end{vmatrix} \cup \begin{vmatrix} 40 & 3 \\ 20 & 2 \end{vmatrix} \right\}$$

$$+ \left\{ \begin{vmatrix} 1 & 2 & 3 \\ 0 & 0 & 1 \\ 1 & 0 & 0 \end{vmatrix} \cup \begin{vmatrix} 10 & 3 \\ 5 & 2 \end{vmatrix} \right\}$$

$$= \{(20, 20) + (20, 20) + (2, 5)\}$$

$$= (42, 45).$$

i.e., $|A_B| = (42, 45)$.

Thus the division of number of any column in the bimatrix does not affect the value of it.

It is important to note that in general if $A_B$ and $C_B$ are any two square bimatrices of same order we have

$$|A_B + C_B| \neq |A_B| + |C_B|.$$

***Example 1.5.5:*** Let

$$A_B = \begin{vmatrix} 3 & 1 & 2 \\ 0 & 1 & 5 \\ 8 & 9 & 4 \end{vmatrix} \cup \begin{vmatrix} 4 & 0 & 1 \\ 2 & 5 & 2 \\ 1 & 1 & 4 \end{vmatrix}$$

and

$$B_B = \begin{vmatrix} 0 & 2 & 0 \\ 1 & 0 & 1 \\ 4 & 0 & 0 \end{vmatrix} \cup \begin{vmatrix} 4 & 0 & 1 \\ 0 & 0 & 5 \\ 0 & -4 & 0 \end{vmatrix}$$

$A_B + B_B =$

$$\left\{\begin{vmatrix} 3 & 1 & 2 \\ 0 & 1 & 5 \\ 8 & 9 & 4 \end{vmatrix} \cup \begin{vmatrix} 4 & 0 & 1 \\ 2 & 5 & 2 \\ 1 & 1 & 4 \end{vmatrix}\right\} + \left\{\begin{vmatrix} 0 & 2 & 0 \\ 1 & 0 & 1 \\ 5 & 0 & 0 \end{vmatrix} \cup \begin{vmatrix} 4 & 0 & 1 \\ 0 & 0 & 5 \\ 0 & -4 & 0 \end{vmatrix}\right\}$$

$$= \left\{\begin{vmatrix} 3 & 1 & 2 \\ 0 & 1 & 5 \\ 8 & 9 & 4 \end{vmatrix} + \begin{vmatrix} 0 & 2 & 0 \\ 1 & 0 & 1 \\ 5 & 0 & 0 \end{vmatrix}\right\} \cup \left\{\begin{vmatrix} 4 & 0 & 1 \\ 2 & 5 & 2 \\ 1 & 1 & 4 \end{vmatrix} + \begin{vmatrix} 4 & 0 & 1 \\ 0 & 0 & 5 \\ 0 & -4 & 0 \end{vmatrix}\right\}$$

$$= \left\{\begin{vmatrix} 3 & 3 & 2 \\ 1 & 1 & 6 \\ 18 & 9 & 4 \end{vmatrix} + \begin{vmatrix} 8 & 0 & 2 \\ 2 & 5 & 7 \\ 1 & -3 & 4 \end{vmatrix}\right\}$$

= (64, 306).

Now consider $|A_B| \cup |C_B|$

= (–99, 69) + (10, 80)
= (–89, 149).
So $|A_B + C_B| \neq |A_B| + |C_B|$.

Hence our claim. Now in any matrix if two rows are identical or 2 columns are identical or the same we have then determinant value to be zero. But we see in case of bimatrix $A_B = A_1 \cup A_2$ when two rows or columns in $A_1$ or $A_2$ are the same we still may not have the determinant of the bimatrix to be zero.

***Example 1.5.6:*** Let

$$A_B = \begin{vmatrix} 3 & 1 & 3 & 4 \\ 0 & 0 & 2 & 1 \\ 1 & 1 & 4 & 2 \\ 5 & 0 & 1 & 0 \end{vmatrix} \cup \begin{vmatrix} 0 & 0 & 2 & 1 \\ 3 & -1 & 3 & -4 \\ 4 & 4 & 0 & 0 \\ 5 & 0 & 1 & 0 \end{vmatrix}$$

be a bimatrix. This has two rows in $A_1$ and $A_2$ to be the same yet $|A_B| \neq (0, 0)$.

## 1.6 Laplace equation for bimatrix

In this section we shall consider yet another, more general technique of expanding a bideterminant using a method analogous to Laplace expansion, which includes as a special case, the expansion by bicofactors.

Now the bideterminant of $|A_B|$ is written as the sum of terms each of which is the product of two bideterminants. We begin by considering the first m rows $(m < n)$ and columns of $A_B = A_1 \cup A_2$ (i.e., m rows of $A_1$ and $A_2$ are taken and both $A_1$ and $A_2$ are $n \times n$ matrices so $A_B$ is a $n \times n$ square bimatrix). If we collect the terms from $A_1$ and $A_2$ containing $a^1_{11}\, a^1_{22} \ldots a^1_{mm}$ and $a^2_{11}\, a^2_{22} \ldots a^2_{mm}$ in the expansion of $|A_B|$ and factor out this quantity we are left with

$$\Sigma \pm a^1_{m+1i}\, a^1_{m+2j} \ldots a^1_{nr} + \sum \pm a^2_{m+1i} + a^2_{m+2j} \ldots a^2_{nr} \qquad \text{(a)}$$

the sum being taken over all permutations of the second scripts. This sum however represents the determinant of the subbimatrix formed from $A_B$ by crossing out the first m rows and columns of $A_1$ and $A_2$. Note that the above sum (a) will be obtained also if we collect all terms containing $a^1_{1\upsilon}\, a^1_{1\nu} \ldots a^1_{mw}$ and $a^2_{1\upsilon}\, a^2_{1\nu} \ldots a^2_{mw}$ where $[(\upsilon,\ \nu,\ \ldots,\ w)$ represents a permutation of (1, 2, 3,…, m)] and factor out $a^1_{1\upsilon}\ a^1_{1\nu} \ldots a^1_{mw}$ and $a^2_{1\upsilon}\ a^2_{1\nu} \ldots a^2_{mw}$ related to matrices $A^1$ and $A^2$ (where $A^1 = (a^1_{ij})$ and $A^2 = (a^2_{ij})$). The sign will alternate depending on whether the permutation (υ,ν, …, w) is even or odd. Thus the following terms appear in the expansion of $|A_B|$,

$$\left(\sum(\pm)\, a^1_{1\upsilon}\ a^1_{2\nu} \ldots a^1_{nw}\right)\left(\sum(\pm)\, a^1_{1m+1i}\ a^1_{m+2i} \ldots a^1_{nr}\right)$$

$$\cup\ \left(\sum(\pm)\, a^1_{1\upsilon}\ a^1_{2\nu} \ldots a^1_{mw}\right)\left(\sum(\pm)\, a^2_{m+1i}\ a^2_{m+2j} \ldots a^2_{nr}\right) \qquad \text{(b)}$$

where $\sum(\pm)\, a^1_{1\upsilon}\ a^1_{2\nu} \cdots a^1_{mw}$ and $\sum(\pm)\, a^2_{1\upsilon}\ a^2_{2\nu} \cdots a^2_{mw}$ is the bideterminant of the subbimatrix formed from the first m rows and n columns of $A_1$ and $A_2$ respectively where $A_B = A_1 \cup A_2$.

Thus the above expression (b) is the product of the bideterminant of the subbimatrix formed by crossing out the first m rows and columns. We have the correct sign since in the expression of $|A_B|$ the terms $a^1_{11}\, a^1_{22} \cdots a^1_{nn}$. and $a^2_{11}\, a^2_{22} \cdots a^2_{nn}$ has a plus sign.

Next we shall consider the m × m subbimatrix formed from rows $\left(i^1_1, i^1_2, \ldots, i^1_m\right)$ and $\left(i^2_1, i^2_2, \ldots, i^2_m\right)$ and columns $j^1_1, j^1_2, \ldots, j^1_m$ and $j^2_1, j^2_2, \cdots, j^2_m$. Except for the sign the expansion of $|A_B|$ will contain the product of the bideterminant of this subbimatrix and the determinant of the subbimatrix formed by crossing out rows $\left(i^1_1, i^1_2, \cdots, i^1_m\right)$ and $\left(i^2_1, i^2_2, \cdots, i^2_m\right)$ and columns $j^1_1, j^1_2, \ldots, j^1_m$ and $j^2_1, j^2_2, \cdots, j^2_m$. The sign of the product is determined by the method used in the expansion of bicofactors. Thus the m × m subbimatrix is moved so that it occupies the first m rows and m columns in both $A_1$ and $A_2$.

Let us assume quite logically

$$i^1_1 < i^1_2 \cdots < i^1_m,\ i^2_1 < i^2_2 \cdots < i^2_m \text{ and}$$

$$j_1^1 < j_2^1 \cdots < j_m^1 \text{ and}$$
$$j_1^2 < j_2^2 \cdots < j_m^2 .$$

Then after

$$\left(i_1^1 - 1\right) + \left(i_2^1 - 2\right) + \cdots + \left(i_m^1 - m\right) \text{ and}$$
$$\left(i_1^2 - 1\right) + \left(i_2^2 - 2\right) + \cdots + \left(i_m^2 - m\right)$$

interchange of rows of $A_1$ and $A_2$ (where $A_B = A_1 \cup A_2$) and

$$\left(j_1^1 - 1\right) + \left(j_2^1 - 2\right) + \cdots + \left(j_m^1 - m\right) \text{ and}$$
$$\left(j_1^2 - 1\right) + \left(j_2^2 - 2\right) + \cdots + \left(j_m^2 - m\right)$$

interchanges of columns the $m \times m$ subbimatrix formed from

rows $\left(i_1^1, i_2^1, \cdots, i_m^1\right)$ and $\left(i_1^2, i_2^2, \cdots, i_m^2\right)$ and

columns $j_1^1, j_2^1, ..., j_m^1$ and $j_1^2, j_2^2, ..., j_m^2$

lies in rows 1, …, m and columns 1, …, m (for both $A_1$ and $A_2$). Further the order of the remaining columns has not changed in $A_B = A_1 \cup A_2$. Once this requirement is completed we are back at the case already considered.

The sign depends on

$$\sum_{K=1}^{m}\left(i_k^1 + j_k^1\right) - 2\,(1 + 2 + ... + m) \text{ and}$$
$$\sum_{K=1}^{m}\left(i_k^2 + j_k^2\right) - 2\,(1 + 2 + ... + m) .$$

However 1 + 2 + … + m = ½ m (m + 1).

Since 2m (m + 1) is always even the sign attached to the product of the determinants is

$$(-1)^{\sum_{K=1}^{m}\left(i_k^1+j_k^1\right)} \text{ and } (-1)^{\sum_{K=1}^{m}\left(i_k^2+j_k^2\right)}.$$

Now we proceed on to define the notion of complementary biminor for the square n × n bimatrix. $A_B^{n \times n} = A_1 \cup A_2$. Given the $n^{th}$ order bimatrix $A_B$. The bideterminant of $(n – m)^{th}$ order subbimatrix $P_B$ formed ($P_B = P_1 \cup P_2$, $A_B = A_1 \cup A_2$) is the $(n – m)^{th}$ order subbimatrix of $A_B$ by crossing out

rows $i_1^1,...,i_m^1;\ i_1^2, i_2^2,..., i_m^2$ and
columns $j_1^1, j_2^1,..., j_m^1, j_1^2, j_2^2,..., j_m^2$

is called the complementary biminor of the $m^{th}$ order subbimatrix $N_B$ formed from the
rows $i_1^1, i_2^1,...,i_m^1;\ i_1^2, i_2^2,..., i_m^2$
and
columns $j_1^1, j_2^1,..., j_m^1, j_1^2, j_2^2,..., j_m^2$.

Now we proceed on to define the notion of complementary bicofactor. However $N_B$ and $P_B$ would be as in the above definition.

**DEFINITION 1.6.1:** *With $N_B$ and $P_B$ as defined above the bideterminant*

$$|M_B| = (-1)^{\sum_{i=1}^{m}\left(i_k^t+j_k^t\right)} |P_B|,$$

*is called the complementary bicofactor of $N_B$ in $A_B$ i.e.,*

$$|M_B| = \left( (-1)^{\sum_{i=1}^{m}\left(i_k^1+j_k^1\right)} |P_1|, \sum (-1)^{\sum_{i=1}^{m}\left(i_k^2+j_k^2\right)} |P_2| \right).$$

We illustrate this by the following example:

***Example 1.6.1:*** Consider $A_B = A_1 \cup A_2$
where

$$|A_B| = \begin{vmatrix} a^1_{11} & a^1_{12} & \cdots & a^1_{15} \\ a^1_{21} & a^1_{22} & \cdots & a^1_{25} \\ \vdots & & & \\ a^1_{51} & a^1_{52} & \cdots & a^1_{55} \end{vmatrix} \cup \begin{vmatrix} a^2_{11} & a^2_{12} & \cdots & a^2_{15} \\ a^2_{21} & a^2_{22} & \cdots & a^2_{25} \\ \vdots & & & \\ a^2_{51} & a^2_{52} & \cdots & a^2_{55} \end{vmatrix}.$$

The bideterminant of the subbimatrix $N_B$ formed from columns 2 and 5 and rows 1 and 3 is

$$|N_B| = \begin{vmatrix} a^1_{12} & a^1_{15} \\ a^1_{32} & a^1_{35} \end{vmatrix} \cup \begin{vmatrix} a^2_{12} & a^2_{15} \\ a^2_{32} & a^2_{35} \end{vmatrix}.$$

The complementary biminor is

$$|P_B| = \begin{vmatrix} a^1_{21} & a^1_{23} & a^1_{24} \\ a^1_{41} & a^1_{43} & a^1_{44} \\ a^1_{51} & a^1_{53} & a^1_{54} \end{vmatrix} \cup \begin{vmatrix} a^2_{21} & a^2_{23} & a^2_{24} \\ a^2_{41} & a^2_{43} & a^2_{44} \\ a^2_{51} & a^2_{53} & a^2_{54} \end{vmatrix}.$$

Hence
$\left(\sum i^1_k + j^1_k,\ \sum i^2_k + j^2_k\right)$ = (2 + 5+1+3 = 1, 2 + 5 + 1 + 3 = 11) and the complementary bicofactor is $|M_B| = -|P_B|$.

These method lead us to derive a new method of expanding the bideterminant of $|A_B|$. Select any m rows of $A_B = A_1 \cup A_2$. From these m rows we can form (n! / m! (n – m)!, n!/m! (n – m)!) different $m^{th}$ order subbimatrices where (n!/m! (n – m)!, n!/m! (n – m)! ) defines the number of combinations of n columns of the pair of matrices $A_1$ and $A_2$ of $A_B$ taken m at a time. In choosing these submatrices $N_B$ we always keep the order of the columns in $A_B$ unchanged. For each subbimatrices $N_B$ we find $|N_B|$ and the corresponding complementary bicofactor $|M_B|$ then we form the product $|N_B|\ M_B|$. From $|N_B|$ we obtain a pair of m! terms and from $|M_B|$ a pair of (n – m) ! terms (each one for $A_1$ and $A_2$). Hence each product $|N_B||M_B|$ yields with the correct sign the pair m! (n – m)! terms of the expansion of $|A_B| = |A_1| \cup |A_2|$. In all these there are n! / m! (n – m)! product $|N_1|\ |M_1|$, $|N_2|\ |M_2|$ (where $N_B = N_1 \cup N_2$ and $M_B = M_1 \cup M_2$) which yield a total of

$$\left(\frac{n!m!(n-m)!}{m!(n-m)!},\ \frac{n!m!(n-m)!}{m!(n-m)!}\right) = (n!, n!)$$

terms in the expansion of $|A_B|$.

Thus we have obtained all terms in $|A_B|$ since our method of selecting the m$^{th}$ order subbimatrices from the m rows eliminates any possible repetition of terms.

Hence we can write $\left|A_B\right| = \left|A_1\right| \cup \left|A_2\right|$ = $(d_1, d_2)$ =

$$\left(\sum_{j_1^1 < j_2^1 < \cdots < j_m^1} \left|N_1\left(i_1^1,\ldots,i_m^1 \left| j_1^1,\ldots,j_m^1\right.\right)\right| \left|M_1\right|,\right.$$
$$\left.\sum_{j_1^2 < j_2^2 < \cdots < j_m^2} \left|N_2\left(i_1^2,\ldots,i_m^2 \left| j_1^2,\ldots,j_m^2\right.\right)\right| \left|M_2\right|\right)$$

where

$$\left|N_1\left(i_1^1, i_2^1,\ldots,i_m^1 \left| j_1^1,\ldots,j_m^1\right.\right)\right| \text{ and } \left|N_2\left(i_1^2, i_2^2,\ldots,i_m^2 \left| j_1^2,\ldots,j_m^2\right.\right)\right|$$

are the m$^{th}$ order bideterminant of the subbimatrix formed from

| | |
|---|---|
| rows | $i_1^1, i_2^1,\ldots,i_m^1,\ i_1^2, i_2^2,\ldots,i_m^2$ and |
| columns | $j_1^1, j_2^1,\ldots,j_m^1,\ j_1^2, j_2^2,\ldots,j_m^2$ . |

The sum is taken over the n! / m! (n – m)! choices for $j_1^1, j_2^1,\ldots, j_m^1,\ j_1^2, j_2^2,\ldots, j_m^2$. The notation $j_1^r < j_2^r < \ldots < j_m^r$ (r = 1, 2) indicates that the sum is taken over all choices of the columns such that the column order is maintained in both

$A_1$ and $A_2$. We call this analogue, technique of expansion as BiLaplace method:

We select any m rows from ($A_B = A_1 \cup A_2$) $A_1$ and $A_2$ (note they need not be adjacent). From these m rows we form the n! / m! (n – m)! possible $m^{th}$ order bi determinants and find their individual complementary bicofactors. We then multiply the bideterminant by its complementary bicofactor and add the n! / m! (n – m)! terms to obtain

$$|A_B| = (|A_1| \cup |A_2|).$$

To find $|A_B|$ we can also expand by any m-columns in an analogous way.

$$|A_B| = \Bigg( \sum_{i_1^1 < i_2^1 < \ldots < i_m^1} \left| N_1 \left( i_1^1, i_2^1, \ldots, i_n^1 \,\middle|\, j_1^1, j_2^1, \ldots, j_m^1 \right) \right| \times |M_1| \cup$$

$$\sum_{i_1^2 < i_2^2 < \ldots < i_m^2} \left| N_2 \left( i_1^2, i_2^2, \ldots, i_n^2 \,\middle|\, j_1^2, j_2^2, \ldots, j_m^2 \right) \right| \, |M_2| \Bigg).$$

Now we illustrate the expansion using a bimatrix of order 4.

***Example 1.6.2:*** Let

$$A_B = \begin{vmatrix} a_{11}^1 & a_{12}^1 & a_{13}^1 & a_{14}^1 \\ a_{21}^1 & a_{22}^1 & a_{23}^1 & a_{24}^1 \\ a_{31}^1 & a_{32}^1 & a_{33}^1 & a_{34}^1 \\ a_{41}^1 & a_{42}^1 & a_{43}^1 & a_{44}^1 \end{vmatrix} \cup \begin{vmatrix} a_{11}^2 & a_{12}^2 & a_{13}^2 & a_{14}^2 \\ a_{21}^2 & a_{22}^2 & a_{23}^2 & a_{24}^2 \\ a_{31}^2 & a_{32}^2 & a_{33}^2 & a_{34}^2 \\ a_{41}^2 & a_{42}^2 & a_{43}^2 & a_{44}^2 \end{vmatrix} = A_1 \cup A_2$$

Let us expand this bimatrix by the first and last rows. There will be 4! / 2! 2! = 6 terms.

The bideterminants of order 2 which can be formed from rows 1 and 4 are

$$|N_B^1| = \begin{vmatrix} a_{11}^1 & a_{12}^1 \\ a_{41}^1 & a_{42}^1 \end{vmatrix} \cup \begin{vmatrix} a_{11}^2 & a_{12}^2 \\ a_{41}^2 & a_{42}^2 \end{vmatrix}$$

$$\left|N_B^2\right| = \begin{vmatrix} a_{11}^1 & a_{13}^1 \\ a_{41}^1 & a_{43}^1 \end{vmatrix} \cup \begin{vmatrix} a_{11}^2 & a_{13}^2 \\ a_{41}^2 & a_{43}^2 \end{vmatrix}$$

$$\left|N_B^3\right| = \begin{vmatrix} a_{11}^1 & a_{14}^1 \\ a_{41}^1 & a_{44}^1 \end{vmatrix} \cup \begin{vmatrix} a_{11}^2 & a_{14}^2 \\ a_{41}^2 & a_{44}^2 \end{vmatrix}$$

$$\left|N_B^4\right| = \begin{vmatrix} a_{12}^1 & a_{13}^1 \\ a_{42}^1 & a_{43}^1 \end{vmatrix} \cup \begin{vmatrix} a_{12}^2 & a_{13}^2 \\ a_{42}^2 & a_{43}^2 \end{vmatrix}$$

$$\left|N_B^5\right| = \begin{vmatrix} a_{12}^1 & a_{14}^1 \\ a_{42}^1 & a_{44}^1 \end{vmatrix} \cup \begin{vmatrix} a_{12}^2 & a_{14}^2 \\ a_{42}^2 & a_{44}^2 \end{vmatrix}$$

and

$$\left|N_B^6\right| = \begin{vmatrix} a_{13}^1 & a_{14}^1 \\ a_{42}^1 & a_{44}^1 \end{vmatrix} \cup \begin{vmatrix} a_{13}^2 & a_{14}^2 \\ a_{43}^2 & a_{44}^2 \end{vmatrix}.$$

The corresponding complementary bicofactors are

$$\left|M_B^1\right| = \begin{vmatrix} a_{23}^1 & a_{24}^1 \\ a_{33}^1 & a_{34}^1 \end{vmatrix} \cup \begin{vmatrix} a_{23}^2 & a_{24}^2 \\ a_{33}^2 & a_{34}^2 \end{vmatrix}$$

$$\left|M_B^2\right| = \begin{vmatrix} a_{22}^1 & a_{24}^1 \\ a_{32}^1 & a_{34}^1 \end{vmatrix} \cup \begin{vmatrix} a_{22}^2 & a_{24}^2 \\ a_{32}^2 & a_{34}^2 \end{vmatrix}$$

$$\left|M_B^3\right| = \begin{vmatrix} a_{22}^1 & a_{23}^1 \\ a_{32}^1 & a_{33}^1 \end{vmatrix} \cup \begin{vmatrix} a_{22}^2 & a_{23}^2 \\ a_{32}^2 & a_{33}^2 \end{vmatrix}$$

$$\left|M_B^4\right| = \begin{vmatrix} a_{21}^1 & a_{24}^1 \\ a_{31}^1 & a_{34}^1 \end{vmatrix} \cup \begin{vmatrix} a_{21}^2 & a_{24}^2 \\ a_{31}^2 & a_{34}^2 \end{vmatrix}$$

$$\left|M_B^5\right| = \begin{vmatrix} a_{21}^1 & a_{23}^1 \\ a_{31}^1 & a_{33}^1 \end{vmatrix} \cup \begin{vmatrix} a_{21}^2 & a_{24}^2 \\ a_{31}^2 & a_{33}^2 \end{vmatrix}$$

and

$$\left|M_B^6\right| = \begin{vmatrix} a_{21}^1 & a_{22}^1 \\ a_{31}^1 & a_{32}^1 \end{vmatrix} \cup \begin{vmatrix} a_{21}^2 & a_{22}^2 \\ a_{31}^2 & a_{32}^2 \end{vmatrix}$$

$$\left|A_B\right| = \sum_{K=1}^{6} \left|N_B^K\right| \; \left|M_B^K\right| .$$

If we denote by $N_B^K = N^{K1} \cup N^{K2}$ and $M_B^K = M^{K1} \cup M^{K2}$, K = 1, 2, 3, …, 6 ;
then

$$\left|A_B\right| = \sum_{K=1}^{6} \left|N^{K1}\right| \; \left|M^{K1}\right| \cup \sum_{K=1}^{6} \left|N^{K2}\right| \; \left|M^{K2}\right| .$$

Now we proceed on to define multiplication of bideterminants. This is done analogous to multiplication of determinants i.e., the simple multiplication relation for determinant of the product of square matrices.
If $A_B$ and $B_B$ be square bimatrices of order n then their product $C_B = A_B B_B$
Here

$A_B = A_1 \cup A_2$
$B_B = B_1 \cup B_2$

so

$C_B = A_B B_B = A_1 B_1 \cup A_2 B_2$

$$\left|C_B\right| = \left|A_B\; B_B\right| = \left|A_1\right|\left|B_1\right| \;\cup\; \left|A_2\right| \;\left|B_2\right| .$$

i.e., the determinant of the product is the product of the determinants.

Now suppose $D_B$ is a $2n \times 2n$ bimatrix which is partitioned as a bimatrix $D_B = D_1 \cup D_2$

$$= \begin{vmatrix} A_1 & O \\ -I_n & B_1 \end{vmatrix} \cup \begin{vmatrix} A_2 & O \\ -I_n & B_2 \end{vmatrix} .$$

Applying the biLaplace expansion by the last n rows to the bideterminant $|D_B|$ we obtain

$$|D_B| = \begin{vmatrix} A_1 & O \\ -I_m & B_1 \end{vmatrix} \cup \begin{vmatrix} A_2 & O \\ -I_n & B_2 \end{vmatrix}$$

$$= \quad (-1)^{2n^2+n(n+1)} |A_1||B_1| \cup (-1)^{2n^2+n(n+1)} |A_2|\ |B_2| .$$

Since the complementary biminor of any subbimatrix including one of the first n rows will have a column of zeros. The bideterminant is completely independent of the bimatrix appearing in the lower left of the above equation.

Bimatrix – $I_n$ ∪ – $I_n$ was placed there for the following special reasons. For if we consider $B_B$ = $b_j^1 \cup b_j^2$,
then consider

$$\begin{bmatrix} A_1 \\ -I_n \end{bmatrix} b_1 \cup \begin{bmatrix} A_2 \\ -I_n \end{bmatrix} b_2 = |A_1\ b_1| \cup |A_2\ b_2| \qquad (*)$$

So the equation (*) is a linear combinations of the first n columns of $D_B$. We can add (*) to n + 1 columns of $D_B$ without changing the value of the bideterminant which yields

$$\begin{bmatrix} A_1 b_1 \\ 0 \end{bmatrix} \cup \begin{bmatrix} A_2 b_2 \\ 0 \end{bmatrix}$$

as the new n + 1 columns. Continuing this process and adding

$$\begin{bmatrix} A_1 \\ -I_n \end{bmatrix} b_j^1 \cup \begin{bmatrix} A_2 \\ -I_n \end{bmatrix} b_j^2 = \begin{bmatrix} A_1 b_j^1 \\ -b_j^1 \end{bmatrix} \cup \begin{bmatrix} A_2 b_j^2 \\ -b_j^2 \end{bmatrix}$$

to the n + j$^{th}$ column of $D_B$ we get

$$|D_B| = \begin{vmatrix} A_1 & A_1B_1 \\ -I_n & O \end{vmatrix} \cup \begin{vmatrix} A_2 & A_2B_2 \\ -I_n & O \end{vmatrix}$$

We can expand by the last n rows using the biLaplace expansion and obtain

$$\begin{aligned} |D_B| &= (-1)^{n^2+n(n+1)} |-I_n| |A_1B_1| \cup (-1)^{n^2+n(n+1)} |-I_n| |A_2B_2| \\ &= (-1)^{2n(n+1)} |A_1B_1| \cup (-1)^{2n(n+1)} |A_2B_2| \\ &= |A_1B_1| \cup |A_2B_2|. \end{aligned}$$

Thus $|D_B| = |A_1||B_1| \cup |A_2||B_2| = |A_1B_1| \cup |A_2B_2|$.

We illustrate this by the following simple example.

***Example 1.6.3:*** Let

$$A_B = \begin{vmatrix} 2 & 3 \\ 1 & 4 \end{vmatrix} \cup \begin{vmatrix} 3 & 6 \\ 1 & 1 \end{vmatrix}$$

and

$$B_B = \begin{vmatrix} 1 & 6 \\ 3 & 2 \end{vmatrix} \cup \begin{vmatrix} 5 & 2 \\ 1 & 3 \end{vmatrix}.$$

$$C_B = |A_B\ B_B| = \begin{bmatrix} 11 & 18 \\ 13 & 14 \end{bmatrix} \cup \begin{bmatrix} 21 & 24 \\ 6 & 5 \end{bmatrix}$$

$|A_1| = 5$
$|A_2| = -3$
$|B_1| = -16$
$|B_2| = 13$

$$\begin{aligned} |C_B| &= (-80, -39) \\ |C_B| &= |A_B|\ |B_B| \\ &= (5, -3)\ (-16, 13) \\ &= (-80, -39). \end{aligned}$$

So $|C_B| = |A_B B_B| = |A_B| . |B_B|$.

Now we shall define the bideterminant of the product of rectangular bimatrices.

Let $A_B$ be an $m \times n$ bimatrix and $B_B$ an $n \times m$ bimatrix with $m < n$. If $C_B = A_B B_B$, $|C_B|$ is an $m^{th}$ order bideterminant. We shall show how to express $|C_B|$ in terms of the bideterminants of order m which can be found from $A_B$ and $B_B$. Consider the product.

$$\begin{vmatrix} I_n & A_1 \\ O & I_n \end{vmatrix} \begin{vmatrix} A_1 & O \\ -I_n & B_1 \end{vmatrix} \cup \begin{vmatrix} I_n & A_2 \\ O & I_n \end{vmatrix} \begin{vmatrix} A_2 & O \\ -I_n & B_2 \end{vmatrix}$$

$$= \begin{vmatrix} O & A_1B_1 \\ -I_n & B_1 \end{vmatrix} \cup \begin{vmatrix} O & A_2B_2 \\ -I_n & B_2 \end{vmatrix}.$$

BiLaplace expansion indicates immediately that the first bideterminant has the value unity and hence

$$\begin{vmatrix} A_1 & O \\ -I_n & B_1 \end{vmatrix} \cup \begin{vmatrix} A_2 & O \\ -I_n & B_2 \end{vmatrix}$$

$$= (-1)^{n(m+1)} \left| A_1 B_1 \right| \cup (-1)^{n(m+1)} \left| A_2 B_2 \right|.$$

To evaluate the bideterminant of order n + m on the left we shall use biLaplace expansion by the first m rows. First note that non vanishing bideterminants of order m can be formed only from the columns of $A_B$, since the use of any other columns would introduce a column of zeros and hence a vanishing bideterminant.

Thus there are no more than $n! /m! (n - m)!$, $n! /m! (n - 1)!$ non-vanishing terms occurring in pairs in the expansion. The complementary biminor to any bideterminant $A_B$ of order m formed from $A_B$ will have $n - m$, $n - m$ columns from $-I_n$. The complementary biminor of order n will be of the form.

$$\left|-e^1_{u_1} - e^1_{u_2}, ..., -e_{u_{n-m}}, B_1\right| \cup \left|-e^2_{u_1} - e^2_{u_2}, ..., -e_{u_{n-m}}, B_2\right|$$

where $u^1_1, u^1_2, ..., u^1_{n-m}$; $u^2_1, u^2_2, ..., u^2_{n-m}$ refer to it columns of $A_1$ and $A_2$ not in $A_1$ and $A_2$ respectively ($A_B = A_1 \cup A_2$, $\Delta_B = \Delta_1 \cup \Delta_2$).

We can immediately expand by bicofactors proceeding from the first column to the second etc, … to the $(n - m)^{th}$ column. Note that aside from sign this expansion crosses out rows $u^1_1, u^1_2, ..., u^1_{n-m}$ of $B_1$ and $u^2_2, u^2_2, ..., u^2_{n-m}$ of $B_2$ so that in the end we obtain a bideterminant of order m formed from $B_1$ and $B_2$ which contains the same rows as the corresponding columns chosen from $A_1$ and $A_2$ to be in $\Delta_1$ and $\Delta_2$ respectively.

The similar form of separately working will yield the result, which is left as an exercise to the reader.

Now we represent this by an example.

***Example: 1.6.4:*** Let

$$A_B = \begin{bmatrix} 1 & 4 & 5 \\ 2 & 0 & 3 \end{bmatrix} \cup \begin{bmatrix} 0 & 1 & 1 \\ 2 & 0 & 1 \end{bmatrix}$$

and

$$B_B = \begin{bmatrix} 3 & 0 \\ 9 & 2 \\ 1 & 7 \end{bmatrix} \cup \begin{bmatrix} 1 & 1 \\ 0 & 2 \\ 5 & -1 \end{bmatrix}$$

$$A_B\ B_B = \begin{bmatrix} 1 & 4 & 5 \\ 2 & 0 & 3 \end{bmatrix} \begin{bmatrix} 3 & 0 \\ 9 & 2 \\ 1 & 7 \end{bmatrix} \cup \begin{bmatrix} 0 & 1 & 1 \\ 2 & 0 & 1 \end{bmatrix} \begin{bmatrix} 1 & 1 \\ 0 & 2 \\ 5 & -1 \end{bmatrix}$$

$$= \begin{bmatrix} 44 & 43 \\ 9 & 21 \end{bmatrix} \cup \begin{bmatrix} 5 & 1 \\ 7 & 1 \end{bmatrix}$$

$$\left|A_B\ B_B\right| = (537,\ -2)\ .$$

However $|A_B\ B_B|$ can be expressed as the sum of three terms each of which is the product of a 2 × 2 bideterminant from $A_B$ and 2 × 2 determinant from $B_B$. The three pairs of bideterminant which can be formed from $A_B$ $B_B$ respectively are: from columns 1, 2 of $A_B$ and rows of 1,2 of $B_B$.

$$|A_B^1| = \begin{vmatrix} 1 & 4 \\ 2 & 0 \end{vmatrix} \cup \begin{vmatrix} 0 & 1 \\ 2 & 0 \end{vmatrix} = (-8, -2)$$

$$|B_B^1| = \begin{vmatrix} 3 & 0 \\ 9 & 2 \end{vmatrix} \cup \begin{vmatrix} 1 & 1 \\ 0 & 2 \end{vmatrix} = (6, 2)$$

$$|A_B^2| = \begin{vmatrix} 1 & 5 \\ 2 & 3 \end{vmatrix} \cup \begin{vmatrix} 0 & 1 \\ 2 & 1 \end{vmatrix} = (-7, -2)$$

$$|B_B^2| = \begin{vmatrix} 3 & 0 \\ 1 & 7 \end{vmatrix} \cup \begin{vmatrix} 1 & 1 \\ 5 & -1 \end{vmatrix} = (21, -6)$$

$$|A_B^3| = \begin{vmatrix} 4 & 5 \\ 0 & 3 \end{vmatrix} \cup \begin{vmatrix} 1 & 1 \\ 0 & 1 \end{vmatrix} = (12, 1)$$

$$|B_B^3| = \begin{vmatrix} 9 & 2 \\ 1 & 7 \end{vmatrix} \cup \begin{vmatrix} 0 & 2 \\ 5 & -1 \end{vmatrix} = (61, -10).$$

$$|A_B\ B_B| = |\ A_B^1\ B_B^1\ | + |\ A_B^2\ B_B^2\ | + |\ A_B^3\ B_B^3\ |$$

= (–8 (6), –2 (2))
+ (–7 (21), –2 (–6))
+ (12 (61) 1 (–10))
= (– 48 – 147 + 732, – 4 + 12 – 10)
= (537, –2).

The results are indeed identical!

Now we proceed on to define the notion of bimatrix inverse. We also know given a number $a \neq 0$ there exists a number $a^{-1}$ such that $aa^{-1} = 1$.

Now if $A_B = A_1 \cup A_2$ be the given bimatrix does there exist a bimatrix $B_B = B_1 \cup B_2$ such that $A_B B_B = B_B A_B = I_B$? If such a bimatrix $B_B$ exists it is called as the biinverse of $A_B$. The biinverse of $A_B$ is usually written as $A_B^{-1} = A_1^{-1} \cup A_2^{-1}$ where $A_B = A_1 \cup A_2$.

It is most important to note that

$$A_B^{-1} \neq \frac{1}{A_B} \text{ or } \frac{I}{A_B}.$$

Note that for the given bimatrix $A_B = A_1 \cup A_2$ the biinverse of $A_B$ denote by $A_B^{-1} = A_1^{-1} \cup A_2^{-1}$ i.e., $A_B^{-1} = A_1^{-1} \cup A_2^{-1}$ is merely a symbol given to the matrix $B_B$ such that $A_B B_B = B_B A_B = I_B$.

It is very important to note that $A_B$ the bimatrix will have an biinverse only if $A_B$ is a square bimatrix. Hence the biinverse of $A_B$ will be only a square bimatrix.

Thus we define the bimatrix biinverse as; given a square bimatrix $A_B = A_1 \cup A_2$, if there exists a square bimatrix $A_B^{-1} = A_1^{-1} \cup A_2^{-1}$ which satisfies the identity

$$A_B \; A_B^{-1} = A_B^{-1} A_B = A_1 A_1^{-1} \cup A_2 A_2^{-1} = I \cup I$$

then $A_B^{-1}$ is called the biinverse or bireciprocal of $A_B$.

It is important to note that even if $A_B$ is a square mixed bimatrix then also $A_B^{-1}$ exists by $I_1 \cup I_2 = I_B$ will be such that $I_1 \neq I_2$. Just biinverse are introduced for the first time we illustrate it by examples.

***Example 1.6.5:*** Let

$$A_B = \begin{bmatrix} 1 & 0 \\ 2 & 3 \end{bmatrix} \cup \begin{bmatrix} 0 & 2 \\ -1 & 1 \end{bmatrix}$$

be the square bimatrix.

$$A_B^{-1} = \begin{bmatrix} 1 & 0 \\ -2/3 & 1/3 \end{bmatrix} \cup \begin{bmatrix} 1/2 & -1 \\ 1/2 & 0 \end{bmatrix}$$

$$A_B A_B^{-1} = \begin{bmatrix} 1 & 0 \\ 0 & 1 \end{bmatrix} \cup \begin{bmatrix} 1 & 0 \\ 0 & 1 \end{bmatrix}.$$

Now we consider an example of a mixed square bimatrix.

***Example 1.6.6:*** Given

$$A_B = \begin{bmatrix} 3 & 1 \\ 7 & 5 \end{bmatrix} \cup \begin{bmatrix} 1 & 2 & 2 \\ 2 & 1 & 2 \\ 2 & 2 & 1 \end{bmatrix}$$

this is a mixed square matrix. The biinverse of $A_B$ denoted by

$$A_B^{-1} = \begin{bmatrix} 5/8 & -1/8 \\ -7/8 & 3/8 \end{bmatrix} \cup \begin{bmatrix} -3/5 & 2/5 & 2/5 \\ 2/5 & -3/5 & 2/5 \\ 2/5 & 2/5 & -3/5 \end{bmatrix}.$$

Some of the properties of inverses which hold good in case of matrices also hold true is case of bimatrices. For it can be verified that

$$(A_B\ B_B)^{-1} = B_B^{-1}\ A_B^{-1}\ ;\ B_B^{-1}\ A_B^{-1}\ A_B\ \ B_B = B_B^{-1}\ B_B = I_B\,.$$

Just we will show by an example that

$$(A_B\ B_B)^{-1} = B_B^{-1}\ A_B^{-1}\,.$$

***Example 1.6.7:*** Let

$$A_B = \begin{bmatrix} 1 & 0 \\ 2 & 3 \end{bmatrix} \cup \begin{bmatrix} 2 & 1 \\ 5 & 3 \end{bmatrix}$$

$$B_B = \begin{bmatrix} 2 & 5 \\ 2 & 1 \end{bmatrix} \cup \begin{bmatrix} 0 & 5 \\ 6 & 4 \end{bmatrix}$$

be two 2 × 2 bimatrices

$$A_B^{-1} = \begin{bmatrix} 1 & 0 \\ -2/3 & 1/3 \end{bmatrix} \cup \begin{bmatrix} 3 & -1 \\ -5 & 2 \end{bmatrix}$$

$$B_B^{-1} = \begin{bmatrix} -1/8 & 5/8 \\ 1/4 & -1/4 \end{bmatrix} \cup \begin{bmatrix} -37/30 & 7/15 \\ 3/5 & -1/5 \end{bmatrix}$$

$$(A_B\, B_B)^{-1} = \left\{\begin{bmatrix} 1 & 0 \\ 2 & 3 \end{bmatrix}\begin{bmatrix} 2 & 5 \\ 2 & 1 \end{bmatrix}\right\}^{-1} \cup \left\{\begin{bmatrix} 2 & 1 \\ 5 & 3 \end{bmatrix}\begin{bmatrix} 0 & 5 \\ 6 & 4 \end{bmatrix}\right\}^{-1}$$

$$= \left\{\begin{bmatrix} 2 & 5 \\ 10 & 13 \end{bmatrix}\right\}^{-1} \cup \left\{\begin{bmatrix} 6 & 14 \\ 18 & 37 \end{bmatrix}^{-1}\right\}$$

$$= \begin{bmatrix} -13/24 & 5/24 \\ 5/12 & -1/12 \end{bmatrix} \cup \begin{bmatrix} -37/30 & 7/15 \\ 3/5 & -1/5 \end{bmatrix}.$$

Thus $(A_B\, B_B)^{-1} = B_B^{-1}\, A_B^{-1}$.

Now we have, even in case of square bimatrix $A_B$; $\left(A_B^{-1}\right)^{-1} = A_B$. We illustrate this by the simple example.

***Example 1.6.8:*** Let

$$A_B = \begin{bmatrix} 2 & 1 \\ 5 & 3 \end{bmatrix} \cup \begin{bmatrix} 1 & 0 \\ 2 & 3 \end{bmatrix}$$

$$A_B^{-1} = \begin{bmatrix} 3 & -1 \\ -5 & 2 \end{bmatrix} \cup \begin{bmatrix} 1 & 0 \\ -2/3 & 1/3 \end{bmatrix}$$

$$\left(A_B^{-1}\right)^{-1} = \begin{bmatrix} 2 & 1 \\ 5 & 3 \end{bmatrix} \cup \begin{bmatrix} 1 & 0 \\ 2 & 3 \end{bmatrix} = A_B .$$

The square bimatrix $A_B$ is non bisingular if $|A_B| \neq 0$. If $|A_B| = 0$ then the bimatrix is bisingular. If $A_B = A_1 \cup A_2$ if one of $A_1$ or $A_2$ is non singular we call $A_B$ to be semi bisingular. The semi bisingular concept is present only in case of bimatrices we have seen non bisingular bimatrices.

Now we proceed on to define singular and semi singular bimatrices.

***Example 1.6.9:*** Let

$$A_B = \begin{bmatrix} 0 & 7 \\ 0 & 5 \end{bmatrix} \cup \begin{bmatrix} 3 & 8 \\ 6 & 16 \end{bmatrix}$$

be the bimatrix clearly $A_B$ is bisingular as $|A_B| = (0, 0)$.

***Example 1.6.10:*** Let

$$= \begin{bmatrix} 1 & 5 \\ 5 & 25 \end{bmatrix} \cup \begin{bmatrix} 1 & 0 \\ 2 & 3 \end{bmatrix}$$

be the square bimatrix $A_B$ is semi bisingular for

$$\begin{vmatrix} 1 & 5 \\ 5 & 25 \end{vmatrix} = 0$$

and

$$\begin{vmatrix} 1 & 0 \\ 2 & 3 \end{vmatrix} = 3$$

and

$$\begin{bmatrix} 1 & 0 \\ 2 & 3 \end{bmatrix}^{-1} = \begin{bmatrix} 1 & 0 \\ -2/3 & 1/3 \end{bmatrix}.$$

The inverse of square bimatrices can be also carried out using partitioning method.

## 1.7 Special Types of Bimatrices

Now we proceed on to define some special types of bimatrices called overlap row (row overlap) and column overlap bimatrices for they will find their application in future.

**DEFINITION 1.7.1:** *Let $A_B = A_1 \cup A_2$ be a rectangular bimatrix. This bimatrix $A_B$ is said to be a row overlap rectangular bimatrix if the matrices $A_1$ and $A_2$ have at least one row with same entries.*

***Example 1.7.1:*** Let

$$A_B = \begin{bmatrix} 3 & 0 & 1 & 1 & 2 \\ 0 & 1 & 1 & 0 & -1 \\ 4 & 4 & 2 & 2 & 3 \\ 0 & 1 & 1 & 1 & 1 \end{bmatrix} \cup \begin{bmatrix} 2 & 0 & 1 & 1 & 1 \\ 0 & 1 & 1 & 1 & 1 \\ 4 & -1 & 2 & 2 & 3 \\ 4 & 4 & 2 & 2 & 3 \end{bmatrix}$$

be a bimatrix which is a $4 \times 5$ bimatrix. This is a row over lap bimatrix for we have the rows (0 1 1 1 1) and (4 4 2 2 3) to be rows with same entries of both the matrices.

Likewise we define column overlap bimatrix.

**DEFINITION 1.7.2:** *A $m \times n$ bimatrix $A_B = A_1 \cup A_2$ is said to be column overlap bimatrix if the matrices $A_1$ and $A_2$ has at least a column with same entries.*

***Example 1.7.2:*** Let

$$A_B = \begin{bmatrix} 3 & 1 & 0 & 1 & 3 \\ 0 & 1 & 4 & 0 & 0 \\ 0 & 1 & 1 & -1 & 4 \end{bmatrix} \cup \begin{bmatrix} -4 & 1 & 2 & 0 & 1 \\ 4 & 1 & 0 & 0 & 3 \\ 2 & 1 & 0 & 0 & 5 \end{bmatrix}$$

be the $3 \times 5$ bimatrix. This is a column overlap bimatrix for

$$\begin{bmatrix} 1 \\ 1 \\ 1 \end{bmatrix}$$

is the column which is both in $A_1$ and $A_2$.

Now we proceed on to define row overlap and column overlap rectangular bimatrix (row column overlap).

**DEFINITION 1.7.3:** *Let $A_B = A_1 \cup A_2$ be a $m \times n$ (rectangular) bimatrix where $A_1$ and $A_2$ have columns which have same entries and $A_1$ and $A_2$ have rows which have same entries, then the bimatrix is called as the row-column overlap bimatrix.*

We illustrate this by the following example.

***Example 1.7.3:*** Let

$$A_B = \begin{bmatrix} 3 & 1 & 0 & 2 & 4 & 5 \\ 6 & 3 & 6 & 3 & 6 & 3 \\ 1 & 2 & 0 & 0 & 1 & 1 \\ 0 & 1 & 0 & 1 & 0 & 1 \end{bmatrix} \cup \begin{bmatrix} 6 & -2 & 0 & 1 & 3 & 5 \\ 6 & 3 & 6 & 3 & 6 & 3 \\ 7 & 5 & 0 & 0 & 1 & 1 \\ 0 & 1 & 0 & 1 & 0 & 1 \end{bmatrix}.$$

The bimatrix $A_B$ has

$$( 6\ 3\ 6\ 3\ 6\ 3) \text{ and } (0\ 1\ 0\ 1\ 0\ 1)$$

as common rows and

$$\begin{bmatrix} 0 \\ 6 \\ 0 \\ 0 \end{bmatrix} \text{ and } \begin{bmatrix} 5 \\ 3 \\ 1 \\ 1 \end{bmatrix}$$

are common columns. Thus $A_B$ is a row-column overlap bimatrix.

Thus the concept of row overlap bimatrix, column overlap bimatrix, and row column overlap bimatrix would be helpful in future applications.

**Chapter Two**

# BIVECTOR SPACES AND THEIR PROPERTIES

Vector space study exploit the properties of matrices, like linear transformation, linear operator, eigen values etc. So we show in this chapter how the concept of bimatrices are used in bivector spaces.

This chapter has six sections; we first recall the definition of bivector spaces and some of its basic properties in section one. In section two for the first time the notion of systems of linear biequations are introduced. Section three deals with bitransformation in bivector space, linear bioperator is introduced in section four. Section Five introduces and analysis the notions of bieigen values and bieigen vectors. In section six we define bidiagonalization and give some its properties. Several other analogous result can be derived in case of bivector spaces.

## 2.1 Introduction to Bivector spaces

In this section we introduce the notion of bivector spaces and illustrate them with examples.

**DEFINITION 2.1.1:** *Let $V = V_1 \cup V_2$ where $V_1$ and $V_2$ are two proper subsets of V and $V_1$ and $V_2$ are vector spaces over the same field F that is V is a bigroup, then we say V is a bivector space over the field F.*

*If one of $V_1$ or $V_2$ is of infinite dimension then so is V. If $V_1$ and $V_2$ are of finite dimension so is V; to be more precise if $V_1$ is of dimension n and $V_2$ is of dimension m then we define dimension of the bivector space $V = V_1 \cup V_2$ to be of dimension m + n. Thus there exists only m + n elements which are linearly independent and has the capacity to generate $V = V_1 \cup V_2$.*

*The important fact is that same dimensional bivector spaces are in general not isomorphic.*

***Example 2.1.1:*** Let $V = V_1 \cup V_2$ where $V_1$ and $V_2$ are vector spaces of dimension 4 and 5 respectively defined over rationals where $V_1 = \{(a_{ij}) / a_{ij} \in Q\}$, collection of all $2 \times 2$ matrices with entries from Q. $V_2$ = {Polynomials of degree less than or equal to 4}.

Clearly V is a finite dimensional bivector space of dimension 9. In order to avoid confusion we in time of need follow the following convention. If $v \in V = V_1 \cup V_2$ then $v \in V_1$ or $v \in V_2$ if $v \in V_1$ then v has a representation of the form $(x_1, x_2, x_3, x_4, 0, 0, 0, 0, 0)$ where $(x_1, x_2, x_3, x_4) \in V_1$ if $v \in V_2$ then $v = (0, 0, 0, 0, y_1, y_2, y_3, y_4, y_5)$ where $(y_1, y_2, y_3, y_4, y_5) \in V_2$.

Thus we follow the notation.

***Notation***: Let $V = V_1 \cup V_2$ be the bivector space over the field F with dimension of V to be m + n where dimension of $V_1$ is m and that of $V_2$ is n. If $v \in V = V_1 \cup V_2$, then $v \in V_1$ or $v \in V_2$ if $v \in V_1$ then $v = (x_1, x_2, \ldots, x_m, 0, 0, \ldots, 0)$ if $v \in V_2$ then $v = (0, 0, \ldots, 0, y_1, y_2, \ldots, y_n)$. Contextually we also accept representation of $v \in V$ by $v = v_1 \cup v_2$ where $v_1 = (x_1, \ldots, x_m)$ and $v_2 = (y_1, \ldots, y_n)$.

We never add elements of $V_1$ and $V_2$. We keep them separately as no operation may be possible among them. For in example we had $V_1$ to be the set of all $2 \times 2$ matrices with entries from Q where as $V_2$ is the collection of all polynomials of degree less than or equal to 4.

So no relation among elements of $V_1$ and $V_2$ is possible. Thus we also show that two bivector spaces of same dimension need not be isomorphic by the following example:

***Example 2.1.2:*** Let $V = V_1 \cup V_2$ and $W = W_1 \cup W_2$ be any two bivector spaces over the field F. Let V be of dimension 8 where $V_1$ is a vector space of dimension 2, say $V_1 = F \times F$ and $V_2$ is a vector space of dimension 6 say all polynomials of degree less than or equal to 5 with coefficients from F. W be a bivector space of dimension 8 where $W_1$ is a vector space of dimension 3 i.e.

$$W_1 = \{\text{all polynomials of degree less than or equal to 2}\} \text{ with coefficients from F}$$

and

$$W_2 = F \times F \times F \times F \times F$$

a vector space of dimension 5 over F. Thus any vector in V is of the form $(x_1, x_2, 0, 0, 0, \ldots, 0)$ or $(0, 0, y_1, y_2, \ldots, y_6)$ and any vector in W is of the form $(x_1, x_2, x_3, 0, \ldots, 0)$ or $(0, 0, 0, y_1, y_2, \ldots, y_5)$. Hence no isomorphism can be sought between V and W in this set up.

This is one of the marked difference between the vector spaces and bivector spaces. Thus we have the following theorems, the proof of which is left for the reader to prove.

**THEOREM 2.1.1:** *Bivector spaces of same dimension defined over same fields need not in general be isomorphic.*

**THEOREM 2.1.2:** *Let $V = V_1 \cup V_2$ and $W = W_1 \cup W_2$ be any two bivector spaces of same dimension over the same field F. Then V and W are isomorphic as bivector spaces if and only if the vector space $V_1$ is isomorphic to $W_1$ and the vector space $V_2$ is isomorphic to $W_2$, that is dimension of $V_1$ is equal to dimension $W_1$ and the dimension of $V_2$ is equal to dimension $W_2$.*

**THEOREM 2.1.3:** *Let $V = V_1 \cup V_2$ be a bivector space over the field F. W any non empty set of V. $W = W_1 \cup W_2$ is a sub-bivector space of V if and only if $W \cap V_1 = W_1$ and $W \cap V_2 = W_2$ are sub spaces of $V_1$ and $V_2$ respectively.*

*Proof*: Direct; left for the reader to prove.

**DEFINITION 2.1.2:** *Let $V = V_1 \cup V_2$ and $W = W_1 \cup W_2$ be two bivector spaces defined over the field F of dimensions p = m + n and q = $m_1 + n_1$ respectively.*

*We say the map $T: V \rightarrow W$ is a bilinear transformation of the bivector spaces if $T = T_1 \cup T_2$ where $T_1 : V_1 \rightarrow W_1$ and $T_2 : V_2 \rightarrow W_2$ are linear transformations from vector spaces $V_1$ to $W_1$ and $V_2$ to $W_2$ respectively satisfying the following two rules.*

i. *$T_1$ is always a linear transformation of vector spaces whose first co ordinates are non-zero and $T_2$ is a linear transformation of the vector space whose last co ordinates are non zero.*

ii. *$T = T_1 \cup T_2$ '$\cup$' is just only a notational convenience.*

iii. *$T(v) = T_1 (v)$ if $v \in V_1$ and $T (v) = T_2 (v)$ if $v \in V_2$.*

*Yet another marked difference between bivector spaces and vector spaces are the associated matrix of an operator of bivector spaces which has $m_1 + n_1$ rows and m + n columns where dimension of V is m + n and dimension of W is $m_1 + n_1$ and T is a linear transformation from V to W. If A is the associated matrix of T then.*

$$A = \begin{bmatrix} B_{m_1 \times m} & O_{n_1 \times m} \\ O_{m_1 \times n} & C_{n_1 \times n} \end{bmatrix}$$

*where A is a $(m_1 + n_1) \times (m + n)$ matrix with $m_1 + n_1$ rows and m + n columns. $B_{m_1 \times m}$ is the associated matrix of $T_1 : V_1 \rightarrow W_1$ and $C_{n_1 \times n}$ is the associated matrix of $T_2 : V_2 \rightarrow W_2$ and $O_{n_1 \times m}$ and $O_{m_1 \times n}$ are non zero matrices.*

***Example 2.1.3:*** Let $V = V_1 \cup V_2$ and $W = W_1 \cup W_2$ be two bivector spaces of dimension 7 and 5 respectively defined over the field F with dimension of $V_1 = 2$, dimension of $V_2 = 5$, dimension of $W_1 = 3$ and dimension of $W_2 = 2$. T be a linear transformation of bivector spaces V and W. The associated matrix of $T = T_1 \cup T_2$ where $T_1 : V_1 \rightarrow W_1$ and $T_2 : V_2 \rightarrow W_2$ given by

$$A = \begin{bmatrix} 1 & -1 & 2 & 0 & 0 & 0 & 0 & 0 \\ -1 & 3 & 0 & 0 & 0 & 0 & 0 & 0 \\ 0 & 0 & 0 & 2 & 0 & 1 & 0 & 0 \\ 0 & 0 & 0 & 3 & 3 & -1 & 2 & 1 \\ 0 & 0 & 0 & 1 & 0 & 1 & 1 & 2 \end{bmatrix}$$

where the matrix associated with $T_1$ is given by

$$\begin{bmatrix} 1 & -1 & 2 \\ -1 & 3 & 0 \end{bmatrix}$$

and that of $T_2$ is given by

$$\begin{bmatrix} 2 & 0 & 1 & 0 & 0 \\ 3 & 3 & -1 & 0 & 1 \\ 1 & 0 & 1 & 1 & 2 \end{bmatrix}.$$

We call $T : V \rightarrow W$ a linear operator of both the bivector spaces if both V and W are of same dimension. So the matrix A associated with the linear operator T of the bivector spaces will be a square matrix. Further we demand that the spaces V and W to be only isomorphic bivector spaces. If we want to define eigen bivalues and eigen bivectors associated with T.

The eigen bivector values associated with are the eigen values associated with $T_1$ and $T_2$ separately. Similarly the eigen bivectors are that of the eigen vectors associated with $T_1$ and $T_2$ individually. Thus even if the dimension of the bivector spaces V and W are equal still we may not have eigen bivalues and eigen bivectors associated with them.

***Example 2.1.4:*** Let T be a linear operator of the bivector spaces – V and W. $T = T_1 \cup T_2$ where $T_1 : V_1 \to W_1$ dim $V_1$ = dim $W_1$ = 3 and $T_2 : V_2 \to W_2$ where dim $V_2$ = dim $W_2$ = 4. The associated matrix of T is

$$A = \begin{bmatrix} 2 & 0 & -1 & 0 & 0 & 0 & 0 \\ 0 & 1 & 0 & 0 & 0 & 0 & 0 \\ -1 & 0 & 3 & 0 & 0 & 0 & 0 \\ 0 & 0 & 0 & 2 & -1 & 0 & 6 \\ 0 & 0 & 0 & -1 & 0 & 2 & 1 \\ 0 & 0 & 0 & 0 & 2 & -1 & 0 \\ 0 & 0 & 0 & 6 & 1 & 0 & 3 \end{bmatrix}.$$

The eigen bivalues and eigen bivectors can be calculated.

**DEFINITION 2.1.3:** *Let T be a linear operator on a bivector space V. We say that T is diagonalizable if $T_1$ and $T_2$ are diagonalizable where $T = T_1 \cup T_2$.*

*The concept of symmetric operator is also obtained in the same way, we say the linear operator $T = T_1 \cup T_2$ on the bivector space $V = V_1 \cup V_2$ is symmetric if both $T_1$ and $T_2$ are symmetric.*

**DEFINITION 2.1.4:** *Let $V = V_1 \cup V_2$. be a bivector space over the field F. We say $\langle\,,\rangle$ is an inner product on V if $\langle\,,\rangle = \langle\,,\rangle_1 \cup \langle\,,\rangle_2$ where $\langle\,,\rangle_1$ and $\langle\,,\rangle_2$ are inner products on the vector spaces $V_1$ and $V_2$ respectively.*

*Note that in $\langle\,,\rangle = \langle\,,\rangle_1 \cup \langle\,,\rangle_2$ the '$\cup$' is just a conventional notation by default.*

**DEFINITION 2.1.5:** *Let $V = V_1 \cup V_2$ be a bivector space on which is defined an inner product $\langle\,,\rangle$. If $T = T_1 \cup T_2$ is a* *linear operator on the bivector spaces V we say $T^*$ is an adjoint of T if $\langle T\alpha / \beta \rangle = \langle \alpha / T^*\beta \rangle$ for all $\alpha, \beta \in V$ where $T^* = T_1^* \cup T_2^*$ are $T_1^*$ is the adjoint of $T_1$ and $T_2^*$ is the adjoint of $T_2$.*

The notion of normal and unitary operators on the bivector spaces are defined in an analogous way. T is a unitary operator on the bivector space $V = V_1 \cup V_2$ if and only if $T_1$ and $T_2$ are unitary operators on the vector space $V_1$ and $V_2$ .

Similarly T is a normal operator on the bivector space if and only if $T_1$ and $T_2$ are normal operators on $V_1$ and $V_2$ respectively. We can extend all the notions on bivector spaces $V = V_1 \cup V_2$ once those properties are true on $V_1$ and $V_2$.

The primary decomposition theorem and spectral theorem are also true is case of bivector spaces. The only problem with bivector spaces is that even if the dimension of bivector spaces are the same and defined over the same field still they are not isomorphic in general.

Now we are interested in the collection of all linear transformation of the bivector spaces $V = V_1 \cup V_2$ to $W = W_1 \cup W_2$ where V and W are bivector spaces over the same field F.

We denote the collection of linear transformation by B-$Hom_F(V, W)$.

**THEOREM 2.1.4:** *Let V and W be any two bivector spaces defined over F. Then B-$Hom_F(V, W)$ is a bivector space over F.*

*Proof:* Given $V = V_1 \cup V_2$ and $W = W_1 \cup W_2$ be two bivector spaces defined over the field F. B-$Hom_F(V, W) = \{T_1 : V_1 \to W_1\} \cup \{T_2 : V_2 \to W_2\} = Hom_F(V_1, W_1) \cup Hom_F(V_2, W_2)$. So clearly B- $Hom_F(V,W)$ is a bivector space as $Hom_F(V_1, W_1)$ and $Hom_F(V_2, W_2)$ are vector spaces over F.

**THEOREM 2.1.5:** *Let $V = V_1 \cup V_2$ and $W = W_1 \cup W_2$ be two bivector spaces defined over F of dimension $m + n$ and $m_1 + n_1$ respectively. Then B-$Hom_F(V,W)$ is of dimension $mm_1 + nn_1$.*

*Proof*: Obvious by the associated matrices of T.

Thus it is interesting to note unlike in other vector spaces the dimension of $Hom_F(V, W)$ is mn if dimension of the vector space V is m and that of the vector space W is n. But in case of bivector spaces of dimension m + n and $m_1 + n_1$ the dimension of B-$Hom_F$ (V, W) is not $(m + n)(m_1 + n_1)$ but $mm_1 + nn_1$, which is yet another marked difference between vector spaces and bivector spaces.

Further even if bivector space V and W are of same dimension but not isomorphic we may not have B-$Hom_F(V, W)$ to be a bilinear algebra analogous to linear algebra. Thus B-$Hom_F(V, W)$ will be a bilinear algebra if and only if the bivector spaces V and W are isomorphic as bivector spaces.

Now we proceed on to define the concept of pseudo bivector spaces.

**DEFINITION 2.1.6:** *Let V be an additive group and $B = B_1 \cup B_2$ be a bifield if V is a vector space over both $B_1$ and $B_2$ then we call V a pseudo bivector space.*

***Example 2.1.5:*** Let V = R the set of reals,

$$B = Q\left(\sqrt{3}\right) \cup Q\left(\sqrt{2}\right)$$

be the bifield. Clearly R is a pseudo bivector space over B. Also if we take $V_1 = R \times R \times R$ then $V_1$ is also a pseudo bivector space over B.

Now how to define dimension, basis etc of V, where V is a pseudo bivector space.

**DEFINITION 2.1.7:** *Let V be a pseudo bivector space over the bifield $F = F_1 \cup F_2$ . A proper subset $P \subset V$ is said to be*

*a pseudo sub-bivector space of V if P is a vector space over $F_1$ and P is a vector space over $F_2$ that is P is a pseudo vector space over F.*

***Example 2.1.6:*** Let V = R × R × R be a pseudo bivector space over

$$F = Q(\sqrt{3}) \cup Q(\sqrt{2}).$$

P = R × {0} × {0} is a pseudo sub-bivector space of V as P is a pseudo bivector space over F.

## 2.2 Systems of linear biequations

In this section we start to view the concept of bimatrices as systems of equation and develop the notion is bieigen values, bieigen vectors, bitransformations, bioperators etc. We give only the basis and in several places we leave the reader to supply proofs for our main motivation is introduce several new concepts and the proofs which are left for the reader is more a matter of routine for they can be derived easily and this will make the reader more at home with the subject. Any way all caution is taken in the earlier chapter to make the reader feel at home with the concept of bimatrices.

Now we first introduce the new notion of biequation

$$\begin{array}{lc} A^1_{11}\, x^1_1 \cup A^2_{11}\, x^2_1 + \ldots + A^1_{1n}\, x^1_n \cup A^2_{1n}\, x^2_n = & y^1_1 \cup y^2_1 \\ A^1_{21}\, x^1_1 \cup A^2_{21}\, x^2_1 + \ldots + A^1_{2n}\, x^1_n \cup A^2_{2n}\, x^2_n & y^1_2 \cup y^2_2 \\ \vdots & \vdots \\ A^1_{m1}\, x^1_1 \cup A^2_{m2}\, x^2_1 + \ldots + A^1_{mn}\, x^1_n \cup A^2_{mn}\, x^2_n = & y^1_m \cup y^2_m \end{array}$$

where $y^1_1, \ldots, y^1_m$ and $y^2_1, \ldots, y^2_m$ and $A^1_{ij}$ and $A^2_{ij}$ are given elements $1 \le i \le m$, $1 \le j \le n$. $x^1_1, \ldots, x^1_n$ and $x^2_1, x^2_2, \ldots, x^2_n$ satisfies the sets of equations. We call the set of equations as a system of m linear biequations in 2n or (n + n) unknowns.

Any n-tables $\left(x_1^1, x_2^1, ..., x_n^1\right)$ and $\left(x_1^2, x_2^2, ..., x_n^2\right)$ of elements of a field F which satisfies each of the above equations is called a bisolution of the system.

If

$$y_1^1 = y_2^1 = \cdots = y_m^1 = 0$$

and

$$y_1^2 = y_2^2 = \cdots = y_m^2 = 0 .$$

We say the system is homogeneous biequation. If only one set is zero
i.e.,

$$y_1^1 = y_2^1 = \cdots = y_m^1 = 0$$

or

$$y_1^2 = y_2^2 = \cdots = y_m^2 = 0$$

then we call the system of linear biequations to be semi homogeneous. Note the symbol union is just a symbol.

Perhaps the most basic way of finding solutions to these system of linear biequations is the method of elimination. We can just illustrate this technique by the following example.

***Example 2.2.1:***

$$2x_1^1 \cup 3x_1^2 + x_2^1 \cup 5x_2^2 + x_3^1 \cup x_3^2 = 0$$
$$-x_1^1 \cup x_1^2 + x_2^1 \cup 3x_2^2 + x_3^1 \cup (-x_3^2) = 0$$
$$x_1^1 \cup x_1^2 + x_2^1 \cup 5x_2^2 + x_3^1 \cup 5x_3^2 = 0$$

this is set of homogeneous biequations in the bivariables

$$\left(x_1^1\ x_2^1\ x_3^1\right) \text{and} \left(x_1^2\ x_2^2\ x_3^2\right) .$$

This set of linear biequations can be written as

$$\left.\begin{aligned} 2x_1^1 - x_2^1 + x_3^1 &= 0 \\ -x_1^1 + x_2^1 + x_3^1 &= 0 \\ x_1^1 + x_2^1 + x_3^1 &= 0 \end{aligned}\right\} \qquad (I)$$

$$\left.\begin{aligned} 3x_1^2 + 5x_2^2 + x_3^2 &= 0 \\ x_1^2 + 3x_2^2 - x_3^2 &= 0 \\ x_1^2 + 5x_2^2 + 5x_3^2 &= 0 \end{aligned}\right\} \qquad (II)$$

Eliminating set of equations (I).

$$x_2^1 = -\,x_3^1$$

$$2x_1^1 + x_2^1 + x_3^1 = 0$$

so

$$x_1^1 = 0$$

$\left(0,\, x_2^1, -x_2^1\right)$ is a section of solution to the biequation.

Eliminating set of equations (II)

$$4x_1^2 + 8x_2^2 = 0$$

so

$$x_1^2 = -\,2x_2^2$$

$$2x_2^2 + 6x_3^2 = 0$$

$$x_2^2 = -\,3x_3^2$$

so

$$2x_1^2 - 4x_3^2 = 0$$

$$x_1^2 = -\,2x_3^2$$

so, (0, 0, 0) is a set of solution. That is this set of homogeneous biequation has

{(0, a, –a), (0, 0, 0)}

as a bisolution.

It is to be noted that the two system of linear biequations are equivalent if each biequation in each system is a linear combination of the biequations in the other system.

It is left as an exercise for the reader to prove the following theorem.

**THEOREM 2.2.1:** *Equivalent systems of linear biequations have exactly the same bisolutions.*

Using the analogous method of forming linear combinations of linear equations there is no need to continue writing the unknowns $x_1,\ldots, x_n$ since one actually computes only with the coefficients $A_{ij}$ and the scalars $y_i$, here we form the matrix of bicoefficients for the system of linear biequations in the biunknowns.

$$\left(x_1^1, x_2^1, \ldots, x_n^1\right) \text{ and } \left(x_1^2, x_2^2, \ldots, x_n^2\right)$$

with the bi-coefficients

$$\left(A_{ij}^1 \cup A_{ij}^2\right)$$

and scalar bipairs

$$y_i^1 \text{ and } y_i^2 .$$

Suppose we abbreviate the system by

$$A\left(X^1 \cup X^2\right) = Y^1 \cup Y^2$$

where

$$A = \begin{bmatrix} A_{11}^1 \cup A_{11}^2 & \cdots & A_{1n}^1 \cup A_{1n}^2 \\ \vdots & & \\ A_{m1}^1 \cup A_{m1}^2 & \cdots & A_{mn}^1 \cup A_{mm}^2 \end{bmatrix}$$

$$X_1 = \begin{bmatrix} x_1^1 \\ x_2^1 \\ \vdots \\ x_n^1 \end{bmatrix}, X_2 = \begin{bmatrix} x_1^2 \\ x_2^2 \\ \vdots \\ x_n^2 \end{bmatrix}$$

and

$$y^1 = \begin{bmatrix} y_1^1 \\ y_2^1 \\ \vdots \\ y_m^1 \end{bmatrix} \text{ and } y^2 = \begin{bmatrix} y_1^2 \\ y_2^2 \\ \vdots \\ y_m^2 \end{bmatrix}.$$

$A = A_1 \cup A_2$ is called the bimatrix of bicoefficients of the system, strictly speaking the rectangular biarray described as bimatrix $A = A^1 \cup A^2 = \left(A_{ij}^1\right) \cup \left(A_{ij}^2\right)$ is not a bimatrix but is a representation of a bimatrix.

An $m \times n$ bimatrix over the field F is a function A from the sets of pairs of integers $(i^1, j^1)$, $(i^2, j^2)$, $1 \leq i^1, i^2 \leq m$ and $1 \leq j^1, j^2 \leq n$ in the field F. The entries of the bimatrix are the scalar $A(i^1, j^1) = A^1_{ij}$, $A(i^2, j^2) = A^2_{ij}$ and this is an appropriate method of describing the bimatrix by displaying its entries in a rectangular biarrays having m rows and n columns.

The $X^1$ and $X^2$ defines a $n \times 1$ bimatrix $X = X^1 \cup X^2$ and $Y = Y^1 \cup Y^2$ is a $m \times 1$ bimatrix.

Now we define the linear biequations $AX = Y$.

$$\left(A^1 \cup A^2\right) \left(X^1 \cup X^2\right) = Y^1 \cup Y^2$$
$$A^1 X^1 \cup A^2 X^2 = Y^1 \cup Y^2$$

gives the bimatrix representation of the system of linear biequations.

Now we can define elementary row operations on the bimatrix $AX = Y$ i.e., $A^1 X^1 \cup A^2 X^2 = Y^1 \cup Y^2$.

Now as we have to do pairs of operations they can be carried out in the following ways.

(1) Uniform elementary row operations means the operations is done similarly on both the matrices $A^1 \cup A^2$. Type I elementary operations mean only elementary row operations on $A^1$ and Type II elementary operations mean only elementary row operations on $A^2$.

**DEFINITION 2.2.1:** *Let $A = A_1 \cup A_2$ be a $m \times n$ bimatrix with entries of $A_1$ and $A_2$ from the same field F then this bimatrix A is called the ordinary bimatrix or just a bimatrix. If $A_1$ has entries from the field K and $A_2$ has entries from the field F, F and K distinct then we call the bimatrix as a strong bimatrix. If in the bimatrix $A = A_1 \cup A_2$ if one of $A_1$ or $A_2$ is in the subfield of the field F and $A_1$ or $A_2$ is in the field K, $(F \subset K)$ then we call A the weak bimatrix.*

This three types of bimatrix; weak bimatrix strong bimatrix and ordinary bimatrix or bimatrix, all of them will have a role to play when we work on eigen bivalues and eigen bivectors.

Now we proceed on to define row biequivalence of two ordinary bimatrices.

**DEFINITION 2.2.2:** *If $A = A_1 \cup A_2$ and $B = B_1 \cup B_2$ be $m \times n$ bimatrices over the field F, we say that B is row biequivalent to A if B can be obtained from A by a finite sequence of elementary row operations.*

Let $A = A_1 \cup A_2$ be a bimatrix A is called row bireduced if the first non entry in each non zero rows of $A_1$ and $A_2$ is equal to 1.

Each column of $A_1$ and $A_2$ which contains the leading non zero entry of some row has all its other entries 0.

The identity bimatrix $I = I_1 \cup I_2$

$$= \delta^1_{ij} = \begin{cases} 0 \text{ if } i \neq j \\ 1 \text{ if } i = j \end{cases}$$

and

$$\delta^2_{ij} = \begin{cases} 0 \text{ if } i \neq j \\ 1 \text{ if } i = j \end{cases}$$

where

$$I_1 = \left(\delta^1_{ij}\right) \text{ and } I_2 = \left(\delta^2_{ij}\right).$$

If the simple row bioperation which makes the bimatrix a row bireduced one is such that the operations are done simultaneously for both matrices $A_1$ and $A_2$ then row bireduced bimatrix is said to be a strongly row bireducible bimatrix, other wise it is said to be just row bireducible or at times a weakly row bireducible bimatrix.

***Example 2.2.2:*** Let

$$A = \begin{bmatrix} 3 & -2 & 1 \\ 3 & 2 & 5 \\ 1 & 0 & 1 \end{bmatrix} \cup \begin{bmatrix} 6 & 7 & 1 \\ 0 & -7 & 2 \\ 1 & 0 & 2 \end{bmatrix}$$

$$\rightarrow \begin{bmatrix} 6 & 0 & 6 \\ 3 & 2 & 5 \\ 1 & 0 & 1 \end{bmatrix} \cup \begin{bmatrix} 6 & 0 & 3 \\ 0 & -7 & 2 \\ 1 & 0 & 2 \end{bmatrix}$$

$$\rightarrow \begin{bmatrix} 2 & 0 & 2 \\ 3 & 2 & 5 \\ 1 & 0 & 1 \end{bmatrix} \cup \begin{bmatrix} 2 & 0 & 1 \\ 0 & -7 & 2 \\ 1 & 0 & 2 \end{bmatrix}$$

$$\rightarrow \begin{bmatrix} 0 & 0 & 0 \\ 3 & 2 & 5 \\ 1 & 0 & 1 \end{bmatrix} \cup \begin{bmatrix} 0 & 0 & -3 \\ 0 & -7 & 2 \\ 1 & 0 & 2 \end{bmatrix}$$

$$\rightarrow \begin{bmatrix} 0 & 0 & 0 \\ 3 & 2 & 5 \\ 1 & 0 & 1 \end{bmatrix} \cup \begin{bmatrix} 0 & 0 & 1 \\ 0 & -7 & 2 \\ 1 & 0 & 2 \end{bmatrix}$$

$$\rightarrow \begin{bmatrix} 0 & 0 & 0 \\ 3 & 2 & 5 \\ 1 & 0 & 1 \end{bmatrix} \cup \begin{bmatrix} 0 & 0 & 1 \\ 0 & -7 & 2 \\ 1 & 0 & 0 \end{bmatrix}$$

$$\rightarrow \begin{bmatrix} 0 & 0 & 0 \\ 1 & 0 & 1 \\ 3 & 2 & 5 \end{bmatrix} \cup \begin{bmatrix} 0 & 0 & 1 \\ 1 & 0 & 0 \\ 0 & -7 & 2 \end{bmatrix}.$$

This is the strongly row bireduced bimatrix.

Suppose the same bimatrix is row bireduced by following the method of weakly row bireducibility. We see how best the result is

***Example 2.2.3:*** Let

$$A = \begin{bmatrix} 3 & -2 & 1 \\ 3 & 2 & 5 \\ 1 & 0 & 1 \end{bmatrix} \cup \begin{bmatrix} 6 & 7 & 1 \\ 0 & -7 & 2 \\ 1 & 0 & 2 \end{bmatrix}$$

$$\rightarrow \begin{bmatrix} 6 & 0 & 6 \\ 3 & 2 & 5 \\ 1 & 0 & 1 \end{bmatrix} \cup \begin{bmatrix} 6 & 0 & 3 \\ 0 & -7 & 2 \\ 1 & 0 & 2 \end{bmatrix}$$

$$\rightarrow \begin{bmatrix} 1 & 0 & 1 \\ 3 & 2 & 5 \\ 1 & 0 & 1 \end{bmatrix} \cup \begin{bmatrix} 2 & 0 & 1 \\ 0 & -7 & 2 \\ 1 & 0 & 2 \end{bmatrix}$$

$$\rightarrow \begin{bmatrix} 0 & 0 & 0 \\ 3 & 2 & 5 \\ 1 & 0 & 1 \end{bmatrix} \cup \begin{bmatrix} 0 & 0 & 1 \\ 0 & -7 & 2 \\ 1 & 0 & 2 \end{bmatrix}$$

$$\rightarrow \begin{bmatrix} 0 & 0 & 0 \\ 0 & 2 & 2 \\ 1 & 0 & 1 \end{bmatrix} \cup \begin{bmatrix} 0 & 0 & 1 \\ 0 & -7 & 2 \\ 1 & 0 & 2 \end{bmatrix}$$

$$\rightarrow \begin{bmatrix} 0 & 0 & 0 \\ 0 & 1 & 1 \\ 1 & 0 & 1 \end{bmatrix} \cup \begin{bmatrix} 0 & 0 & 1 \\ 0 & 1 & 0 \\ 1 & 0 & 0 \end{bmatrix}.$$

We see by weakly row bireducibility we see the bimatrix is more simple.

However we cannot say that weakly row bireducibility is better than strong row bireducibility for we see that when in some problems the reduction is to be taken relatively we cannot shift to weakly row bireducibility method.

When no condition is imposed then only we can use any method. For if the modeled equations by these matrices are interdependent then only strongly row bireducible method is adopted.

In case of strongly row bireduced bimatrices we do not expect every rule to be followed by both the matrices we expect only one of the component matrix to satisfy all the necessary conditions, but in case of weakly row bireducibility both the component matrices will satisfy all the conditions.

The following theorem can be easily proved as in case of bimatrices.

**THEOREM 2.2.2:** *If $A = A_1 \cup A_2$ is a rectangular m x n bimatrix and m < n, then the homogeneous system of linear biequations AX = 0 has a non trivial bisolution.*

*Proof:* As $A = A_1 \cup A_2$ and $X = X^1 \cup X^2$ we have AX = 0 is the same as

$$A_1 X^1 \cup A_2 X^2 = 0.$$

So both the components $A_1 X^1 = 0$ and $A_2 X^2 = 0$ has non trivial solution so the bimatrix equation has a nontrivial solution.

The next theorem can also be analogously proved.

**THEOREM 2.2.3:** *If $A = A_1 \cup A_2$ is a $n \times n$ square bimatrix, then A is row biequivalent to the $n \times n$ identity bimatrix if* *and only if the system of biequations AX = 0 has only the trivial bisolution.*

Note we accept for these proof only weakly row bireducibility of the bimatrix $A = A_1 \cup A_2$.

**DEFINITION 2.2.3:** *Let $A = A_1 \cup A_2$ be a $m \times n$ bimatrix over the field F and let $B = B_1 \cup B_2$be an $n \times p$ bimatrix over F. The product AB is a $m \times p$ bimatrix C where $(i^1, j^1, i^2, j^2)$ entry is*

$$C_{ij} = C_{i^1 j^1} \cup C_{i^2 j^2} = \sum_{r_1=1}^{n} A_{i^1 r_1} B_{r_1 j^1} \cup \sum_{r_2=1}^{n} A_{i^2 r_2} B_{r_2 j^2} .$$

## 2.3 Bitransformation in bivector spaces

Linear bitransformation of bivector spaces result in associated bimatrices analogous to linear transformation associated with matrices. We assume throughout $V = V_1 \cup V_2$ is only a finite dimensional bivector space.

Now we proceed on to define linear Bitransformation in a different way for this way of representation will use bimatrices and is both time saving and economic for instead of using $(m + n)^2$ entries we use only $m^2 + n^2$ entries.

**DEFINITION 2.3.1:** *$V = V_1 \cup V_2$ and $W = W_1 \cup W_2$ be bivector spaces over the bifield. A linear bitransformation V into W is a bifunction T from V into W such that*

$$T_1 : V_1 \to W_1$$

*and*

$$T_2 : V_2 \to W_2$$

*are linear transformations from $V_1$ to $W_1$ and $V_2$ to $W_2$ respectively and $T = T_1 \cup T_2$ is the way by which the linear bitransformation is represented.*

*The bitransformation from V to W i.e., from $V_1 \cup V_2$ to $W_1 \cup W_2$ can be always represented by bimatrices.*

We first represent this notion by an example.

***Example 2.3.1:*** $V = V_1 \cup V_2$ and $W = W_1 \cup W_2$ where $V_1 = R \times R \times R$, $V_2 = Q \times Q$ and $W_1 = Q \times Q \times Q$ and $W_2 = R \times R$ be bivector spaces over Q.$T = T_1 \cup T_2 : V \to W$ is given by $T_1 (x, y, z) = (3x + 2y + x, y, x + y)$ and $T_2 (x', y') = (2x', x' + 9y')$. We see the associated bimatrix is a square mixed bimatrix in fact an ordinary bimatrix for both have their entries in the field Q.

***Example 2.3.2:*** Let $V = V_1 \cup V_2$ and $W = W_1 \cup W_2$ be two bivector spaces over Q here

$$V_1 = R \times R \times R \times R$$

and $$V_2 = Q \times Q \times Q,$$

$$W_1 = R \times R \times R$$

and $$W_2 = Q \times Q \times Q \times Q \times Q.$$

Let T: V $\to$ W be a linear bitransformation whose associated bimatrix A is given by

$$\begin{bmatrix} 0 & 1 & 2 & 1 \\ 0 & 0 & -1 & 2 \\ 1 & 2 & 0 & 0 \end{bmatrix} \cup \begin{bmatrix} 0 & 2 & 1 \\ 1 & 1 & 0 \\ 0 & 1 & 2 \\ 2 & 0 & 0 \\ 5 & 1 & 2 \end{bmatrix}$$

$$T = T_1 \cup T_2 : V \to W$$

where

$$T_1 (x, y, z, w) = (y + 2z + w, -z + 2w, x + 2y)$$

and

$$T_2 (x_1, x_2, x_3) = (2x_2 + x_3, x_1 + x_2, x_2 + 2x_3, 2x_1, 5x_1 + x_2 + 2x_3).$$

Clearly the bimatrix associated with the transformation is a mixed rectangular bimatrix. This is also an ordinary bimatrix for its entries are only from the same field Q. It is important to note by using the method of bimatrix we have reduced from 56 entries to just only 27 entries.

It is left as an exercise for the reader to prove if $V = V_1 \cup V_2$ and $W = W_1 \cup W_2$ are two bivector spaces such that dim $V_1$ = m and dim $V_2$ = n, dim $W_1$ = p and dim $W_2$ = q. Both W and V are defined over Q.

If T: V $\to$ W be a linear bitransformation then the bimatrix associated with T is a rectangular bimatrix $A = A_1 \cup A_2$ where $A_1$ is a $p \times m$ matrix and $A_2$ is a $q \times n$ matrix with entries from Q.

**DEFINITION 2.3.2:** *Let V and W be bivector spaces over the field F and let T be a linear bitransformation from V into W. The null bispace of T is the set of all bivectors $\alpha = (\alpha^1, \alpha^2)$ in $V = V_1 \cup V_2$ such that*

$$
\begin{aligned}
T(\alpha) &= (T_1 \cup T_2)(\alpha^1, \alpha^2) \\
&= T_1(\alpha^1) \cup T_2(\alpha^2) \\
&= (0, 0).
\end{aligned}
$$

*If V is a finite dimensional the birank of T is the dimension of the range of T and the binullity of T is the dimension of the null bispace of T.*

***Example 2.3.3:*** Let $V = V_1 \cup V_2$ and $W = W_1 \cup W_2$ be two bivector spaces T : V $\to$ W be a linear bitransformation to find the null space of T.

Let

$$
\begin{aligned}
V_1 &= R \times R \\
V_2 &= Q \times Q \times Q \\
W_1 &= Q \times Q \\
W_2 &= R \times R \times R
\end{aligned}
$$

$$T = T_1 \cup T_2$$

where

$$
\begin{aligned}
T_1(x, y) &= (x - y, x + y) \\
T_2(x\ y\ z) &= (x + y, x - z, x + y + z)
\end{aligned}
$$

any $\alpha \in V = \{(x, y), (x, y, z)\}$ where $(x, y) \in V_1$ and $(x, y, z) \in V_2$.

T (x, y) = (0, 0) i.e., $x + y = 0$

$x - y = 0$

$x = 0$
$y = 0.$

The null space of $T_1$ is 0.

$T_2 (x, y, z) = 0$ so

$$x + y = 0$$
$$x – z = 0$$
$$x + y + z = 0$$

so $z = 0$

$$z = 0$$

so $x = 0$

$$T_2 (x\ y\ z) = 0.$$

Thus the null bispace is T is the {(0, 0)} subbispace of $V = V_1 \cup V_2$.

If $V = V_1 \cup V_2$ and $W = = W_1 \cup W_2$ be bivector spaces over F. Let $T = T_1 \cup T_2$ and $U = U_1 \cup U_2$ be linear bitransformations of V in to W. The function T + U defined by

$$[(T_1 \cup T_2) + (U_1 \cup U_2)]\ \alpha$$
$$= [(T_1 + U_1) \cup (T_2 + U_2)]\ \alpha$$
$$= [(T_1 + U_1) \cup (T_2 + U_2)]\ (\alpha^1, \alpha^2)$$

(where $\alpha^1 \in V_1$ and $\alpha^2 \in V_2$, $\alpha = \alpha^1 \cup \alpha^2$)

$$= (T_1 + U_1)\ (\alpha^1) \cup (T_2 + U_2)\ (\alpha^2)$$
$$= (T_1\ \alpha^1 + U_1\ \alpha^1 ) \cup (T_2\ \alpha^2 + U_2\ \alpha^2 )$$

(using the theorem for when V and W are just vector spaces and T and U are linear transformations from V into W).

Also

$$(T + U)\alpha = (T\alpha + U\alpha).$$
$$= (T_1\ \alpha^1 \cup T_2\ \alpha^2 ) + (U_1\ \alpha^1 \cup U_2\ \alpha^2 )$$
$$= (T_1 \cup T_2)\ (\alpha) + (U_1 \cup U_2)\ \alpha$$
$$= T\ \alpha + U\ \alpha.$$

On similar lines one can prove if T is a linear bitransformation, c any scalar

$$cT\ (\alpha) = c\ (T_1 \cup T_2)\ (\alpha) = c\ (T\alpha).$$

Let $V = V_1 \cup V_2$ and $W = W_1 \cup W_2$ and $Z = Z_1 \cup Z_2$ be three bivector spaces over F. If $T = T_1 \cup T_2$, $U = U_1 \cup U_2$ are linear bitransformation from V into W and W into Z i.e.,

$$T: V \to W$$

and

$$U : W \to Z$$

then UT defined by

$$\begin{aligned}
(UT)(\alpha) &= U(T(\alpha)) \\
&= U[(T_1 \cup T_2)(\alpha^1, \alpha^2) \\
&= U[T(\alpha^1) \cup T(\alpha^2)] \\
&= (U_1 \cup U_2)(T_1(\alpha^1) \cup T_2(\alpha^2)) \\
&= U_1(T_1(\alpha^1)) \cup U_2 T_2(\alpha^2)) \\
&= (U_1 T_1 \cup U_2 T_2)(\alpha^1, \alpha^2) \\
&= (U_1 T_1 \cup U_2 T_2)\alpha \\
&= UT(\alpha)
\end{aligned}$$

is a linear bitransformation from V in to Z.

Now we proceed on to define linear bioperator on the bivector space V defined over the field F.

Let $V = V_1 \cup V_2$ be a bivector space a linear bitransformation $T = T_1 \cup T_2$ is a called linear bioperator from V into V defined by

$$\begin{aligned}
T(\alpha) &= (T_1 \cup T_2)\alpha \\
&= (T_1 \cup T_2)(\alpha^1, \alpha^2) \\
&= T(\alpha^1) \cup T(\alpha^2).
\end{aligned}$$

i.e., a linear Bitransformation is called a linear bioperator if the range bispace and the domain bispace are one and the same.

Now we proceed on to define invertible linear bitransformations. The function $T = T_1 \cup T_2$ from V into W (where V and W are bivector spaces over F, $V = V_1 \cup V_2$ and $W = W_1 \cup W_2$) is called invertible if there exists a function $U = U_1 \cup U_2$ from W into V such that $UT = U_1 T_1 \cup U_2 T_2$ is the identity function i.e., $UT = I_1 \cup I_2$, where $I_1$ is the identity function on $W_1$ and $I_2$ is the identity function on $W_2$.

If T is invertible the function U is unique and is denoted by $T^{-1} = T_1^{-1} \cup T_2^{-1}$. Further more T is invertible if and only if $T = T_1 \cup T_2$ is one to one i.e.,

$$T_1\alpha^1 = T_1 \beta^1 \Leftrightarrow \alpha^1 = \beta^1 \text{ in } V_1$$
$$T_2\alpha^2 = T_2 \beta^2 \Leftrightarrow \alpha^2 = \beta^2 \text{ in } V_2$$

i.e., $(T_1 \cup T_2)(\alpha^1, \alpha^2) = (T_1 \cup T_2)(\beta^1 \beta^2)$
$\Leftrightarrow (T_1 \cup T_2)\alpha = (T_1 \cup T_2)(\beta)$
i.e., $T(\alpha) = T(\beta) \Leftrightarrow \alpha = \beta$.

T is on to, that is if range of T is all of $W = W_1 \cup W_2$.

It is left as an exercise for the reader to prove.

If $V = V_1 \cup V_2$ and $W = W_1 \cup W_2$ be bivector spaces over the field F and let $T = T_1 \cup T_2$ be a linear bitransformation from V into W. If T is invertible, then prove the inverse function $T^{-1} = T_1^{-1} \cup T_2^{-1}$ is a linear bitransformation from W into V.

This T is one to one if and only if T is non singular.

Now it is important to recollect that unlike vector spaces which states same dimension vector spaces over the same field are isomorphic we cannot say in case of bivector spaces over same field but of same dimension are isomorphic.

Now we just state that as in case of vector spaces we can in case of bivector spaces also represent the bitransformation by bimatrices.

Let $V = V_1 \cup V_2$ be a (n, p) dimensional bivector space over the field F and let $W = W_1 \cup W_2$ be a (m, q) dimensional bivector space over F.

Let

$$B_1 = \{\alpha_1^1, ..., \alpha_n^1\}$$

be an ordered basis of $V_1$

$$B_2 = \{\alpha_1^2, ..., \alpha_p^2\}$$

be an ordered basis of $V_2$.

$$B'_1 = \{\beta_1^1, ..., \beta_m^1\}$$

be an ordered basis of $W_1$ and

$$B'_2 = \{\beta_1^2, \ldots, \beta_q^2\}$$

be an ordered basis of $W_2$. $B_1 \cup B_2$ is also known as bibasis of $V = V_1 \cup V_2$. If T is a linear bitransformation from V into W then T is determined by its actions on the bivector $\alpha_j^1 \cup \alpha_j^2$. Each of the (n, p) vectors T $(\alpha_j^1, \alpha_j^2)$ is uniquely expressible as a linear combination of the bibasis $B_1 \cup B_2$.

$$T(\alpha_j^1, \alpha_j^2) = (T_1 \cup T_2)(\alpha_j^1, \alpha_j^2)$$

$$= \sum_{i=1}^{n} \left(A_{ij}^1 \beta_i^1 \cup A_{ij}^2 \beta_i^2\right) \qquad \text{(I)}$$

of the $(\beta_i^1, \beta_j^2)$, the scalars $(A_{ij}^1, A_{2j}^1, \ldots, A_{mj}^1)$ and $(A_{ij}^2, A_{2j}^2, \ldots, A_{mj}^2)$ being coordinates of T $(\alpha_j^1, \alpha_j^2)$ in the order basis $B'_1$ and $B'_2$. Accordingly the bitransformation T = $T_1 \cup T_2$ is determined by the mn $\cup$ pq scalars $(A_{ij}^1, A_{ij}^2)$.

The bimatrix A = $A_1^{m \times n} \cup A_2^{p \times q}$ defined by $A^1(ij) = A_{ij}^1$ $A^2(ij) = A_{ij}^2$ defined by (I). Is called the bimatrix of T = $T_1 \cup T_2$ relative to the pair of ordered basis $(B_1, B_2)$, $(B'_1, B'_2)$.

Our immediate task is to understand explicitly how the bimatrix A = $A_1 \cup A_1$ determines the linear bitransformation T = $T_1 \cup T_2$.

If $\alpha = \left(x_1^1 \alpha_1^1 + \ldots + x_n^1 \alpha_n^1\right) \cup \left(x_1^2 \alpha_1^2 + \ldots + x_p^2 \alpha_p^2\right)$ is a bivector in V = $V_1 \cup V_2$, then

$$\begin{aligned}
T\alpha &= (T_1 \cup T_2)(\alpha^1 \cup \alpha^2) \\
&= T\left(\sum_{j=1}^{n} x_j^1 \alpha_{ij}^1 \cup \sum_{j=1}^{p} x_j^2 \alpha_{ij}^2\right) \\
&= \sum_{j=1}^{m} x_j^1 T_1(\alpha_j^1) \cup \sum_{j=1}^{1} x_j^2 T_2(\alpha_j^2) \\
&= \sum_{j=1}^{n} x_j^1 \sum_{i=1}^{m} A_{ij}^1 \beta_i^1 + \sum_{j=1}^{p} x_j^2 \sum_{i=1}^{q} A_{ij}^2 \beta_i^2
\end{aligned}$$

$$= \sum_{i=1}^{n} \left( \sum_{j=1}^{n} A_{ij}^{1}\ x_{i}^{1} \right) \beta_{i}^{1} + \sum_{i=1}^{q} \left( \sum_{j=1}^{p} A_{ij}^{2}\ x_{j}^{2} \right) \beta_{i}^{2}$$

If X is the co ordinate bimatrix of $\alpha = (\alpha^1\ \alpha^2)$ in the ordered basis $(B_1, B_2)$ then the computation above shows A X is the co ordinate bimatrix of the vector $T(\alpha) = (T_1 \cup T_2)(\alpha^1 \cup \alpha^2)$ in the ordered basis $B'_1, B'_2$ because the scalar.

$$\sum_{j=1}^{n} A_{ij}^{1}\ x_{j}^{1} \ \cup\ \sum_{j=1}^{p} A_{ij}^{2}\ x_{j}^{2}$$

is the $i^{th}$ row of the column bimatrix $A X = (A_1 \cup A_2)(X_1 \cup X_2)$. Let us also observe that if $A = A_1 \cup A_2$ where $A_1$ is a $m \times n$ matrix and $A_2$ is a $p \times q$ matrix over the field F then

$$T \left( \sum_{j=1}^{n} x_{j}^{1}\ \alpha_{j}^{1} \ \cup\ \sum_{j=1}^{p} x_{j}^{2}\ \alpha_{j}^{2} \right) =$$

$$\sum_{i=1}^{m} \left( \sum_{j=1}^{n} A_{ij}^{1}\ x_{i}^{1} \right) \beta_{i}^{1} \cup \sum_{i=1}^{q} \left( \sum_{j=1}^{P} \left( A_{ij}^{2}\ x_{i}^{2} \right) \right) \beta_{i}^{2}$$

defines a linear bitransformations $T = T_1 \cup T_2$ from $V = V_1 \cup V_2$ into $W = W_1 \cup W_2$, the bimatrix which is $A = A_1 \cup A_2$ relative to $(B_1, B_2)$, $(B'_1, B'_2)$.

This result can be formulated formally as, let $V = V_1 \cup V_2$ be a (n, p) dimensional bivector space over the field F and W be an (m, q) dimensional bivector space over F. Let $(B_1, B_2)$ be an ordered basis for $V_1 \cup V_2 = V$ and $(B'_1, B'_2)$ be on ordered basis for $W_1 \cup W_2 = W$, for each linear bitransformation $T = T_1 \cup T_2$ from $V = V_1 \cup V_2$ into $W = W_1 \cup W_2$ there is an $m \times n$, $p \times q$ bimatrix $A = A_1^{m \times n} \cup A_2^{p \times q}$ with entries in F such that

$$[T(\alpha)]_{(B'_1, B'_2)} = A\ [\alpha]_{(B'_1, B'_2)}$$

for any vector $\alpha = (\alpha^1, \alpha^2)$ in $V = V_1 \cup V_2$.

Further more $T \mapsto A$ is a one to one correspondence between the set of all linear bitransformations from $V = V_1 \cup V_2$ into $W = W_1 \cup W_2$ and the set of all bimatrices $(m \times n, p \times q)$ over the field F.

The matrix $A = A_1 \cup A_2$ associated with $T = T_1 \cup T_2$ is called the bimatrix of T relative to the ordered basis $(B_1, B_2)$, $(B'_1, B'_2)$.

Note that the equation.

$$T\left(\alpha_j^1\ \alpha_j^2\right) = \sum_{i=1}^{m} A_{ij}^1 \beta_j^1 \cup \sum_{i=1}^{q} A_{ij}^2 \beta_j^2$$

says that A is a bimatrix whose bicolumns $A_1^1, \ldots, A_n^1$ and $A_1^2, \ldots, A_p^2$ are given by

$$A_j^1 = \left[T_1\ \alpha_j^1\right] B'_1;\ j = 1, 2, \ldots, n$$

$$A_j^2 = \left[T_2\ \alpha_j^2\right] B'_2;\ j = 1, 2, \ldots, p.$$

If U is another linear bitransformation from $V = V_1 \cup V_2$ into $W = W_1 \cup W_2$ and $B^1 = \left[B_1^1, \ldots, B_n^1\right]$ and $B^2 = \left[B_1^2, \ldots, B_p^2\right]$ is the bimatrix U relative to the ordered basis $(B_1, B_2)$, $(B'_1, B'_2)$ then CA + B is a bimatrix of C(T) + U relative to $(B_1, B_2)$, $(B'_1, B'_2)$.

This is clear because

$$\begin{aligned} CA_j + B_j &= C_1\ A_j^1 \cup C_2\ A_j^2 + B_j^1 \cup B_j^2 \\ &= C_1 \left[T_1\ \alpha_j^1\right]_{B'_1} \cup C_2 \left[T_2\ \alpha_j^2\right]_{B'_2} \\ &\quad + \left[U_1\ \alpha_j^1\right]_{B'_1} \cup \left[U_2\ \alpha_j^2\right]_{B'_2} \\ &= \left[CT\left(\alpha_j^1\ \alpha_j^2\right) + \cup\left(\alpha_j^1\ \alpha_j^2\right)\right]_{(B_1^1, B_2^1)} \\ &= \left[\left(CT + U\right)\left(\alpha_j\right)\right]_{[B'_1 B'_2]}. \end{aligned}$$

## 2.4. Linear Bioperators on Bivector spaces

Now we shall discuss when the linear bitransformation is a linear bioperator. We shall be particularly interested in the representation by bimatrices of linear bitransformations of a bispace to itself; i.e., linear bioperators on the bivector space $V = V_1 \cup V_2$. In this case it is most convenient to use the same ordered basis in each case i.e., $B_1 = B'_1$ and $B_2 = B'_2$. We shall then call the representing bimatrix simply the bimatrix of T relative to the ordered basis $(B_1, B_2)$.

Since this concept will be so important to us, we shall review the definition. If $T = T_1 \cup T_2$ is a linear bioperator on the finite dimensional bivector space $V = V_1 \cup V_2$ and $B = (B_1, B_2)$.

$$= \left\{\left(\alpha_1^1, \alpha_2^1 ..., \alpha_n^1\right), \left(\alpha_1^2, \alpha_2^2, ..., \alpha_n^2\right)\right\}$$

is an ordered basis for $V = V_1 \cup V_2$, the bimatrix of $T = T_1 \cup T_2$ relative to $B = (B_1, B_2)$ (or the bimatrix of $T = T_1 \cup T_2$ in the ordered basis $B = (B_1, B_2)$) is the $n \times n$ bimatrix $A = A_1 \cup A_2$ whose entries $A_{ij}^1 \cup A_{ij}^2$ are defined by the equations;

$$T\alpha_j = T\left(\alpha_j^1, \alpha_j^2\right)$$

$$= \sum_{i=1}^{n}\left(A_{ij}^1 \cup A_{ij}^2\right)\left(\alpha_j^1, \alpha_j^2\right) \quad (j = 1, 2, \ldots, n)$$

$$= \sum_{i=1}^{n} A_{ij}^1 \alpha_i^1 \cup \sum_{i=1}^{n} A_{ij}^2 \alpha_i^2 \qquad (j = 1, 2, \ldots, n).$$

Thus

$$T_1 \alpha_j^1 \cup T_2 \alpha_j^2 = \sum_{i=1}^{n} A_{ij}^1 \alpha_i^1 \cup \sum_{i=1}^{n} A_{ij}^2 \alpha_i^2 \text{ ; } (j = 1, 2, \ldots, n).$$

It is to be noted that the bimatrix representing T depends upon the ordered basis $B = (B_1, B_2)$ and that there is a representing bimatrix for $T = T_1 \cup T_2$ in each ordered basis for $V = V_1 \cup V_2$. Here we shall not forget this dependence. We shall use the notation $[T]_B = [T_1]_{B_1} \cup [T_2]_{B_2}$ for the bimatrix of linear bioperator $T = T_1 \cup T_2$ in the ordered basis $B = (B_1, B_2)$. The manner in which this bimatrix and the basis describe the linear bioperator is such that for each $\alpha = (\alpha^1, \alpha^2)$ in V,

$$[T\alpha]_B = [T(\alpha^1 \cup \alpha^2)]_{[B_1 B_2]}$$

$$= [T_1 \cup T_2 (\alpha^1 \cup \alpha^2)]_{[B_1, B_2]}$$

$$= [T_1]_{B_1} [\alpha^1]_{B_1} \cup [T_2]_{B_2} [\alpha^2]_{B_2} .$$

We have seen that what happens when representing the bimatrices when the bitransformations are added namely that the bimatrices add. We should now like to ask what happens when we compose bitransformation. More specifically let $V = V_1 \cup V_2$ into $W = W_1 \cup W_2$ and $Z = Z_1 \cup Z_2$ be bivector spaces over the field F of respective dimensions $(n, n_1)$, $(m, m_1)$ and $(p, p_1)$. Let T be linear bitransformation from $V = V_1 \cup V_2$ into $W = W_1 \cup W_2$ and $U = U_1 \cup U_2$ a linear bitransformation from $W = W_1 \cup W_2$ into $Z = Z_1 \cup Z_2$.

Suppose we have ordered basis.

$$B = \{B_1, B_2\}$$
$$= \left\{\left(\alpha_1^1, \alpha_2^1, \ldots, \alpha_n^1\right), \left(\alpha_1^2, \alpha_2^2, \ldots, \alpha_{n_1}^2\right)\right\}$$

$$B^1 = \{B'_1, B'_2\}$$
$$= \left\{\left(\beta_1^1, \beta_2^1, \ldots, \beta_m^1\right), \left(\beta_1^2, \beta_2^2, \ldots, \beta_{m_1}^2\right)\right\}$$

and

$$B'' = \{B''_1, B''_2\}$$
$$= \left\{\left(\gamma_1^1, \gamma_2^1, \ldots, \gamma_p^1\right), \left(\gamma_1^2, \gamma_2^2, \ldots, \gamma_{p_1}^1\right)\right\}$$

for the respective bispaces $V = V_1 \cup V_2$ , $W = W_1 \cup W_2$ and $Z = Z_1 \cup Z_2$. Let $A = A_1 \cup A_2$ be the bimatrix of $T = T_1 \cup T_2$ relative to the pair.

$B = (B_1, B_2)$, $B' = (B'_1, B'_2)$ and let $B = B_1 \cup B_2$ be the bimatrix of $U = U_1 \cup U_2$ relative to the pair $B' = (B'_1, B'_2)$, and $B'' = (B''_1, B''_2)$. It is then easy to see that the bimatrix C of the bitransformation UT relative to the pair B, $B''$ ($B = (B_1, B_2)$ and $B'' = (B''_1, B''_2)$) is the product of the bimatrices $B = B_1 \cup B_2$ and $A = A_1 \cup A_2$ for $\alpha$ any vector in $V = V_1 \cup V_2$.

$$[[T\alpha]_{B'} = \left[\left(T_1 \cup T_2\right)\ \left(\alpha^1 \cup \alpha^2\right)\right]_{(B_1,\, B_2)}$$

$$= A\ [\alpha]_B = [A_1 \cup A_2]\ \ \left[\alpha^1 \cup \alpha^2\right]_{(B_1,\, B_2)}.$$

$$\left[U(T\alpha)\right]_{B'} = \left[U\ (T\alpha)\right]_{(B'_1,\, B'_2)}$$

$$= \left[\left(U_1 \cup U_2\right)(T_1 \cup T_2)(\alpha)\right]_{(B'_1,\, B'_2)}$$

$$= BA\ \left[\alpha\right]_{(B_1,\, B_2)}$$

and hence by the definition and uniqueness of the representing bimatrix we must have C = BA.

One can also verify this by carrying out the computation

$$\begin{aligned}
(UT)\ (\alpha_j) &= U\ (T(\alpha_j)) \\
&= U\ (T\ (\alpha_j^1, \alpha_j^2)) \\
&= U\left(\sum_{K=1}^{m} A_{Kj}^1\ \beta_K^1 \cup \sum_{K=1}^{m_1} A_{Kj}^2\ \beta_K^2\right) \\
&= \sum_{K=1}^{m} A_{Kj}^1 \left(U_1\ \beta_K^1\right) \cup \sum_{K=1}^{m} A_{Kj}^2 \left(U_2\ \beta_K^2\right) \\
&= \sum_{K=1}^{m} A_{Kj}^1 \sum_{i=1}^{p} B_{iK}^1\ \gamma_i^1 \cup \sum_{K=1}^{m_1} A_{Kj}^2 \sum_{i=1}^{p_1} B_{iK}^2\ \gamma_i^2\ .
\end{aligned}$$

We can just summarize the linear bitransformation in the following as a result or as a theorem.

*Result:* Let V, W and Z be finite dimensional bivector spaces over the field F; let T be a linear bitransformation from $V = V_1 \cup V_2$ into $W = W_1 \cup W_2$ and $U = U_1 \cup U_2$ a linear bitransformation of W into Z. If $B = (B_1, B_2)$. $B'_1 = (B'_1, B'_2)$, and $B'' = (B''_1, B''_2)$, are ordered basis for the spaces V, W and Z respectively; if A is the bimatrix of $T = T_1 \cup T_2$ relative to the pair $B = (B_1, B_2)$. $B'_1 = (B'_1, B'_2)$, and $B = B_1 \cup B_2$ is the bimatrix of U relative to the pair $B' = (B'_1, B'_2)$, $B'' = (B''_1, B''_2)$, respectively then the bimatrix of the composition UT relative to the pair $B = (B_1, B_2)$ and $B'' = (B''_1, B''_2)$, is the product bimatrix C = AB.

It is important to note that if T and U are linear bioperators on a space $V = V_1 \cup V_2$ and we are representing by a single ordered basis $B = (B_1, B_2)$ then by the above results we have the simple form

$$(UT)_{B = (B_1, B_2)} = [U]_{B = (B_1 B_2)} [T]_{B = (B_1 B_2)}.$$

Thus in this case the correspondence which $B = (B_1, B_2)$ determines between bioperators and bimatrices is not only a bivector space isomorphism but preserves products. A simple consequence of this is that the linear bioperator T is invertible if and only if

$$[T]_{B = (B_1, B_2)} = [T_1 \cup T_2]_{B = (B_1, B_2)} = [T_1]_{B_1} \cup [T_2]_{B_2}$$

is an invertible bimatrix. For the identity operator. I is represented by the identity bimatrix and thus $UT = TU = I_1 \cup I_2$ is equivalent to

$$[U]_{B = (B_1, B_2)} [T]_{B = (B_1, B_2)}$$
$$= [T]_B [U]_B = I_1 \cup I_2 = I$$

of course when the bimatrix T is invertible

$$[T^{-1}]_{B = (B_1, B_2)} = [T]^{-1}{}_{B = (B_1, B_2)}.$$

Now we will study the case when the ordered basis is changed. For the sake of simplicity we shall consider this question only for linear bioperators on a space $V = V_1 \cup V_2$ so that we can use a single pair of ordered basis i.e., $B = (B_1, B_2)$. The specific question is this. Let $T = T_1 \cup T_2$ be a linear bioperator on the finite dimensional space $V = V_1 \cup V_2$ and let $B = (B_1, B_2)$ and $B' = (B'_1, B'_2)$ be two ordered basis of $V = V_1 \cup V_2$

To find out the bimatrices of $[T]_{B=(B_1, B_2)}$ and $[T]_{B'=(B'_1, B'_2)}$.

By definition

$$
\begin{aligned}
[T\alpha]_{B=(B_1,B_2)} &= [T]_{B=(B_1,B_2)}[\alpha] \\
&= \left[(T_1 \cup T_2)(\alpha^1, \alpha^2)\right]_{B=(B_1,B_2)} \\
&= \left[(T_1 \cup T_2)_{B=(B_1,B_2)}(\alpha^1, \alpha^2)\right]_{B=(B_1,B_2)}
\end{aligned}
$$

$$\left[T_1(\alpha^1) \cup T_2(\alpha^2)\right]_{B=(B_1,B_2)}$$

$$\left[[T_1]_{B_1} \cup [T_2]_{B_2}\right]\left[\alpha^1, \alpha^2\right]_{(B_1,B_2)}.$$

i.e., $\left[T_1\alpha^1\right]_{B_1} \cup \left[T_2\alpha^2\right]_{B_2}$

$$
\begin{aligned}
&= \left([T_1]_B \cup [T_2]_{B_2}\right)\left\{\left[P_1\left(\alpha^1\right)_{B'_1}\right] \cup \left[P_2\left(\alpha^2\right)_{B'_2}\right]\right\} \\
&= P_1\left[T_1\alpha_1\right]_{B'_1} \cup P_2\left[T_2\alpha_2\right]_{B^1_2}.
\end{aligned}
$$

Using all these equations we get

$$[T]_{B=\{B_1,B_2\}} P[\alpha]_{B'=\{B'_1,B'_2\}} = P[T\alpha]_{B'=\{B'_1,B'_2\}}$$

i.e., $\left[T_1 \cup T_2\right]_{B=(B_1-B_2)}\left[P\left[\alpha^1 \cup \alpha^2\right]\right]_{B'=\{B'_1,B'_2\}}$

$$
\begin{aligned}
&= P\left[T_1\alpha^1 \cup T_2\alpha^2\right]_{[B'_1 B'_2]} \\
&= P\left[T_1\alpha^1\right]_{B'_1} \cup P\left[T_2\alpha^2\right]_{B'_2}.
\end{aligned}
$$

$$\left(P^{-1}[T]_{B=\{B_1,B_2\}} P\right)[\alpha]_{B'=\{B_1,B_2\}}$$

$$= \quad [T\alpha]_{B' = (B'_1, B'_2)}$$

$$= \quad [T_1\alpha^1 \cup T_2\alpha^2]_{B' = (B'_1\ B'_2)}$$

$$= \quad [T_1\alpha^1]_{B'_1} \cup [T_2\alpha^2]_{B'_2}\,.$$

So we have

$[T]_{B^1 = (B_1^1, B_2^1)}$

$$= \quad [T_1]_{B_1^1} \cup [T_2]_{B_2^1}$$

$$= \quad P^{-1}\left\{[T_1]_{B_1} \cup [T_2]_{B_2}\right\}P\,.$$

It important to observe that there is a unique linear bioperator $U = U_1 \cup U_2$ which maps $B = \{B_1, B_2\}$ into $B' = \{B'_1, B'_2\}$ defined by $U\alpha_j = \alpha_j^1$; j = 1, 2, …, n.

$$\text{i.e., } (U_1 \cup U_2)\left(\alpha_j^1 \cup \alpha_j^2\right) = \left(\alpha'_j\right)^1 \cup \left(\alpha'_j\right)^2$$

$$\text{i.e., } U_1\,\alpha_j^1 \cup U_2\alpha_j^2 = \left(\alpha'_j\right)^1 \cup \left(\alpha'_j\right)^2.$$

This operator U is invertible since it carries a basis for $V = V_1 \cup V_2$ on to a basis for $V = V_1 \cup V_2$. The matrix $P = P_1 \cup P_2$ above is precisely the matrix of the operator U in the ordered basis $B = (B_1, B_2)$. For $P = P_1 \cup P_2$ is defined by

$$\left(\alpha_j^1\right)^1 \cup \left(\alpha_j^1\right)^2 = \sum_{i=1}^{n} P_{ij}^1\alpha_j^1 \cup \sum_{i=1}^{n} P_{ij}^2\ \alpha_i^2$$

and since

$$U_1\ \alpha_j^1 = \left(\alpha_j^1\right)^1$$

and

$$U_2\ \alpha_j^2 = \left(\alpha_j^1\right)^2$$

this equation can be rewritten as

$$U\alpha_j = \sum_{i=1}^{m} P_{ij}\alpha_i$$

So $P = [U]_B$

$$\text{i.e., } P = P_1 \cup P_2 = [U_1 \cup U_2]_{B = (B_1, B_2)}$$

$$\text{i.e., } [P_1 \cup P_2] = \left[U_1\right]_{B_1} \cup \left[U_2\right]_{B_2}.$$

Suppose T = $T_1 \cup T_2$ is a linear bioperator on V = $V_1 \cup V_2$. If P = $[P_1 \cup P_2]$ = $\left[P_1^1,\ldots,P_n^1\right];\left[P_1^2,\ldots,P_m^2\right]$ is a n × n and m × m matrix (i.e., a mixed square bimatrix) with columns

$$P_j^1 = \left[\left(\alpha'_j\right)^1\right]_{B_1}, \; P_j^2 = \left[\left(\alpha'_j\right)^2\right]_{B_2}$$

then

$$\left[T\right]_{B'} = \left[T_1 \cup T_2\right]_{B' = (B'_2, B'_2)}$$

$$= \left[T_1\right]_{B'_1} \cup \left[T_2\right]_{B'_2}$$

$$= P^{-1} [T]_B P$$

$$= P^{-1} \left[T_1 \cup T_2\right]_{B=(B_1 B_2)} P.$$

$$= \left(P_1^{-1}[T_1]_{B_1} P_1\right) \cup \left(P_2^{-1} [T_2]_{B_2} P_2\right).$$

Alternatively if U is the invertible operator on V defined by

$$U\alpha_i = \alpha'_j \text{ i.e., } U_1\alpha_i^1 \cup U_2\alpha_i^2 = \left(\alpha'_j\right)^1 \cup \left(\alpha'_{j_1}\right)^2$$

$$j = 1, 2,\ldots, n$$

$$j_1 = 1, 2, \ldots, m.$$

Then

$$\left[T\right]_{B'=(B'_1, B'_2)} = [U]_B^{-1}[T]_B[U]_B$$

i.e., $\left[T_1\right]_{B'_1} \cup \left[T_2\right]_{B'_2}$

$$= \left[U_1\right]_{B_1}^{-1} [T_1]_{B_1} [U_1]_{B_1} \cup \left[U_2\right]_{B_2}^{-1} [T_2]_{B_2} [U_2]_{B_2}.$$

We illustrate this by the following example:

***Example 2.4.1:*** Let T be a linear bioperator on R × R ∪ {Space of all polynomial functions from R into R which have degree less than or equal to 3} = V

$$T: R^2 \cup V \to R^2 \cup V$$

where

$$T = T_1 \cup T_2$$

$$T_1 (x_1, x_2) = (x_1, 0)$$

and $T_2$ is the differentiation operator on V.

Let B = {$B_1$, $B_2$} be an ordered basis of $R^2 \cup V$ where $B_1 = \{\in_1, \in_2\}$ and $B_2 = \{b_1, b_2, b_3, b_4\}$. Suppose B' = $\left\{B_1^1, B_2\right\}$ where $B'_1$ is an ordered basis for $R^2$ consisting of vectors $\in_1^1 = (1, 1), \in_2^1 = (2, 1)$

Then
$$\in_1^1 = \in_1 + \in_2$$
$$\in_2^1 = 2\in_1 + \in_2$$

so

$$[T_1]_{B_1} = \begin{bmatrix} 1 & 0 \\ 0 & 0 \end{bmatrix}.$$

And the $P_1$ matrix

$$\begin{bmatrix} 1 & 2 \\ 1 & 1 \end{bmatrix}$$

$$P_1^{-1} = \begin{bmatrix} -1 & 2 \\ 1 & -1 \end{bmatrix}$$

$$T_{B^1} = P_1^{-1}\, T_{1B}\, P_1$$

$$= \begin{bmatrix} -1 & 2 \\ 1 & -1 \end{bmatrix} \begin{bmatrix} 1 & 0 \\ 0 & 0 \end{bmatrix} \begin{bmatrix} 1 & 2 \\ 1 & 2 \end{bmatrix}$$

$$= \begin{bmatrix} -1 & 2 \\ 1 & -1 \end{bmatrix} \begin{bmatrix} 1 & 2 \\ 0 & 0 \end{bmatrix}$$

$$= \begin{bmatrix} -1 & -2 \\ 1 & 2 \end{bmatrix}.$$

So
$$T_1 \in_1^1 = (1,\ 0) = -\in_1^1 + \in_2^1$$
$$T_1 \in_2^1 = (2, 0) = -2\in_1^1 + 2\in_2^1$$

Now for the other part V the ordered basis $B_2 = \{f_1, f_2, f_3, f_4\}$ is defined by $f_i(x) = x^{i-1}$. Let t be a real number and define
$g_i(x) = (x + t)^{i-1}$ i.e.

$g_1 = f_1$
$g_2 = t f_1 + f_2$
$g_3 = t^3 f_1 + 2t f_2 + f_3$
$g_4 = t^2 f_1 + 3t^2 f_2 + 3t f_3 + f_4$.

Since the matrix

$$P_2 = \begin{bmatrix} 1 & t & t^2 & t^3 \\ 0 & 1 & 2t & 3t^2 \\ 0 & 0 & 1 & 3t \\ 0 & 0 & 0 & 1 \end{bmatrix}$$

is easily seen to be invertible with

$$P_2^{-1} = \begin{bmatrix} 1 & -t & t^2 & -t^3 \\ 0 & 1 & -2t & 3t^2 \\ 0 & 0 & 1 & -3t \\ 0 & 0 & 0 & 1 \end{bmatrix}$$

it follows $B'_2 = \{g_1, g_2, g_3, g_4\}$ is an ordered basis for V.

$$\left[T_2\right]_{B_2=\{f_1,f_2,f_3,f_4\}}$$

$$= \begin{bmatrix} 0 & 1 & 0 & 0 \\ 0 & 0 & 2 & 0 \\ 0 & 0 & 0 & 3 \\ 0 & 0 & 0 & 0 \end{bmatrix}.$$

The matrix of $T_2$ in the ordered basis $B'_2$ is

$$P_2^{-1}\left[T_2\right]_{B_2} P_2 = \begin{bmatrix} 0 & 1 & 0 & 0 \\ 0 & 0 & 2 & 0 \\ 0 & 0 & 0 & 3 \\ 0 & 0 & 0 & 0 \end{bmatrix}.$$

Thus $\left[T\right]_{B'=[B'_1, B'_2]} = P_B^{-1}\left[T\right]_B P_B,\ \left(B = (B_1, B_2)\right)$

$= \left[T_1 \cup T_2\right]_{B'=[B'_1 B'_2]}$

$= \left[T_1\right]_{B'_1} \cup \left[T_2\right]_{B'_2}$

$$= \begin{bmatrix} -1 & -2 \\ 1 & 2 \end{bmatrix} \cup \begin{bmatrix} 0 & 1 & 0 & 0 \\ 0 & 0 & 2 & 0 \\ 0 & 0 & 0 & 3 \\ 0 & 0 & 0 & 0 \end{bmatrix}$$

$\left(P_1 \cup P_2\right)^{-1}_{B=(B_1 B_2)} \left[T\right]_{B=(B_1 B_2)} \left(P_1 \cup P_2\right)_B$

$= \left(P_1 \cup P_2\right)^{-1}_{B=(B_1 B_2)} \left[T_1 \cup T_2\right]_{B=(B_1 B_2)} \left(P_1 \cup P_2\right)_{B=(B_1 B_2)}$

$= \left[P_1\right]^{-1}_{B_1} \left[T_1\right]_{B_1} \left[P_1\right]_{B_1} \cup \left[P_2\right]^{-1}_{B_2} \left[T_2\right]_{B_2} \left[P_2\right]_{B_2}$

$$= \begin{bmatrix} -1 & 2 \\ 1 & -1 \end{bmatrix} \begin{bmatrix} 1 & 0 \\ 0 & 0 \end{bmatrix} \begin{bmatrix} 1 & 2 \\ 1 & 1 \end{bmatrix} \cup \begin{bmatrix} 1 & -t & t^2 & t^3 \\ 0 & 1 & -2t & 3t^2 \\ 0 & 0 & 1 & -3t \\ 0 & 0 & 0 & 1 \end{bmatrix}$$

$$= \begin{bmatrix} 0 & 1 & 0 & 0 \\ 0 & 0 & 2 & 0 \\ 0 & 0 & 0 & 3 \\ 0 & 0 & 0 & 0 \end{bmatrix} \cup \begin{bmatrix} 1 & t & t^2 & t^3 \\ 0 & 1 & 2t & 3t^2 \\ 0 & 0 & 1 & 3t \\ 0 & 0 & 0 & 1 \end{bmatrix}.$$

We just define when are two bimatrices similar. Let A = $A_1^{m \times m} \cup A_2^{n \times n}$ and B = $B_1^{m \times m} \cup B_2^{n \times n}$ be two mixed square bimatrix (or just square bimatrix) over to field F. We say B is similar to A over F if there is an invertible (n × n, m × m) bimatrix P i.e., P = $P_1^{m \times m} \cup P_2^{n \times n}$ such that B = $P^{-1}AP$ i.e., B = $B_1^{m\times m} \cup B_2^{n\times n} = P_1^{-1}A_1P_1 \cup P_2^{-1}A_2P_2$ .

Unlike in usual matrices in case of bimatrices we can define the notion of semi similar. Let A = $A_1^{n\times n} \cup A_2^{m\times m}$ and B = $B_1^{n\times n} \cup B_2^{m\times m}$ bimatrices if one of $A_1$ (or $A_2$) is similar with $B_1$ (or $B_2$) i.e.,

$$B = B_1 \cup B_2$$
$$= P_1^{-1}A_1P_1 \cup I_1A_2I_2$$

($A_2$ and $B_2$ are not similar)
or

$$B = B_1 \cup B_2$$
$$= A_1 \cup P_2^{-1}A_2P_2$$

(or $A_1$ and $B_1$ are not similar) then we say A and B semi similar bimatrices.

Now we see how the concept of similar bimatrices is helpful in transformation or relation between basis in bivector spaces. If V = $V_1 \cup V_2$ is a (n, m) dimensional bivector space over the field F and B = ($B_1$, $B_2$) and B' = ($B'_1$, $B'_2$) be two ordered basis for V, then for each linear bioperator T = $T_1 \cup T_2$ of V = $V_1 \cup V_2$ the bimatrix B = $B_1 \cup B_2$.

$$= [T]_{B'=(B'_1,B'_2)}$$
$$= (T_1 \cup T_2)_{B'=(B'_1 B'_2)}$$
$$= [T_1]_{B'_1} \cup [T_2]_{B'_2}$$

is similar to the matrix A = $A_1 \cup A_2$ =

$$= [T]_{B=(B_1,B_2)}$$
$$= (T_1 \cup T_2)_{B=(B_1B_2)}$$
$$= [T_1]_{B_1} \cup [T_2]_{B_2} .$$

The argument also goes in the other direction. Suppose A and B are (m × m, n × n) bimatrices and that B is similar to A, Let V be any (m, n) dimensional bivector space over F and let B = ($B_1$, $B_2$) be an ordered basis for V = $V_1 \cup V_2$. Let T = $T_1 \cup T_2$ be a linear bioperator on V = $V_1 \cup V_2$ represented in the basis B = ($B_1$, $B_2$ ) by A = $A_1 \cup A_2$. If

$$B = B_1 \cup B_2 = P^{-1}AP = P^{-1}A_1 \cup A_2P$$
$$= \left(P_1^{-1} \cup P_2^{-1}\right)(A_1 \cup A_2)(P_1 \cup P_2)$$
$$= P_1^{-1}A_1P_1 \cup P_2^{-1}A_2P_2 .$$

Let B' = ($B_1$, $B_2$) be an ordered basis for V = $V_1 \cup V_2$ obtained from B = ($B_1$, $B_2$) by $P_1 \cup P_2$ = P i.e.,

$$\alpha'_j = \left(\alpha_j^1\right)' \cup \left(\alpha_j^2\right)' = \sum_{i=1}^{m} P_{ij}^1\alpha_i^1 \cup \sum_{i=1}^{n} P_{ij}^2\alpha_i^2 .$$

Then the bimatrix T = $T_1 \cup T_2$ in the ordered basis B' = ($B'_1$ $B'_2$) will be B = $B_1 \cup B_2$. Thus the statement B = $B_1 \cup B_2$ is similar to A = $A_1 \cup A_2$ means that on each n-dimensional space over F the bimatrices A = $A_1 \cup A_2$ and B = $B_1 \cup B_2$ represent the same linear bitransformation in two (possibly) different ordered basis.

Note that each (m × m, n × n) bimatrix A = $A_1^{m\times m} \cup A_2^{n\times n}$ is similar to itself using P = $I_1^{m\times m} \cup I_2^{n\times n}$. If B is similar to A then A is similar to B for B = $B_1 \cup B_2$

$$= P^{-1}AP$$
$$= P^{-1}(A_1 \cup A_2)\, P$$
$$= P_1^{-1}A_1P_1 \cup P_2^{-1}A_2P_2 .$$

implies that

$$A = \left(P^{-1}\right)^{-1} B\, P^{-1}$$
$$= \left(P^{-1}\right)^{-1} \left(B_1 \cup B_2\right) P^{-1}$$

$$= \quad P_1B_1P_1^{-1} \cup P_2B_2P_2^{-1}.$$

So A is similar to B as bimatrix implies B is similar to A. If A, B and C are bimatrices and if A is similar B and B is similar to C then it is easily verified that C is similar to A.

## 2.5 Bieigen values and Bieigen vectors

The concept of eigen values or character values have been very extensively studied by several researchers. Now we have just introduced the notion of bimatrices and bivector spaces and linear bioperators. Using these notions we will define the new notions like character bivalues, bieigen vectors or bicharacteristic vectors etc. which may be termed also as eigen bivalues, characteristic bivectors or eigen bivectors both will mean the same.

Suppose T is a linear bioperator on an (n, m) dimensional bispace. If we could find an ordered basis B = $(B_1, B_2) = \left\{\left(\alpha_1^1,...,\alpha_n^1\right),\left(\alpha_1^2,...,\alpha_m^2\right)\right\}$ of $V = V_1 \cup V_2$ in which $T = T_1 \cup T_2$ are represented by a diagonal bimatrix D = $D_1 \cup D_2$ we would gain considerable information about T = $T_1 \cup T_2$. For instance simple number associated with T = $T_1 \cup T_2$ such as the rank or determinant of $T = T_1 \cup T_2$ could be determined with little more than a glance at the bimatrix $D = D_1 \cup D_2$. We could describe explicitly the range and null bispace of T. Since

$$[T]_{B=(B_1,B_2)} = D = D_1 \cup D_2$$

$$\text{i.e., } [T_1 \cup T_2]_{(B=B_1 \cup B_2)} = D_1 \cup D_2 .$$

$$\text{i.e. } [T_1]_{B_1} \cup [T_2]_{B_2} = D_1 \cup D_2$$

if and only if

$$T\alpha_i = C_i\, \alpha_i$$

$$\text{i.e., } (T_1 \cup T_2)\left(\alpha_{i_1}^1 \cup \alpha_{i_2}^2\right) = \left(C_i^1 \cup C_i^2\right)\left(\alpha_{i_1}^1 \cup \alpha_{i_2}^2\right)$$

$$T_1\alpha_i^1 = C_i^1\alpha_i^1$$

and

$$T_2\alpha^2_{i_2} = C^2_{i_2}\alpha^2_{i_2}, \ i = 1, 2, \ldots, n, \ i_2 = 1, 2, \ldots, m.$$

The range would be a subbivector space spanned by $\alpha_k$'s for which $C_k$'s are non zero and the null space would be spanned by reaming $\alpha_k$'s. Indeed it seems fair to say that if we knew a basis B = ($B_1$, $B_2$) and a diagonal bimatrix D = $D_1 \cup D_2$ such that

$$[T]_B = [T_1]_{B_1} \cup [T_2]_{B_2} = D = D_1 \cup D_2$$

we would answer readily any question about T = $T_1 \cup T_2$ which might arise.

Can such linear bioperator T = $T_1 \cup T_2$ be represented by a diagonal bimatrix in some ordered basis? If not for which bioperators T = $T_1 \cup T_2$ does such a basis exist? How to find the pair of basis for the bivector space V = $V_1 \cup V_2$ if there is one?

If no such basis exists what is the simplest type of bimatrices which we can represent the linear bioperator T = $T_1 \cup T_2$. To have a view we introduce certain new notions like characteristic bivalues characteristic bivector and characteristic bispace.

**DEFINITION 2.5.1:** *Let V = $V_1 \cup V_2$ be a bivector space over the field F and let T = $T_1 \cup T_2$ be a linear bioperator on V = $V_1 \cup V_2$.*

*A characteristic bivalue of T = $T_1 \cup T_2$ is a pair of scalars C = ($C_1$, $C_2$) $C_1$ $C_2$ in F such that there is a non zero bivector $\alpha = \alpha^1 \cup \alpha^2$ in V = $V_1 \cup V_2$ with $T_1 \alpha^1 = C_1 \alpha^1$ and $T_2 \alpha^2 = C_2 \alpha^2$.*

*If C = ($C_1$, $C_2$) is the characteristic bivalue of T = $T_1 \cup T_2$ then*

*(a) for any $\alpha = \alpha^1 \cup \alpha^2$ such that $T\alpha = C\alpha$ i.e., $T_1 \alpha^1 = C^1 \alpha^1$ and $T_2 \alpha^2 = C_2 \alpha^2$*

$$[(T_1 \cup T_2)] \ \alpha = C(\alpha^1 \cup \alpha^2) = C_1 \alpha^1 \cup C_2\alpha^2$$

*is called a characteristic bivector of T associated with the characteristic bivalue C.*

*(b) The collection of all $\alpha$ such that $T\alpha = C\alpha$ is called the characteristic bispace associated with C. Characteristic bivalues are often called as characteristic biroots, latent biroots, eigen bivalues, proper bivalues, or spectral bivalues.*

We would use in this book either the term characteristic bivalues or eigen bivalues only.

If T is any linear bioperator and C is any scalar, the set of bivectors $\alpha = \alpha^1 \cup \alpha^2$ such that $T\alpha = C\alpha$ is a subbispace of $V = V_1 \cup V_2$. It is the null bispace of the linear bitransformation $(T - CI) = (T_1 - C_1 I_1) \cup (T_2 - C_2 I_2)$ we call $C = C_1 \cup C_2$ a characteristic bivalue of $T = T_1 \cup T_2$ if this sub bispace is different from the zero subbispace i.e., if $(T_1 - C_1 I_1) \cup (T_2 - C_2 I_2)$ fails to be 1 : 1. If the underlying bispace V is finite dimensional, $(T - CI) = (T_1 - C_1 I_1) \cup (T_2 - C_2 I_2)$ fails to be 1 to 1 precisely when its determinant is different from 0.

**THEOREM 2.5.1:** *Let $T = T_1 \cup T_2$ is a linear bioperator on a finite dimensional bivector space $V = V_1 \cup V_2$ and let $C = C_1 \cup C_2$ be a scalar. The following are equivalent.*

- *(i) $C = C_1 \cup C_2$ is the characteristic bivalue of T.*
- *(ii) The operator $(T - CI) = (T_1 - C_1 I_1) \cup (T_2 - C_1 I_2)$ is singular (not invertible).*
- *(iii) $\det (T - CI) = (0, 0)$ i.e., $\det (T_1 - C_1 I_1) = 0$ and $\det (T_2 - C_2 I_2) = 0$.*

So we proceed on to define characteristic bivalues of a bimatrix $A = A_1 \cup A_2$ in F.

**DEFINITION 2.5.2:** *If $A = A_1^{m \times m} \cup A_2^{n \times n}$ be a square mixed bimatrix over the field F, a characteristic bivalue of A in F is a pair of scalar $C = (C_1, C_2)$ in F such that the bimatrix*

$$(A - CI) = (A_1 - C_1 I^{m \times m}) \cup (A_2 - C_2 I^{n \times n})$$

*is singular (not invertible) since C = $C_1 \cup C_2$ is a characteristic bivalue of A = $A_1 \cup A_2$ if and only if*

$$det\ (A - CI) = (0, 0)$$

*i.e.,* $det\ (A_1 - C_1 I^{m \times m}) \cup det\ (A_2 - C_2 I^{n \times n}) = (0, 0).$

*Or equivalently if and only if*

$$det\ (CI - A) = 0$$

$$i.e.,\ det\ \left(C_1 I_1^{m\times m} - A_1\right) = 0$$

*and*

$$det\ \left(C_2 I_2^{n\times n} - A_2\right) = 0,$$

*we form the bimatrix*

$$(XI - A) = \left(x_1\ I_1^{m\times m} - A_1\right) \cup \left(x_2\ I_2^{n\times n} - A_2\right)$$

*with bipolynomial entries, and consider the bipolynomial*

$$\begin{aligned} f &= f_1 \cup f_2 \\ &= det\ \left(x_1\ I_1^{m\times m} - A_1\right) \cup\ det\left(x_2\ I_2^{m\times m} - A_2\right) \\ &= (f_1\ (C_1), f_2\ (C_2)), \end{aligned}$$

*clearly the characteristic bivalue of A in F are just the scalars $C_1$, $C_2$ in F such that $f_1(C_1) = 0$ and $f_2\ (C_2) = 0$. For this reason $f = f_1 \cup f_2$ is called the characteristic bipolynomial of A = $A_1 \cup A_2$. It is important to note that $f = f_1 \cup f_2$ is a monic bipolynomial which has degree exactly (m, n).*

This is easily seen from the formula for the bideterminant of a bimatrix in terms of its entries.
Now we proceed on to prove the following result.

***Result:*** Similar bimatrices have the same characteristic bipolynomial.

*Proof:* If B = $P^{-1}$ AP where B = $B_1 \cup B_2$, A = $A_1 \cup A_2$ and P = $P_1 \cup P_2$ then B = $P^{-1}$ A P.

$B_1 \cup B_2 \quad = \quad P^{-1} (A_1 \cup A_2)\, P$

$\quad\quad\quad = \quad (P_1 \cup P_2)^{-1} (A_1 \cup A_2)\, P_1 \cup P_2$

$\quad\quad\quad = \quad P_1^{-1} A_1 P_1 \cup P_2^{-1} A_2 P_2 .$

$\det (xI - B) = (\det (x_1\ I_1^{m\times m} - B_1)\ \det (x_2 I_2^{n\times n} - B_2))$

$= \quad \det (P^{-1} (xI - A)\, P)$

$= \quad \det (P^{-1} (x\, ( I_1^{m\times m} \cup I_2^{n\times n} - A_1 \cup -A_2))P)$

$= \quad \det P^{-1} \det (x\ I_1^{m\times m} \cup I_2^{n\times n} - A_1 \cup -A_2)P)\ \det P$

$= \quad \det (x\ I_1^{m\times m} \cup I_2^{n\times n} - A_1 \cup -A_2)$

$= \quad \det ( x_1 I_1^{m\times m} - A_1) , \det ( x_2 I_2^{n\times n} - A_2) ,$

i.e.,

$$\det ( x_1 I_1^{m\times m} - B_1) = \det ( x_1 I_1^{m\times m} - A_1)$$

and

$$\det ( x_2 I_2^{m\times m} - B_2) = \det ( x_2 I_2^{m\times m} - A_2) .$$

This result enables us to define sensibly the characteristic bipolynomial of the bioperator $T = T_1 \cup T_2$ as the characteristic bipolynomial of any $(m \times m, n \times n)$ bimatrix which represents $T = T_1 \cup T_2$ in some ordered basis for $V = V_1 \cup V_2$. Just as for bimatrices the characteristic bivalues of $T = T_1 \cup T_2$ will be the biroots of the characteristic bipolynomial for T. In particular $T = T_1 \cup T_2$ cannot have more than (n, m) distinct characteristic bivalues. It is important to point out that $T = T_1 \cup T_2$ may not have any characteristic bivalues.

***Example 2.5.1:*** Let $T = T_1 \cup T_2$ be a linear bioperator on the vector space $V = V_1 \cup V_2$ of dimension (2, 3) on the standard basis.

$$T: V \to V$$
$$\text{i.e., } T_1 \cup T_2 : V_1 \cup V_2 \to V_1 \cup V_2$$
$$\text{i.e., } (T_1 : V_1 \to V_1) \cup (T_2 : V_2 \to V_2)$$

the related bimatrix is given by

$$A = \begin{bmatrix} 0 & -1 \\ 0 & 1 \end{bmatrix} \cup \begin{bmatrix} 3 & 1 & -1 \\ 2 & 2 & -1 \\ 2 & 2 & 0 \end{bmatrix}.$$

The characteristics bipolynomial is given by the bideterminant;

$$\begin{aligned} |\, xI - A \,| &= \left| xI_1^{2\times 2} - A_1 \right| \cup \left| xI_2^{3\times 3} - A_2 \right| \\ &= \begin{vmatrix} x & 1 \\ -1 & x \end{vmatrix} \cup \begin{vmatrix} x-3 & -1 & 1 \\ -2 & x-2 & 1 \\ -2 & -2 & x \end{vmatrix} \\ &= \{(x^2 + 1), (x^3 - 5x^2 + 8x - 4)\} \\ &= \{(x^2 + 1), (x - 1)(x - 2)^2\}. \end{aligned}$$

Now we face with a very odd situation. That is we have the bipolynomial which is such that one of the bipolynomial $x^2 + 1$ has no real roots i.e., $T_1$ has no characteristic values, where as $T_2$ has characteristic values 1 and 2 so we see $T = T_1 \cup T_2$ is a linear bioperator such that only one of the component operators has characteristic values and the other has no characteristic values. How to distinguish such cases we do this by the following ways.

**DEFINITION 2.5.3:** *Let $T = T_1 \cup T_2 : V \to V$ be a linear bioperator on the bivector space $V = V_1 \cup V_2$ i.e.,*

$$(T_1 : V_1 \to V_1) \cup (T_2 : V_2 \to V_2),$$

*if T is a linear bitransformation from V to V..*

*Using the standard basis for both the linear operators $T_1$ and $T_2$ we get the bimatrix $A = A_1 \cup A_2$. The characteristic bipolynomial is a pair given by (det $(xI_1 - A_1)$, det $(xI_2 - A_2)$). If both the polynomials in the bipolynomial has roots in the field over which the bivector space $V = V_1 \cup V_2$ is defined we say the linear bioperator $T = T_1 \cup T_2$ has characteristic bivalues.*

*If only one of $T_1$ or $T_2$ in the linear bioperator $T = T_1 \cup T_2$ has characteristic values we say T is a linear bioperator which has semi characteristic bivalues. If in the linear bioperator $T = T_1 \cup T_2$ both $T_1$ and $T_2$ has no characteristic values then we say the linear bioperator has no characteristic bivalues in F.*

We say it has characteristic bivalues in $\mathbb{C}$ (For $\mathbb{C}$ is the algebraically closed field and all polynomials are linearly reducible over $\mathbb{C}$, the field of complex numbers.)

If for the linear bioperator we have characteristic bivalues than we, as in case of matrices calculate the characteristic bivectors associated with the characteristic bivalues, we say the characteristic bivectors exists for these if characteristic bivalues exist. We say in case of semi characteristic bivalues the related bivectors as semi characteristic bivectors.

In the example we have just given the linear bioperator $T = T_1 \cup T_2$ is such that semi characteristic bivalues exists over the reals. So we have only semi characteristic bivectors associated with it is given by (1 0 2), (1 1 2) which are the only semi characteristic bivectors associated with it. We just denote it by {ϕ, (1 0 2), (1 1 2)}.

We now give an example of a linear bioperator, which has characteristic bivalues associated with it.

***Example 2.5.2:*** Let V be a bivector space give by $V = V_1 \cup V_2 = (R \times R) \cup R \times Q \times Q \times Q$ over the field Q.

Let $T = T_1 \cup T_2 : V \to V$ be a linear bioperator. Let the bimatrix associated with T relative to the standard basis be given by

$$A = A_1 \cup A_2 = \begin{bmatrix} 2 & 0 \\ 1 & -2 \end{bmatrix} \cup \begin{bmatrix} 1 & 1 & 0 & 0 \\ -1 & -1 & 0 & 0 \\ -2 & -2 & 2 & 1 \\ 1 & 1 & -1 & 0 \end{bmatrix}.$$

The characteristic bipolynomial associated with $T = T_1 \cup T_2$ (or for $A = A_1 \cup A_2$) is

$$\left(\det\left(xI_1^{2\times 2} - A_1\right) \cup \det\left(xI_2^{4\times 4} - A_2\right)\right) = ((x^2 - 4) \cup x^2(x - 1)^2).$$

Both the polynomials has real roots so $T = T_1 \cup T_2$ (or $A = A_1 \cup A_2$) has characteristics bivalues.

Now we will find the associated characteristic bivectors for the characteristic bivalues {(2, –2), (0, 0, 1, 1)}.

The characteristic bivalues are ({(0, k), (4k, k)}, {k, – k, 0, 0)}).

Now we have to study whether the associated characteristic bivectors are linearly independent or not or whether they generated the whole subbispace.

The concept of linear dependence, linear independence are discussed in the earlier section. For more about these notions please refer any book on linear algebra.

## 2.6 Bidiagonalizable linear bioperator and its properties

Now we proceed onto define the notions of bidiagonalizable, semi bidiagonalizable not bidiagonalizable for linear bioperator $T = T_1 \cup T_2$.

**DEFINITION 2.6.1:** *Let $T = T_1 \cup T_2$ be a linear bioperator on the finite-dimensional bivector space $V = V_1 \cup V_2$. We say T is bidiagonalizable if there is pair of basis $B = (B_1, B_2)$ for $V = V_1 \cup V_2$ each vector of which is the characteristic bivector of T.*

It is important to note that given any linear bioperator T from V to V, T in general need not always be a bidiagonalizable bioperator. We saw in the examples given in section 2.5 both the linear bioperators are not bidiagonalizable linear bioperators as in one case the characteristic bivalues did not exist and so we cannot imagine of characteristic bivectors. In the second case we saw the linear bioperator $T = T_1 \cup T_2$ has characteristic bivalues and in fact we could find the related characteristic bivectors but unfortunately those bivectors did not generate the bivector space V. So they will not be forming a basis for $V = V_1 \cup V_2$. Hence the linear bioperator in both examples given in section 2.5 are not bidiagonalizable. In fact in the first example the characteristic bivectors don't exist as the bivalues do not exist i.e., one of the characteristic bipolynomial has no solution over the field. In the second case characteristic bivector existed but did not form the basis for the bivector space $V = V_1 \cup V_2$.

*Note:* As all n-dimensional vector spaces V defined over the field F is such that

$$V \cong \underbrace{F \times \ldots \times F}_{n-\text{times}},$$

without loss of generality we denote it as (m, n) dimensional bivector space over F by

$$V \cong \underbrace{F \times \ldots \times F}_{m-\text{times}} \cup \underbrace{F \times \ldots \times F}_{n-\text{times}}.$$

***Example 2.6.1:*** Let $V = F \times F \cup F \times F \times F$ be a bivector space defined over the field F. Let $T = T_1 \cup T_2$ be a linear bioperator on $V = V_1 \cup V_2$. Let $A = A_1 \cup A_2$ be the associated bimatrix of T related to the standard basis of V.

$$A = \begin{bmatrix} 1 & 0 \\ 5 & 3 \end{bmatrix} \cup \begin{bmatrix} 2 & 0 & 0 \\ 9 & 1 & 0 \\ 0 & 0 & 3 \end{bmatrix}$$

det (x I – A )

$$= \quad (\det\left( xI_1^{2\times 2} - \begin{bmatrix} 1 & 0 \\ 5 & 3 \end{bmatrix}\right), \det\left( xI_2^{3\times 3} - \begin{bmatrix} 2 & 0 & 0 \\ 9 & 1 & 0 \\ 0 & 0 & 3 \end{bmatrix}\right)$$

$$= \quad \left(\det\begin{bmatrix} x-1 & 0 \\ -5 & x-3 \end{bmatrix}\right), \det\begin{bmatrix} x-2 & 0 & 0 \\ -9 & x-1 & 0 \\ 0 & 0 & x-3 \end{bmatrix}$$

$= \{(x-1)(x-3), (x-2)(x-1)(x-3)\}.$

Thus the related characteristic bivalues are ([1, 3], [1, 2, 3]) Now we find the relative characteristic bivectors associated with $T = T_1 \cup T_2$. The characteristic bivectors are given by [{(1, –5/2 ) (0, k)} {(0, 1, 0) (1, –9, 0) (0, 0, 1)}]. Clearly this forms a basis for the bivector space. Hence the linear bioperator $T = T_1 \cup T_2$ is bidiagonalizable.

Suppose that $T = T_1 \cup T_2$ is a bidiagonalizable linear bioperator. Let

$$\left\{\left(C_1^1, C_2^1, \ldots, C_k^1\right), \left(C_1^2, C_2^2, \ldots, C_{k'}^2\right)\right\}$$

be the distinct characteristic bivalues of the linear bioperator T.

Then there is an ordered basis $B = (B_1, B_2)$ in which T is represented by a diagonal bimatrix which has for its diagonal entries the scalar $C_i = \left(C_i^1, C_i^2\right)$, each repeated a certain number of times. If $C_i$ is repeated $\left(d_i^1, d_i^2\right)$ times then the bimatrix has the block form.

$$[T]_B = [T_1]_{B_1} \cup [T_2]_{B_2}$$

$$= \begin{bmatrix} C_1^1 I_{11} & 0 & \cdots & 0 \\ 0 & C_2^1 I_{21} & 0 & 0 \\ \vdots & & & \\ 0 & \cdots & 0 & C_k^1 I_{k1} \end{bmatrix} \cup \begin{bmatrix} C_1^2 I_{12} & 0 & \cdots & 0 \\ 0 & C_2^2 I_{22} & 0 & 0 \\ \vdots & & & \\ 0 & \cdots & 0 & C_k^2 I_{k2} \end{bmatrix}$$

where $I_{i_1} \cup I_{i_2}$ is the $\left(d_i^1 \times d_i^1, d_i^2 \times d_i^2\right)$ identity bimatrix. From the bimatrix we make the following two observations.

First the characteristic bipolynomial for T is the biproduct of linear factors

$f = f_1 \cup f_2$

$$= \left(\left(x - c_1^1\right)^{d_1^1} \cdots \left(x - c_k^1\right)^{d_{k'}^1}, (x - c_1^2)^{d_1^2} \cdots (x - c_{k'}^2)^{d_{k'}^2}\right).$$

If the field F is algebraically closed for example the field of complex numbers every polynomial over F can be factored, however if F is not algebraically closed we are citing a special property of T = $T_1 \cup T_2$ when we say that its characteristic bipolynomial has such a factorization.

The second is that $d_i = \left(d_i^1, d_i^2\right)$ are number of times in which $C_i = \left(C_i^1, C_i^2\right)$ is repeated as root of f = $f_1 \cup f_2$, equal to the dimension of the bispace of characteristic bivectors associated with the characteristic bivalues $C_i = \left(C_i^1, C_i^2\right)$. That is because the nullity of a diagonal bimatrix is equal to the number of zeros which it has on the main diagonals and the bimatrix

$$\left|T - C_i I\right|_B = \left|T - C_i^1 I^1\right|_{B_1} \cup \left|T - C_i^2 I^2\right|_{B_2}$$

has $\left(d_i^1, d_i^2\right)$ zeros on its main diagonal. (i associated with $d^1{}_i$ varies from 1 to k and I associated with $d_i{}^2$ varies from 1 to k').

The relation between dimension of the characteristic bispace and the multiplicity of it characteristic bivalues as a root of f does not seem interesting first, however it will help us with a simpler way of determining whether a given bioperator is bidiagonalizable. It is left as an exercise the reader to prove.

Suppose that T $\alpha$ = c$\alpha$ here $\alpha = (\alpha^1, \alpha^2)$. T = $T_1 \cup T_2$ linear bioperator. If f = $f_1 \cup f_2$ is any bipolynomial then prove f(T) $\alpha$ = f (c) $\alpha$ i.e., $f_1 (T_1) \alpha^1 \cup f_2 (T_2) \alpha^2 = f_1 (c_1) \alpha^1 \cup f_2 (c_2) \alpha^2$.

Now we proceed on to prove a nice result about, the linear bioperator and its relation to the subspace generated by the characteristic bivectors.

**LEMMA 2.6.1:** *Let T be a linear bioperator on the finite dimensional bivector space V = $V_1 \cup V_2$.*

*Let* $\left\{\left(C_1^1, C_2^1, \ldots, C_k^1\right), \left(C_1^2, C_2^2, \ldots, C_{k'}^2\right)\right\}$ *be the distinct characteristic bivalues of T = $T_1 \cup T_2$. Let $W_i = W_i^1 \cup W_i^2$* *be the bispace of characteristic bivectors associated with the characteristic bivalue $C_i$ = $C_i^1 \cup C_i^2$. If W = $W_1 \cup W_2$ then*

$$W = \left(W_1^1 + W_2^1 + \ldots + W_k^1\right) \cup \left(W_1^2 + W_2^2 + \ldots + W_{k'}^2\right)$$

*and*

$$\dim W = \dim (W_1 \cup W_2) =$$

$$\dim W_1^1 + \ldots + \dim W_k^1 \cup \dim W_1^2 + \ldots + \dim W_{k'}^2$$

*In fact if $B_i$ = $B_i^1 \cup B_i^2$ is an ordered basis of Wi = $W_i^1 \cup W_i^2$.*

*then B = $\left(\left\{B_1^1, \ldots, B_k^1\right\}, \left\{B_1^2, \ldots, B_{k'}^2\right\}\right)$ is an ordered basis for W.*

*Proof*: The bispace

$$W = W_1 \cup W_2 = \left(W_1^1 + \ldots + W_k^1\right) \cup \left(W_1^2 \cup \ldots \cup W_{k'}^2\right)$$

is a bisubspace spanned by all of the characteristic bivectors of T = $T_1 \cup T_2$. Usual, when one forms the sum W of bisubspaces $W_i = W_i^1 \cup W_i^2$ one expects that

$$\dim W = \dim (W_1 \cup W_2)$$
$$\dim W_1 \cup \dim W_2$$
$$< (\dim W_1^1 + \ldots + \dim W_k^1) \cup (\dim W_1^2 + \ldots + \dim W_{k'}^2),$$

because of their linear relation which may exist between bivectors in the various bispaces. This lemma or result however states that the characteristic bispaces associated with different characteristic bivalues are independent of one another.

Suppose that (for each i) we have a bivector $\beta_i$ in $W_i$ = $W_i^1 \cup W_i^2$ and assume that $\beta_1^1 + \ldots + \beta_k^1 = 0$ and $\beta_1^2 + \ldots + \beta_{k'}^2 = 0$ we shall show that $\beta_i^1 = 0$ and $\beta_i^2 = 0$ for each i. Let f = $f_1 \cup f_2$ be any bipolynomial since

$$T \beta_I = (T_1 \cup T_2)(\beta_i)$$

$$= \quad (T_1 \cup T_2) \left(\beta_i^1 \cup \beta_i^2\right)$$
$$= \quad T_1\ \beta_i^1 \cup T_2\beta_i^2 = C_i\beta_i$$
$$= \quad C_i^1\beta_i^1 \cup C_i^2\beta_i^2 .$$

The preceding result lets us that 0 = f (T) 0

$$= \quad f(T)\ \beta_1 + \ldots + f(T)\ \beta_k$$
$$= \quad \{f_1\ (T_1)\ \beta_1^1 + \ldots + f_1(T_1)\beta_k^1\}$$
$$\cup\ \{f_2\ (T_2)\ \beta_1^2 + \ldots + f_2(T_2)\beta_{k'}^2\}$$
$$= \quad f_1\ \left(C_1^1\right)\beta_1^1 + \ldots + f_1\left(C_k^1\right)\beta_k^1$$
$$\cup\ f_2\ \left(C_1^2\right)\beta_1^2 + \ldots + f_2\left(C_{k'}^2\right)\beta_{k'}^2$$

Choose polynomials $f_1^1 \ldots f_k^1$ and $f_1^2 \ldots f_{k'}^2$, such that

$$f_i^1\,(C_j^1) = \delta_{ij}^1 = \begin{cases} 1 & i = j \\ 0 & i \neq j \end{cases}$$

$$f_i^2\,(C_j^2) = \delta_{ij}^2 = \begin{cases} 0 & i = j \\ 0 & i \neq j \end{cases}$$

$$0 \quad = \quad f_i^1(T_1)\,0 \cup f_i^2\,(T_2)\,0$$
$$= \quad \sum_j \delta_{ij}^1\beta_j^1 \cup \sum_j \delta_{ij}^2\beta_j^2 .$$
$$= \quad \beta_i^1 \cup \beta_i^2 .$$

Now $B_i$ = $\left(B_i^1, B_i^2\right)$ be an ordered basis for $W_i$ = $W_i^1 \cup W_i^2$ and let B be the sequence B = $\left\{\left(B_1^1,\ldots,B_k^1\right),\left(B_1^2,\ldots,B_{k'}^2\right)\right\}$.

Then B spans the subspace W = $W^1 \cup W^2$

$$= \left(W_1^1 + \ldots + W_k^1\right) \cup \left(W_1^2 + \ldots + W_{k'}^2\right).$$

Also B is a linearly independent sequence if bivectors. For any linear combination of these bivectors in $B_i$ =

$\left(B_i^1, B_i^2\right)$, we have proved $B_i = \left(B_i^1, B_i^2\right)$ is linearly independent.

We now prove the most interesting theorem.

**THEOREM 2.6.1:** *Let $T = T_1 \cup T_2$ be a linear bioperator of a finite dimensional bivector space $V = V_1 \cup V_2$. Let*

$$\left\{\left(C_1^1, \ldots, C_k^1\right), \left(C_1^2, \ldots, C_{k'}^2\right)\right\}$$

*be distinct characteristic bivalues of $T = T_1 \cup T_2$ and let $W_i = W_i^1 \cup W_i^2$ be the null space of*

$$(T - C_i I) = (T_1 - C_i^1 I_1) \cup (T_2 - C_i^2 I_2).$$

*The following are equivalent*

*(i) T is bidiagonalizable*

*(ii) The characteristic bipolynomial for $T = T_1 \cup T_2$ is*

$f = f_1 \cup f_2$

$$= \left\{\left(x - C_1^1\right)^{d_1^1} \ldots \left(x - C_k^1\right)^{d_k^1}\right\} \cup \left\{\left(x - C_1^2\right)^{d_1^2} \ldots \left(x - C_{k'}^2\right)^{d_{k'}^2}\right\}$$

*and dim $W_i^1 = d_i^1$, dim $W_i^2 = d_i^2$, i = 1, 2,…, k and i = 1, 2 ,…, k'.*

*(iii)* $\left(dim W_1^1 + \ldots + dim W_k^1\right) \cup \left(dim W_1^2 + \ldots + dim W_{k'}^2\right)$

$= dim\left(V_1 \cup V_2\right)$

$= dim\ V.$

*Proof:* We have observed that T is bidiagonalizable implies there is a characteristic bipolynomial for $T = T_1 \cup T_2$ satisfying (ii) i.e., (i) $\Rightarrow$ (ii). Further if the characteristic bipolynomial $f = f_1 \cup f_2$ is the product of linear bifactors as in (ii) then $d_1^1 + \ldots + d_k^1 = \dim V_1$ and $d_1^2 + \ldots + d_{k'}^2 = \dim V_2$ where $V = V_1 \cup V_2$. For the sum of the $d_i$'s is the degree of the characteristic bipolynomial of that degree is dim V = (dim $V_1$, dim $V_2$). Hence (ii) implies (iii).

Suppose (iii) is true then

$$V = V_1 \cup V_2 = \left(W_1^1 + \ldots + W_k^1\right) \cup \left(W_1^2 + \ldots + W_{k'}^2\right)$$

i.e., characteristic bivectors of $T = T_1 \cup T_2$ span V.

It is left as an exercise for the reader to obtain the bimatrix analogue of the above theorem.

Several other related result can also be obtained in case of bimatrices which are analogous to matrices. We in this book restrain our selves only to a few of the results which have interested us.

Now we give a result analogous to minimal polynomial for a linear operator T. Let $T = T_1 \cup T_2$ be a linear bioperator on a finite dimensional bivector space $V = V_1 \cup V_2$ over the field F. The biminimal polynomial for the linear bioperator $T = T_1 \cup T_2$ is the unique monic generator of the biideal of polynomials over F which annihilate $T = T_1 \cup T_2$. The reader is requested to refer [ ] for the notions of biideals.

This can roughly be restated as if $T = T_1 \cup T_2$ is linear bioperator on $V = V_1 \cup V_2$ then $T_1$ can be realized as the linear operator on $V_1$ and $T_2$ can be realized as a linear operator on $V_2$ so the biminimal polynomial of $T = T_1 \cup T_2$ is the minimal polynomial of $T_1$ and minimal polynomial of $T_2$.

Several interesting results in this direction can be derived. However we just denote in some cases the biminimal polynomial will be the characteristic bipolynomial itself. Further unlike in case of linear operators we may in case of linear bioperators both minimal polynomials may coincide with the characteristic bipolynomial or one may coincide or both may not be coincident with the characteristic bipolynomial. We illustrate by one example how the biminimal polynomial looks in case of a linear bioperator.

***Example 2.6.2:*** Let us consider a bivector space $V = V_1 \cup V_2$ over reals R where both $V_1$ and $V_2$ are vector spaces of same dimension say 4.

i.e., $V = V_1 \cup V_2 = R \times R \times R \times R \cup R \times R \times R \times R$.

Let $T : V \to V$ be a linear bioperator, $T = T_1 \cup T_2 : V_1 \cup V_2 \to V_1 \cup V_2$. Let the associated bimatrix of T with respect to the standard basis is given by $A = A_1 \cup A_2$ where

$$A_1 = \begin{bmatrix} 0 & 1 & 0 & 1 \\ 1 & 0 & 1 & 0 \\ 0 & 1 & 0 & 1 \\ 1 & 0 & 1 & 0 \end{bmatrix}$$

and

$$A_2 = \begin{bmatrix} 1 & 1 & 0 & 0 \\ -1 & -1 & 0 & 0 \\ -2 & -2 & 2 & 1 \\ 1 & 1 & -1 & 0 \end{bmatrix}.$$

The characteristic bipolynomial associated with the linear bioperator is $f = x^2 (x^2 - 4) \cup x^2 (x - 1)^2 = f_1 \cup f_2$.

The biminimal polynomial for T is the unique monic generator of the biideal of bipolynomials over R which annihilate T is given by $p = x (x + 2) (x - 2) \cup x^2 (x - 1)^2 = p_1 \cup p_2$.

Now we define the subbispace of a bivector space V which is invariant under a linear bioperator $T = T_1 \cup T_2$.

**DEFINITION 2.6.2:** *Let $V = V_1 \cup V_2$ be a bivector space over F and let $T = T_1 \cup T_2$ be a linear bioperator on $V = V_1 \cup V_2$. If $W = W_1 \cup W_2$ is a subbivector space of V we say $W = W_1 \cup W_2$ is invariant under $T = T_1 \cup T_2$ if for each bivector $\alpha = \alpha^1 \cup \alpha^2$ in $W_1 \cup W_2$ ($\alpha^1 \in W_1, \alpha^2 \in W_2$) the vector $T(\alpha) = T_1(\alpha^1) \cup T_2(\alpha^2) \in W = W_1 \cup W_2$ i.e., if T(W) is contained in W.*

Now we just give an interesting lemma about invariant bivector subspaces of V relative to a linear bioperator T.

**LEMMA 2.6.2:** *Let $V = V_1 \cup V_2$ be a finite dimensional bivector space over the field F. Let $T = T_1 \cup T_2$ be a linear*

*bioperator on V = $V_1 \cup V_2$ such that the biminimal polynomial for T = $T_1 \cup T_2$ is a product of linear factors.*

$p = \left(x - C_1^1\right)^{r_1^1} \ldots \left(x - C_k^1\right)^{r_k^1} \cup \left(x - C_1^2\right)^{r_1^2} \ldots \left(x - C_{k'}^2\right)^{r_{k'}^2}$

*$C_i$ in F* $\left(C_i = \left(C_i^1, C_i^2\right)\right)$.

Let W be a proper subbispace of V i.e.,

$$W = W_1 \cup W_2 \ (W_1 \neq V_1 \text{ and } W_2 \neq V_2)$$

which is invariant under T. There exists a vector $\alpha = \alpha^1 \cup \alpha^2$ in $V = V_1 \cup V_2$ such that

(a) $\alpha$ is not in W i.e., $\alpha^1 \notin W_1$ and $\alpha^2 \notin W_2$

(b) $(T - CI)\ \alpha = (T_1 - C^1 I_1)\ \alpha^1 \cup (T_2 - C^2 I_2)\ \alpha^2$ is in $W = W_1 \cup W_2$ for some characteristic value C of the bioperator T.

Proof is left as an exercise for the reader.

**DEFINITION 2.6.3:** *Let W be an invariant bisubspace for T = $T_1 \cup T_2$ and let $\alpha = \alpha^1 \cup \alpha^2$ be a bivector in V. The T-biconductor i.e., (($T_1 \cup T_2$) biconductor) of $\alpha$ into W = $W_1 \cup W_2$ is the set $S_r(\alpha, W)$ which consists of all bipolynomials g = $g_1 \cup g_2$ over the scalar field F such that g (T) = $g_1(T_1) \cup g_2(T_2)$ is in $W_1 \cup W_2$).*

Let $g = g_1 \cup g_2$ be the T-biconductor of $\beta$ into W. Then g divides p the minimal bipolynomial for $T = T_1 \cup T_2$ . Since $\beta$ is not in $W = W_1 \cup W_2$ the bipolynomial g is not constant. Therefore

$$g = \left(x - C_1^1\right)^{e_1^1} \cdots \left(x - C_k^1\right)^{e_k^1} \cup \left(x - C_1^2\right)^{e_1^2} \cdots \left(x - C_{k'}^2\right)^{e_{k'}^2}$$

where at least one of the integers $e_i^1$ and $e_i^2$ are positive. Choose j, j[1] so that $e_j^1 > 0$ $e_{j'}^2 > 0$. Then

$$\left(x - C_j^1\right)/g_1 \text{ and} \left(x - C_{j'}^2\right)/g_2$$

where $$g = \left(x - C_j^1\right) h_1 \cup \left(x - C_{j'}^2\right) h_2 .$$

By definition of g the vector $\alpha = h(T)\ \beta = h_1\ (T_1)\ \beta^1 \cup h_2\ (T_2)\ \beta^2$, but by the definition of g the vector h (T)β cannot be in W i.e., $h_1\ (T_1)\ \beta^1 \notin W_1$ and $h_2\ (T_2)\ \beta^2 \notin W_2$.
But

$$\left(T_1 - C_i^1 I^1\right)\alpha^1 \cup \left(T_2 - C_i^2 I^2\right)\alpha^2$$

$$= \left(T_1 - C_i^1 I^1\right) h_1 (T_1)\ \beta^1 \cup \left(T_2 - C_i^2 I^2\right) h_2 (T_2)\ \beta^2$$

$$= g_1\ (T_1)\ \beta^1 \cup g_2\ (T_2)\ \beta^2 \text{ is in W.}$$

Thus we prove another interesting theorem for this we just recall the analogous definition of bitriangulable.

**DEFINITION 2.6.4:** *A linear bioperator $T = T_1 \cup T_2$ is represented by a bitriangular bimatrix i.e., there is an ordered basis for $V_1$ in which $T_1$ is represented by a triangular matrix and there is an ordered basis for $V_2$ in which $T_2$ is represented by a triangular matrix so that $T = T_1 \cup T_2$ is represented by a bitriangular bimatrix. (i.e., a bimatrix is a bitriangular bimatrix if both the matrices are triangular matries) i.e., we say T is bitriangulable.*

**THEOREM 2.6.2:** *Let $V = V_1 \cup V_2$ be a finite dimensional bivector space over the field F and let $T = T_1 \cup T_2$ be a linear bioperator on V. Then T is bitriangulable if and only if the minimal bipolynomial for T is the product of linear bipolynomials over F.*

*Proof:* Suppose that the minimal bipolynomial factors $p = p_1 \cup p_2$

$$= \left(x - C_1^1\right)^{r_1^1} \cdots \left(x - C_k^1\right)^{r_k^1} \cup \left(x - C_1^2\right)^{r_1^2} \cdots \left(x - C_{k'}^2\right)^{r_{k'}^2}$$

by repeated applications of the above lemma 2.6.2 we shall arrive at an ordered basis $B = \left\{\left(\alpha_1^1, ..., \alpha_k^1\right), \left(\alpha_1^2, ..., \alpha_{k'}^2\right)\right\}$ in which the bimatrix representing $T = T_1 \cup T_2$ is upper bitriangular

$$[T]_B = \begin{bmatrix} \alpha_{11}^1 & \alpha_{12}^1 & \cdots & \alpha_{1n}^1 \\ 0 & \alpha_{22}^1 & \cdots & \alpha_{2n}^1 \\ 0 & 0 & \alpha_{33}^1 & \alpha_{3n}^1 \\ \cdots & \cdots & 0 & \alpha_{nn}^1 \end{bmatrix} \cup \begin{bmatrix} \alpha_{11}^2 & \alpha_{12}^2 & \cdots & \alpha_{1m}^2 \\ 0 & \alpha_{22}^2 & \cdots & \alpha_{2m}^2 \\ 0 & 0 & \alpha_{33}^2 & \alpha_{3m}^1 \\ \cdots & \cdots & 0 & \alpha_{mm}^2 \end{bmatrix}$$

$$= \quad [T_1]_{B_1} \cup [T_2]_{B_2}.$$

The above identity merely implies

$$T_1\alpha_j^1 = a_{1j}^1\alpha_1^1 + \cdots + a_{jj}^1\alpha_j^1 \qquad 1 \le j \le n$$

$$\cup\ T_2\alpha_{j_1}^2 = a_{1j_1}^2\alpha_1^2 + \cdots + a_{j_1j_1}^2\alpha_{j_1}^2 \qquad 1 \le j_1 \le m$$

that is $T_1\alpha_j^1 \cup T_2\alpha_{j1}^2 = T\alpha_j$ is the bisubspace spanned by

$$\left\{\left(\alpha_1^1, ..., \alpha_j^1\right), \left(\alpha_1^2, ..., \alpha_{j_1}^2\right)\right\}.$$

To find $\left(\alpha_1^1, ..., \alpha_n^1\right), \left(\alpha_1^2, ..., \alpha_m^2\right)$ we start by applying lemma to the subspace W = $W_1 \cup W_2$ to obtain the vector $\alpha_1^1 \cup \alpha_1^2$, then apply the lemma to $W_1 = W_1^1 \cup W_1^2$, the space spanned by $\alpha_1$ = $\alpha_1^1 \cup \alpha_1^2$ and obtain $\alpha_2$ = $\alpha_1^1 \cup \alpha_1^2$. Next apply lemma to $W_2$ = $W_2^1 \cup W_2^2$ the space spanned by $\alpha_1$ and $\alpha_2$.

Continue in that way. One point deserves attention. After $\left(\alpha_1^1, ..., \alpha_i^1\right) \cup \left(\alpha_1^2, ..., \alpha_i^2\right)$ we have found it as the triangular type relations for j = 1, 2, … i, $j_1$ = 1, 2, …, i which ensure that the bisubspaces spanned by $\left(\alpha_1^1, ... \alpha_i^1\right), \left(\alpha_1^2, ... \alpha_{i'}^2\right)$ is invariant under T = $T_1 \cup T_2$ is bitriangulable, it is evident that the characteristic bipolynomial for T has the form f = $f_1 \cup f_2$

$$= \left(x - C_1^1\right)^{d_1^1} \cdots \left(x - C_k^1\right)^{d_k^1} \cup \left(x - C_1^2\right)^{d_1^2} \cdots \left(x - C_{k'}^2\right)^{d_{k'}^2}$$

$$C_i = C_i^1 \cup C_i^2 \in F.$$

Now we proceed on to define the notion of biprojection of $V = V_1 \cup V_2$. If $V = V_1 \cup V_2$ is a linear bioperator E on V such that $E^2$ = E i.e., E = $E_1 \cup E_2$ where $E_1$ on $V_1$ is such that $E_1^2 = E_1$ and $E_2$ on $V_2$ such that $E_2^2 = E_2$, then E is a biprojection on V. Let R be the range of E and let N be the null space of E.

1. The vector β is in the range R if and only if Eβ = β = (i.e., $E = E_1 \cup E_2$ and $E_1 \beta^1 = \beta^1$ and $E_2 \beta^2 = \beta^2$).
2. $V = R \oplus N$
   i.e., $V = V_1 \cup V_2 = R_1 \oplus N_1 \cup R_2 \oplus N_2$.
3. The unique expression for α as a sum of vector in R and N is $\alpha = E\alpha + (\alpha - E\alpha)$ $\alpha = \alpha^1 \cup \alpha^2 = E_1 \alpha^1 + (\alpha^1 - E_1 \alpha^1) \cup E_2 \alpha^2 + (\alpha^2 - E_2 \alpha^2)$.

We prove the following Theorem:

**THEOREM 2.6.3:** *If*

$$V = V_1 \cup V_2 = W_1^1 \oplus \ldots \oplus W_k^1 \cup W_1^2 \oplus \ldots \oplus W_{k'}^2$$

*then there exists k linear bioperators*

$$E_1 = E_1^1 \cup E_1^2, \ldots, E_k = E_k^1 \cup E_k^2 \text{ on } V_1 \cup V_2.$$

*such that*

(i) *each $E_i$ is a biprojection*
(ii) *$E_i \text{ o } E_j = 0 \quad i \neq j$*
(iii) *$I = E_1 + \ldots + E_k$*
(iv) *The range of $E_i$ is $W_i$ .*

*Proof:* The proof follows as in case of linear operators. For making or realizing each $E_i$ is the biprojection, i.e., $E_i^1 \cup E_i^2$ is projection on $W_i^1 \cup W_i^2$. We can easily verify

$$\begin{aligned} E_i \text{ o } E_j &= \left(E_i^1 \cup E_i^2\right) \text{o} \left(E_j^1 \cup E_j^2\right) \\ &= \left(E_i^1 \text{o} E_j^1\right) \cup E_i^2 \text{o} E_j^2 \\ &= 0 \cup 0 \end{aligned}$$

as $E_i^k$ and $E_j^k$ are projection on $W_i^k$ and $W_j^k$ respectively k = 1, 2.

Thus on similar lines with necessary modifications one can easily prove the following theorem on bidiagonalization.

**THEOREM 2.6.4:** *Let $T = T_1 \cup T_2$ be a linear bioperator on a finite dimensional bivector space $V = V_1 \cup V_2$. If T is diagonalizable and if $C_1, \ldots, C_k$ are distinct characteristic bivalues of T then there exists linear bioperators $E_1 \ldots E_k$ on V such that*

*(1)* $T = C_1 E_1 + \ldots + C_k E_k$
*i.e.,* $T_1 \cup T_2$
$= \left(C_1^1 E_1^1 + \ldots + C_k^1 E_k^1\right) \cup \left(C_1^2 E_1^2 + \ldots + C_{k'}^2 E_{k'}^2\right)$

*(2)* $I = E_1 + \ldots + E_k$
*i.e.,* $I = I^1 \cup I^2$
$= \left(E_1^1 + \ldots + E_k^1\right) \cup \left(E_1^2 + \ldots + E_{k'}^2\right)$

*(3)* $E_i E_j = \in 0 \ i \neq j$
*i.e.,* $E_i^1 E_j^1 = 0 \qquad E_i^2 \ E_j^2 = 0$ *if* $i \neq j$

*(4)* $E_i^2 = E_i$ *i.e.,* $E_i^1 \ E_i^1 = E_i^1$ *and* $E_i^2 \ E_i^2 = E_i^2$

*(5) The range of $E_i$ is the characteristic bispace of T associated with $C_i$* $C_i^1 \cup C_i^2 = (C_i^2, C_i^2)$

*Conversely if there exists (k, k') distinct scalars* $\left(C_1^1, \ldots, C_k^1\right), \left(C_1^2, \ldots, C_{k'}^2\right)$ *and (k, k') non-zero linear operators* $E_1^1, \ldots, E_k^1 \ E_1^2, \ldots, E_{k'}^2$ *which satisfy conditions (1) (2) and (3) then $T = T_1 \cup T_2$ is bidiagonalizable.* $C_1^1, \ldots, C_k^1$ *and* $C_1^2, \ldots, C_{k'}^2$ *are distinct characteristic bivalues of T and conditions (4) and (5) are also satisfied.*

The proof of this is a matter of routine as in case of linear operators with suitable operations.

Now we proceed on to define the notion of nilpotent bioperators or binilpotent operator and show that any linear bioperator can be realized as a sum of a diagonalizable part and a nilpotent part.

**DEFINITION 2.6.5:** *Let $V = V_1 \cup V_2$ be a bivector space over a field F. Let $T : V \to V$ be a linear bioperator on V i.e., $T = T_1 \cup T_2$ and $T_1 : V_1 \to V_1$ is a linear operator on $V_1$ and $T_2 : V_2 \to V_2$ is a linear operator on $V_2$. We say $T = T_1 \cup T_2$ is a nilpotent bioperator (binilpotent operator) on $V = V_1 \cup V_2$ if there is some pair of positive integers (r, s) such that $T_1^r = 0$ and $T_2^r = 0$ i.e., each $T_1$ and $T_2$ are nilpotent operators.*

It is to be noted that in all our discussion the bivector space $V = V_1 \cup V_2$ is such that $V_1$ may be of dimension n and $V_2$ may be of dimension m, m in general need not be equal to n. So even the bimatrix associated with the linear bioperator $T = T_1 \cup T_2$ may be only a square mixed bimatrix.

Now if we make a small deviation. All bivector spaces and bimatrices were defined only on the same field so we have no problem, but when the bivector spaces are defined over different fields or when the bimatrices are defined over different fields i.e., when the bimatrices are strong bimatrices or weak bimatrices we have to make the necessary modifications both in the definitions and in the working to get proper or analogous results. However while using strong bimatrices instead bimatrices the very concept of bidiagonalization or decomposition as submatrices becomes entirely different.

Now we give some problems / exercises in the last chapter which will make the reader to become at home with the notions of bivector spaces linear bioperators, bitransformation characteristic bivalues etc. We wish to state all these problems proposed are simple and can be solved by any reader who has followed these notions.

Chapter Three

# INTRODUCTION TO NEUTROSOPHIC BIMATRICES

Now in this chapter we define the notion of neutrosophic matrices and neutrosophic vector spaces and the related neutrosophic structures used to define it. This Chapter has three sections. The first section introduces some neutrosophic algebraic structures and neutrosophic matrices and vector spaces. In section two we for the first time define the notion of fuzzy bimatrices and neutrosophic bimatrices. Section three introduces the concept of neutrosophic bivector spaces.

In this book we assume all fields to be real fields of characteristic 0 all vector spaces are taken as real spaces over reals and we denote the indeterminacy by '$I$ ' as i will make a confusion as it denotes the imaginary value, viz. $i^2 = -1$ that is $\sqrt{-1} = i$. The indeterminacy $I$ is such that $I . I = I^2 = I$.

## 3.1 Introduction to some neutrosophic algebraic structures

In section we just recall some basic neutrosophic algebraic structures which are essential to define neutrosophic bimatrices and neutrosophic bivector spaces.

**DEFINITION 3.1.1:** *Let K be the field of reals. We call the field generated by $K \cup I$ to be the neutrosophic field for it*

*involves the indeterminacy factor in it. We define $I^2 = I$, $I + I = 2I$ i.e., $I + \ldots + I = nI$, and if $k \in K$ then $k.I = kI$, $0I = 0$. We denote the neutrosophic field by K(I) which is generated by $K \cup I$ that is $K(I) = \langle K \cup I \rangle$. ($\langle K \cup I \rangle$ denotes the smallest field generated by K and I.*

***Example 3.1.1:*** Let R be the field of reals. The neutrosophic field of reals is generated by $\langle R \cup I \rangle$ i.e. R(*I*) clearly $R \subset \langle R \cup I \rangle$.

***Example 3.1.2:*** Let Q be the field of rationals. The neutrosophic field of rationals is generated by $Q \cup I$ denoted by Q(*I*).

**DEFINITION 3.1.2:** *Let K(I) be a neutrosophic field we say K(I) is a prime neutrosophic field if K(I) has no proper subfield which is a neutrosophic field.*

***Example 3.1.3:*** Q(*I*) is a prime neutrosophic field where as R(*I*) is not a prime neutrosophic field for $Q(I) \subset R(I)$.

It is very important to note that all neutrosophic fields used in this book are of characteristic zero. Likewise we can define neutrosophic subfield.

**DEFINITION 3.1.3:** *Let K(I) be a neutrosophic field, $P \subset K(I)$ is a neutrosophic subfield of K if P itself is a neutrosophic field. K(I) will also be called as the extension neutrosophic field of the neutrosophic field P.*

Now we proceed on to define neutrosophic vector spaces, which can be defined over fields or neutrosophic fields. We can define two types of neutrosophic vector spaces one when it is a neutrosophic vector space over ordinary field other being neutrosophic vector space over neutrosophic fields. Both the vector spaces are distinctly different from each other for they also differ in the dimension. To this end we have to define neutrosophic group under addition.

**DEFINITION 3.1.4:** *We know Z is the abelian group under addition. Z(I) denote the additive abelian group generated by the set Z and I, Z(I) is called the neutrosophic abelian group under '+'.*

Thus to define basically a neutrosophic group under addition we need a group under addition. So we proceed on to define neutrosophic abelian group under addition.

Suppose G is an additive abelian group under '+'. G(*I*) = $\langle G \cup I \rangle$, additive group generated by G and *I*, G(*I*) is called the neutrosophic abelian group under '+'.

***Example 3.1.4:*** Let Q be the group under '+'; Q (*I*) = $\langle Q \cup I \rangle$ is the neutrosophic abelian group under addition; '+'.

***Example 3.1.5:*** R be the additive group of reals, R(*I*) = $\langle R \cup I \rangle$ is the neutrosophic group under addition.

***Example 3.1.6:*** $M_{n \times m}(I) = \{(a_{ij}) \mid a_{ij} \in Z(I)\}$ be the collection of all n × m matrices under '+'; with entries from the set Z(I) $M_{n \times m}(I)$ is a neutrosophic group under '+'.

Now we proceed on to define neutrosophic subgroup of a neutrosophic group.

**DEFINITION 3.1.5:** *Let G(I) be the neutrosophic group under addition. P ⊂ G(I) be a proper subset of G(I). P is said to be the neutrosophic subgroup of G(I) if P itself is a neutrosophic group i.e.* $P = \langle P_1 \cup I \rangle$ *where* $P_1$ *is an additive subgroup of G.*

***Example 3.1.7:*** Let Z(*I*) = $\langle Z \cup I \rangle$ be a neutrosophic group under '+'. $\langle 2Z \cup I \rangle$ = 2Z(*I*) is the neutrosophic subgroup of Z (*I*).

In fact Z(*I*) has infinitely many neutrosophic subgroups.

Now we proceed on to define the notion of neutrosophic quotient group.

**DEFINITION 3.1.6:** *Let G (I) = $\langle G \cup I \rangle$ be a neutrosophic abelian group under '+', suppose P(I) be a neutrosophic subgroup of G (I) then the neutrosophic quotient group*

$$\frac{G(I)}{P(I)} = \{a + P(I) \mid a \in G\,(I)\}\,.$$

***Example 3.1.8:*** Let Z ($I$ ) be a neutrosophic group under addition, Z the group of integers under addition, P =2Z($I$ ) is a neutrosophic subgroup of Z($I$ ), the neutrosophic quotient group

$$\frac{Z(I)}{P} = \{a + 2Z(I) \mid a \in Z(I)\}$$
$$= \{(2n+1) + (2n+1)\, I \mid n \in Z\}.$$

Clearly $\frac{Z(I)}{P}$ is a group. For P = 2Z ($I$ ) serves as the additive identity. Take a, b $\in \frac{Z(I)}{P}$. If a, b $\in$ Z($I$ ) \ P then two possibilities occur.

a + b is odd times $I$ or a + b is odd or a + b is even times $I$ or even if a + b is even or even times $I$ then a + b $\in$ P, if a + b is odd or odd times $I$,

$$\text{a + b} \in \frac{Z(I)}{P = 2Z(I)}\,.$$

It is easily verified that P acts as the identity and every element in $\frac{Z(I)}{P}$, that is for every

$$\text{a + 2Z } (I) \in \frac{Z(I)}{2Z(I)}$$

has a unique inverse.

Now we proceed on to define the notion of neutrosophic vector spaces over fields and then we define neutrosophic vector spaces over neutrosophic fields.

**DEFINITION 3.1.7:** *Let G(I) by an additive abelian neutrosophic group. K any field. If G(I) is a vector space over K then we call G(I) a neutrosophic vector space over K.*

Now we give the notion of strong neutrosophic vector space.

**DEFINITION 3.1.8:** *Let G(I) be a neutrosophic additive abelian group. K(I) be a neutrosophic field. If G(I) is a vector space over K(I) then we call G(I) the strong neutrosophic vector space.*

**THEOREM 3.1.1:** *All strong neutrosophic vector space over K(I) are a neutrosophic vector space over K; as $K \subset K(I)$.*

*Proof:* Follows directly by the very definitions.

Thus when we speak of neutrosophic spaces we mean either a neutrosophic vector space over K or a strong neutrosophic vector space over the neutrosophic field K(*I*). By basis we mean a linearly independent set in K(I) which spans the neutrosophic space.

Now we illustrate with an example.

***Example 3.1.9:*** Let $R(I) \times R(I) = V$ be an additive abelian neutrosophic group over the neutrosophic field *R(I).* Clearly V is a strong neutrosophic vector space over *R(I).* A basis of V are {(0,1), (1,0)}.

***Example 3.1.10:*** Let $V = R(I) \times R(I)$ be a neutrosophic abelian group under addition. V is a neutrosophic vector space over R. The neutrosophic basis of V are {(1,0), (0,1), (*I*,0), (0,*I*)}, which is a basis of the neutrosophic vector space V over R.

A study of these basis and its relations happens to be an interesting problems of research.

**DEFINITION 3.1.9:** *Let G(I) be a neutrosophic vector space over the field K. The number of elements in the neutrosophic basis is called the neutrosophic dimension of G(I).*

**DEFINITION 3.1.10:** *Let G(I) be a strong neutrosophic vector space over the neutrosophic field K(I). The number of elements in the strong neutrosophic basis is called the strong neutrosophic dimension of G(I).*

We denote the neutrosophic dimension of G(*I*) over K by $N_k$ (dim) of G (*I*) and that the strong neutrosophic dimension of G (*I*) by $SN_{K(I)}$ (dim) of G(*I*).

Now we define the notion of neutrosophic matrices.

**DEFINITION 3.1.11:** *Let* $M_{nxm} = \{(a_{ij}) \;/\; a_{ij} \in K(I)\}$*, where K(I), is a neutrosophic field. We call* $M_{nxm}$ *to be the neutrosophic rectangular matrix.*

***Example 3.1.11:*** Let Q(*I* ) = ⟨Q ∪ *I*⟩ be the neutrosophic field.

$$M_{4\times 3} = \begin{pmatrix} 0 & 1 & I \\ -2 & 4I & 0 \\ 1 & -I & 2 \\ 3I & 1 & 0 \end{pmatrix}$$

is the neutrosophic matrix, with entries from rationals and the indeterminacy *I*.

We define product of two neutrosophic matrices and the product is defined as follows:

Let

$$A = \begin{pmatrix} -1 & 2 & -I \\ 3 & I & 0 \end{pmatrix}_{2\times 3}$$

and

$$B = \begin{pmatrix} I & 1 & 2 & 4 \\ 1 & I & 0 & 2 \\ 5 & -2 & 3I & -I \end{pmatrix}_{3x4}$$

$$AB = \begin{bmatrix} -6I+2 & -1+4I & -2-3I & I \\ -4I & 3+I & 6 & 12+2I \end{bmatrix}_{2x4}.$$

(we use the fact $I^2 = I$).

Let $M_{n\times n} = \{(a_{ij}) \mid (a_{ij}) \in Q(I)\}$, $M_{n\times n}$ is a neutrosophic vector space over Q and a strong neutrosophic vector space over Q($I$ ).

Now we proceed onto define the notion of fuzzy integral neutrosophic matrices and operations on them, for more about these refer.

**DEFINITION 3.1.12:** *Let N = [0, 1] $\cup$ I where I is the indeterminacy. The m $\times$ n matrices $M_{m\times n} = \{(a_{ij}) / a_{ij} \in [0, 1] \cup I\}$ is called the fuzzy integral neutrosophic matrices. Clearly the class of m $\times$ n matrices is contained in the class of fuzzy integral neutrosophic matrices.*

***Example 3.1.12:*** Let

$$A = \begin{pmatrix} I & 0.1 & 0 \\ 0.9 & 1 & I \end{pmatrix},$$

A is a 2 $\times$ 3 integral fuzzy neutrosophic matrix.

We define operation on these matrices. An integral fuzzy neutrosophic row vector is a 1 $\times$ n integral fuzzy neutrosophic matrix. Similarly an integral fuzzy neutrosophic column vector is a m $\times$ 1 integral fuzzy neutrosophic matrix.

***Example 3.1.13:*** A = (0.1, 0.3, 1, 0, 0, 0.7, *I*, 0.002, 0.01, *I*, 0.12) is a integral row vector or a 1 × 11, integral fuzzy neutrosophic matrix.

***Example 3.1.14:*** B = (1, 0.2, 0.111, *I*, *I*, 0.32, 0.001, *I*, 0, 1)$^T$ is an integral neutrosophic column vector or B is a 9 × 1 integral fuzzy neutrosophic matrix.

We would be using the concept of fuzzy neutrosophic column or row vector in our study.

**DEFINITION 3.1.13:** *Let P = ($p_{ij}$) be a m × n integral fuzzy neutrosophic matrix and Q = ($q_{ij}$) be a n × p integral fuzzy neutrosophic matrix. The composition map P • Q is defined by R = ($r_{ij}$) which is a m × p matrix where* $r_{ij} = \max_k \min(p_{ik} q_{kj})$ *with the assumption max($p_{ij}$, I) = I and min($p_{ij}$, I) = I where* $p_{ij} \in [0, 1]$. *min (0, I) = 0 and max(1, I) = 1.*

***Example 3.1.15:*** Let

$$P = \begin{bmatrix} 0.3 & I & 1 \\ 0 & 0.9 & 0.2 \\ 0.7 & 0 & 0.4 \end{bmatrix}, \quad Q = (0.1, I, 0)^T$$

be two integral fuzzy neutrosophic matrices.

$$P \bullet Q = \begin{bmatrix} 0.3 & I & 1 \\ 0 & 0.9 & 0.2 \\ 0.7 & 0 & 0.4 \end{bmatrix} \bullet \begin{bmatrix} 0.1 \\ I \\ 0 \end{bmatrix} = (I, I, 0.1).$$

***Example 3.1.16:*** Let

$$P = \begin{bmatrix} 0 & I \\ 0.3 & 1 \\ 0.8 & 0.4 \end{bmatrix}$$

and

$$Q = \begin{bmatrix} 0.1 & 0.2 & 1 & 0 & I \\ 0 & 0.9 & 0.2 & 1 & 0 \end{bmatrix}.$$

One can define the max-min operation for any pair of integral fuzzy neutrosophic matrices with compatible operation.

Now we proceed onto define the notion of fuzzy neutrosophic matrices. Let $N_s = [0, 1] \cup nI / n \in (0, 1]\}$; we call the set $N_s$ to be the fuzzy neutrosophic set.

**DEFINITION 3.1.14:** *Let $N_s$ be the fuzzy neutrosophic set. $M_{n\times m} = \{(a_{ij}) / a_{ij} \in N_s\}$ we call the matrices with entries from $N_s$ to be the fuzzy neutrosophic matrices.*

***Example 3.1.17:*** Let $N_s = [0,1] \cup \{nI/ n \in (0,1]\}$ be the set

$$P = \begin{bmatrix} 0 & 0.2I & 0.31 & I \\ I & 0.01 & 0.7I & 0 \\ 0.31I & 0.53I & 1 & 0.1 \end{bmatrix}$$

P is a $3 \times 4$ fuzzy neutrosophic matrix.

***Example 3.1.18:*** Let $N_s = [0, 1] \cup \{nI / n \in (0, 1]\}$ be the fuzzy neutrosophic matrix. A = [0, 0.12*I*, *I*, 1, 0.31] is the fuzzy neutrosophic row vector:

$$B = \begin{bmatrix} 0.5I \\ 0.11 \\ I \\ 0 \\ -1 \end{bmatrix}$$

is the fuzzy neutrosophic column vector.

Now we proceed on to define operations on these fuzzy neutrosophic matrices.

Let $M = (m_{ij})$ and $N = (n_{ij})$ be two $m \times n$ and $n \times p$ fuzzy neutrosophic matrices. $M \bullet N = R = (r_{ij})$ where the entries in the fuzzy neutrosophic matrices are fuzzy indeterminates i.e. the indeterminates have degrees from 0 to 1 i.e. even if some factor is an indeterminate we try to give it a degree to which it is indeterminate for instance $0.9I$ denotes the indeterminacy rate; it is high where as $0.01I$ denotes the low indeterminacy rate. Thus neutrosophic matrices have only the notion of degrees of indeterminacy. Any other type of operations can be defined on the neutrosophic matrices and fuzzy neutrosophic matrices. The notion of these matrices have been used to define neutrosophic relational equations and fuzzy neutrosophic relational equations.

## 3.2. Neutrosophic bimatrices

Here for the first time we define the notion of neutrosophic bimatrix and illustrate them with examples. Also we define fuzzy neutrosophic matrices.

**DEFINITION 3.2.1:** *Let $A = A_1 \cup A_2$ where $A_1$ and $A_2$ are two distinct neutrosophic matrices with entries from a neutrosophic field. Then $A = A_1 \cup A_2$ is called the neutrosophic bimatrix.*

*It is important to note the following:*

1. *If both $A_1$and $A_2$ are neutrosophic matrices we call A a neutrosophic bimatrix.*
2. *If only one of $A_1$ or $A_2$ is a neutrosophic matrix and other is not a neutrosophic matrix then we all $A = A_1 \cup A_2$ as the semi neutrosophic bimatrix. (It is clear all neutrosophic bimatrices are trivially semi neutrosophic bimatrices).*

*It both $A_1$ and $A_2$ are $m \times n$ neutrosophic matrices then we call $A = A_1 \cup A_2$ a $m \times n$ neutrosophic bimatrix or a rectangular neutrosophic bimatrix.*

*If $A = A_1 \cup A_2$ be such that $A_1$ and $A_2$ are both $n \times n$ neutrosophic matrices then we call $A = A_1 \cup A_2$ a square or a $n \times n$ neutrosophic bimatrix. If in the neutrosophic bimatrix $A = A_1 \cup A_2$ both $A_1$ and $A_2$ are square matrices but of different order say $A_1$ is a $n \times n$ matrix and $A_2$ a $s \times s$ matrix then we call $A = A_1 \cup A_2$ a mixed neutrosophic square bimatrix. (Similarly one can define mixed square semi neutrosophic bimatrix).*

*Likewise in $A = A_1 \cup A_2$ if both $A_1$ and $A_2$ are rectangular matrices say $A_1$ is a $m \times n$ matrix and $A_2$ is a $p \times q$ matrix then we call $A = A_1 \cup A_2$ a mixed neutrosophic rectangular bimatrix. (If $A = A_1 \cup A_2$ is a semi neutrosophic bimatrix then we call A the mixed rectangular semi neutrosophic bimatrix).*

Just for the sake of clarity we give some illustration.

**Notation:** We denote a neutrosophic bimatrix by $A_N = A_1 \cup A_2$.

***Example 3.2.1:*** Let

$$A_N = \begin{bmatrix} 0 & I & 0 \\ 1 & 2 & -1 \\ 3 & 2 & I \end{bmatrix} \cup \begin{bmatrix} 2 & I & 1 \\ I & 0 & I \\ 1 & 1 & 2 \end{bmatrix}$$

$A_N$ is the $3 \times 3$ square neutrosophic bimatrix.

***Example 3.2.2:*** Let

$$A_N = \begin{bmatrix} 2 & 0 & I \\ 4 & I & 1 \\ 1 & 1 & 2 \end{bmatrix} \cup \begin{bmatrix} 3 & I & 0 & 1 & 5 \\ 0 & 0 & I & 3 & 1 \\ I & 0 & 0 & I & 2 \\ 1 & 3 & 3 & 5 & 4 \\ 2 & 1 & 3 & 0 & I \end{bmatrix}$$

$A_N$ is a mixed square neutrosophic bimatrix.

***Example 3.2.3:*** Let

$$A_N = \begin{bmatrix} 3 & 1 & 1 & 1 & I \\ I & 0 & 2 & 3 & 4 \end{bmatrix} \cup \begin{bmatrix} I & 2 & 0 & I \\ 3 & 1 & 2 & 1 \\ 4 & 1 & 0 & 0 \\ 3 & 3 & 1 & 1 \\ 1 & I & 0 & I \end{bmatrix}$$

$A_N$ is a mixed rectangular neutrosophic bimatrix. We can denote

$$A_N = A_1 \cup A_2 = A_1^{2\times5} \cup A_2^{5\times4}.$$

***Example 3.2.4:*** Let

$$A_N = \begin{bmatrix} 3 & 1 \\ 1 & 2 \\ I & 0 \\ 3 & I \end{bmatrix} \cup \begin{bmatrix} 3 & 3 \\ 1 & I \\ 4 & I \\ -I & 0 \end{bmatrix}$$

$A_N$ is $4 \times 2$ rectangular neutrosophic bimatrix.

***Example 3.2.5*:** Let

$$A_N = \begin{bmatrix} 3 & 1 & 1 \\ 2 & 2 & 2 \end{bmatrix} \cup \begin{bmatrix} -I & 1 & 2 \\ 0 & I & 3 \end{bmatrix}$$

$A_N$ is a rectangular semi neutrosophic bimatrix for $A = A_1 \cup A_2$ with

$$A_1 = \begin{bmatrix} 3 & 1 & 1 \\ 2 & 2 & 2 \end{bmatrix}$$

is not a neutrosophic matrix, only

$$A_2 = \begin{bmatrix} -I & 1 & 2 \\ 0 & I & 3 \end{bmatrix}$$

is a neutrosophic matrix.

***Example 3.2.6:*** Let

$$A_N = \begin{bmatrix} 3 & 1 & 1 & 1 \\ 0 & I & 1 & 2 \\ 0 & 0 & 0 & 3 \\ I & 1 & 1 & 1 \end{bmatrix} \cup \begin{bmatrix} 0 & 2 & 2 & 2 \\ 1 & 0 & 0 & 0 \\ 2 & 0 & 0 & 1 \\ 5 & 0 & -1 & 2 \end{bmatrix}$$

$A_N$ is a square semi neutrosophic bimatrix.

***Example 3.2.7:*** Let

$$A_N = \begin{bmatrix} 1 & 1 & 1 & 1 \\ 0 & 0 & 0 & 0 \end{bmatrix} \cup \begin{bmatrix} 0 \\ I \\ 1 \\ 2 \\ 3 \end{bmatrix}.$$

$A_N$ is a rectangular mixed semi neutrosophic bimatrix.

Thus as in case of bimatrices we may have square, mixed square, rectangular or mixed rectangular neutrosophic (semi neutrosophic) bimatrices.

Now we can also define the neutrosophic bimatrices or semi neutrosophic bimatrices over different fields. When both $A_1$ and $A_2$ in the bimatrix $A = A_1 \cup A_2$ take its values from the same neutrosophic field K we call it a neutrosophic (semi neutrosophic) bimatrix. If in the bimatrix $A = A_1 \cup A_2$, $A_1$ is defined over a neutrosophic field F and $A_2$ over some other neutrosophic field $F^1$ then we call the neutrosophic bimatrix as strong neutrosophic bimatrix. If on the other hand the neutrosophic matrix $A = A_1 \cup A_2$ is such that $A_1$ takes entries from the neutrosophic field F and $A_2$ takes its entries from a proper subfield of a neutrosophic field then we call A the weak neutrosophic bimatrix. All properties of bimatrices can be carried on to neutrosophic bimatrices and semi neutrosophic bimatrices.

Now we proceed on to define fuzzy bimatrix, fuzzy neutrosophic bimatrix, semi-fuzzy bimatrix, semi fuzzy neutrosophic bimatrix and illustrate them with examples.

**DEFINITION 3.2.2:** *Let $A = A_1 \cup A_2$ where $A_1$ and $A_2$ are two distinct fuzzy matrices with entries from the interval [0, 1]. Then $A = A_1 \cup A_2$ is called the fuzzy bimatrix.*

*It is important to note the following:*

1. *If both $A_1$and $A_2$ are fuzzy matrices we call A a fuzzy bimatrix.*
2. *If only one of $A_1$ or $A_2$ is a fuzzy matrix and other is not a fuzzy matrix then we all $A = A_1 \cup A_2$ as the semi fuzzy bimatrix. (It is clear all fuzzy matrices are trivially semi fuzzy matrices).*

*It both $A_1$ and $A_2$ are $m \times n$ fuzzy matrices then we call $A = A_1 \cup A_2$ a $m \times n$ fuzzy bimatrix or a rectangular fuzzy bimatrix.*

*If $A = A_1 \cup A_2$ is such that $A_1$ and $A_2$ are both $n \times n$ fuzzy matrices then we call $A = A_1 \cup A_2$ a square or a $n \times n$ fuzzy bimatrix. If in the fuzzy bimatrix $A = A_1 \cup A_2$ both $A_1$ and $A_2$ are square matrices but of different order say $A_1$ is a $n \times n$ matrix and $A_2$ a $s \times s$ matrix then we call $A = A_1 \cup A_2$ a mixed fuzzy square bimatrix. (Similarly one can define mixed square semi fuzzy bimatrix).*

*Likewise in $A = A_1 \cup A_2$ if both $A_1$ and $A_2$ are rectangular matrices say $A_1$ is a $m \times n$ matrix and $A_2$ is a $p \times q$ matrix then we call $A = A_1 \cup A_2$ a mixed fuzzy rectangular bimatrix. (If $A = A_1 \cup A_2$ is a semi fuzzy bimatrix then we call A the mixed rectangular semi fuzzy bimatrix).*

Just for the sake of clarity we give some illustration.

**Notation:** We denote a fuzzy bimatrix by $A_F = A_1 \cup A_2$.

***Example 3.2.8:*** Let

$$A_F = \begin{bmatrix} 0 & .1 & 0 \\ .1 & .2 & .1 \\ .3 & .2 & .1 \end{bmatrix} \cup \begin{bmatrix} .2 & .1 & .1 \\ .1 & 0 & .1 \\ .2 & .1 & .2 \end{bmatrix}$$

$A_F$ is the 3 × 3 square fuzzy bimatrix.

***Example 3.2.9:*** Let

$$A_F = \begin{bmatrix} .2 & 0 & 1 \\ .4 & .2 & 1 \\ .3 & 1 & .2 \end{bmatrix} \cup \begin{bmatrix} .3 & 1 & 0 & .4 & .5 \\ 0 & 0 & 1 & .8 & .2 \\ 1 & 0 & 0 & .1 & .2 \\ .1 & .3 & .3 & .5 & .4 \\ .2 & .1 & .3 & 0 & 1 \end{bmatrix}$$

$A_F$ is a mixed square fuzzy bimatrix.

***Example 3.2.10:*** Let

$$A_F = \begin{bmatrix} .3 & 1 & .5 & 1 & .9 \\ .6 & 0 & .2 & .3 & .4 \end{bmatrix} \cup \begin{bmatrix} 1 & .2 & 0 & 0 \\ .3 & 1 & .2 & 1 \\ .4 & 1 & 0 & 0 \\ .3 & .3 & .2 & 1 \\ 1 & .5 & .7 & .6 \end{bmatrix}$$

$A_F$ is a mixed rectangular fuzzy bimatrix. We can denote

$$A_F = A_1 \cup A_2 = A_1^{2\times 5} \cup A_2^{5\times 4}.$$

***Example 3.2.11:*** Let

$$A_F = \begin{bmatrix} .3 & 1 \\ 1 & .2 \\ .5 & 0 \\ .3 & .6 \end{bmatrix} \cup \begin{bmatrix} .3 & .7 \\ 1 & 1 \\ .4 & 1 \\ .2 & 0 \end{bmatrix}$$

$A_F$ is $4 \times 2$ rectangular fuzzy bimatrix.

***Example 3.2.12*****:** Let

$$A_F = \begin{bmatrix} 3 & 1 & 1 \\ 2 & 2 & 2 \end{bmatrix} \cup \begin{bmatrix} .5 & .7 & .2 \\ 0 & .1 & .3 \end{bmatrix}$$

$A_F$ is a rectangular semi fuzzy bimatrix for

$$A = A_1 \cup A_2$$

with

$$A_1 = \begin{bmatrix} 3 & 1 & 1 \\ 2 & 2 & 2 \end{bmatrix}$$

is not a fuzzy matrix, only

$$A_2 = \begin{bmatrix} .5 & .7 & .2 \\ 0 & .1 & .3 \end{bmatrix}$$

is a fuzzy matrix.

***Example 3.2.13:*** Let

$$A_F = \begin{bmatrix} .3 & 1 & 1 & 1 \\ 0 & 0 & .1 & .2 \\ 0 & 0 & 0 & .3 \\ .3 & 1 & 1 & 1 \end{bmatrix} \cup \begin{bmatrix} 0 & 2 & 2 & 2 \\ 1 & 0 & 0 & 0 \\ 2 & 0 & 0 & 1 \\ 5 & 0 & -1 & 2 \end{bmatrix}$$

$A_N$ is a square semi fuzzy bimatrix.

***Example 3.2.14:*** Let

$$A_F = \begin{bmatrix} 1 & 1 & 1 & 1 \\ 0 & 0 & 0 & 0 \end{bmatrix} \cup \begin{bmatrix} 0 \\ 1 \\ 1 \\ 2 \\ 3 \end{bmatrix}.$$

$A_F$ is a rectangular mixed semi fuzzy bimatrix.

Thus as in case of bimatrices we may have square, mixed square, rectangular or mixed rectangular fuzzy (semi fuzzy) bimatrices.

Now we proceed on to define fuzzy integral neutrosophic bimatrix.

**DEFINITION 3.2.3:** *Let $A_{FN} = A_1 \cup A_2$ where $A_1$ and $A_2$ are distinct integral fuzzy neutrosophic matrices. Then $A_{FN}$ is called the integral fuzzy neutrosophic bimatrix. If both $A_1$ and $A_2$ are $m \times m$ distinct integral fuzzy neutrosophic matrix then $A_{FN} = A_1 \cup A_2$ is called the square integral fuzzy neutrosophic bimatrix.*

As in case of neutrosophic bimatrices we can define rectangular integral fuzzy neutrosophic bimatrix, mixed square integral fuzzy neutrosophic matrix and so on.

If in $A_{FN} = A_1 \cup A_2$ one of $A_1$ or $A_2$ is a fuzzy neutrosophic matrix and the other is just a fuzzy matrix or a neutrosophic matrix we call $A_{FN}$ the semi integral fuzzy neutrosophic bimatrix.

Now we will illustrate them with examples.

***Examples 3.2.15:*** Let $A_{FN} = A_1 \cup A_2$ where

$$A_{FN} = \begin{bmatrix} 0 & I & .3 \\ .2I & .4 & 1 \\ 0 & .3 & -.6 \end{bmatrix} \cup \begin{bmatrix} 1 & I & 0 \\ I & 1 & .8 \\ .6 & 1 & .7I \end{bmatrix}$$

then $A_{FN}$ is a square fuzzy neutrosophic bimatrix.

***Example 3.2.16:*** Consider $A_{FN} = A_1 \cup A_2$ where

$$A_1 = \begin{bmatrix} 0 & .2I & 1 \\ I & .7 & 0 \\ 1 & 0 & .1 \end{bmatrix}$$

and

$$A_2 = \begin{bmatrix} .1 & 0 \\ .2 & I \\ 1 & .2I \end{bmatrix}.$$

Clearly $A_{FN}$ is a mixed fuzzy neutrosophic bimatrix.

***Example 3.2.17:*** Let $A_{FN} = A_1 \cup A_2$ where

$$A_1 = \begin{bmatrix} 2 & 0 & 1 \\ 1 & 2 & 3 \\ \mathrm{I} & 0 & \mathrm{I} \end{bmatrix}$$

and

$$A_2 = \begin{bmatrix} I & .2I & .6 & .1 \\ .3 & 1 & 0 & I \\ 0 & 0 & .2 & 1 \end{bmatrix}.$$

Clearly $A_{FN}$ is a mixed semi fuzzy neutrosophic bimatrix.

***Example 3.2.18:*** Let

$$A_{FN} = \begin{bmatrix} I & .3I & 0 \\ 1 & .2 & .6 \end{bmatrix} \cup \begin{bmatrix} .6 & 0 & .3 & 1 \\ 1 & 1 & .6 & .2 \\ .3 & 0 & 0 & .5 \end{bmatrix},$$

$A_{FN}$ is a mixed semi fuzzy neutrosophic bimatrix.

***Example 3.2.19:*** Let

$$A_{FN} = \begin{bmatrix} .3 \\ I \\ .2 \\ 0 \end{bmatrix} \cup \begin{bmatrix} I \\ 7 \\ 2 \\ 1 \end{bmatrix}.$$

$A_{FN}$ is a column semi fuzzy neutrosophic bimatrix.

## 3.3 Neutrosophic bivector spaces

Now we proceed on to define the concept of neutrosophic bivector spaces. For this we define the notion of neutrosophic bigroup. Just we recall, suppose G is an additive abelian group under + . G(I) = $\langle G \cup I \rangle$ i.e., the additive group generated by G and I under '+' . G(I) is called the neutrosophic abelian group.

It is to be noted that we call (G, +, °), a bigroup if G = $G_1 \cup G_2$ and ($G_1$, +) and ($G_2$, °) are groups $G_1$ and $G_2$ are proper subsets of G i.e., $G_1 \not\subset G_2$ and $G_2 \not\subset G_1$. It is very important to note that '+' and 'o' need not in general be two distinct operations.

Now we proceed on to define neutrosophic bigroups.

**DEFINITION 3.3.1:** *The set (G(I), o, ×) is called the neutrosophic bigroup if G(I) = $G_1(I) \cup G_2(I)$ are such that $G_1(I) \not\subset G_2(I)$ or $G_2(I) \not\subset G_1(I)$ where ($G_1(I)$, o) and ($G_2(I)$, ×) are neutrosophic groups.*

***Example 3.3.1:*** $M_{n\times m}(I) = \{(a_{ij}) \mid a_{ij} \in Z(I)\}$ be the collection of all n × m matrices. $M_{n\times m}(I)$ is a group under '+', 2Z(I) be the neutrosophic group under addition G(I) $= M_{n\times m}(I) \cup 2Z(I)$ is a neutrosophic bigroup.

For more about neutrosophic groups, neutrosophic vector spaces and neutrosophic fields please refer [50, 51].

Now having defined a neutrosophic bigroup we proceed on to define neutrosophic bivector spaces.

**DEFINITION 3.3.2:** *Let K(I) be a neutrosophic field or any field. G(I) under addition be a neutrosophic bigroup, G(I) = $G_1$ (I) ∪ $G_2$ (I) is called as a bivector space over K(I) if $G_1$(I) and $G_2$ (I) are neutrosophic vector spaces over the neutrosophic field K(I) or over any field..*

Now we illustrate this by an example.

***Example 3.3.2:*** Let K be any field, V = $M_{n \times m}$ (I) ∪ K(I) × K(I) is a neutrosophic bivector space over the field K where $M_{n \times n}$ (I) are neutrosophic matrices with entries from K(I).

***Example 3.2.3:*** Let V = R(I) × R(I) ∪ $M_{3 \times 3}$ (I). Clearly V is a neutrosophic bivector space over R in fact V is a neutrosophic bivector space over R.

Having defined neutrosophic bivector space one can proceed on to define the notion of linear bitransformation and linear bioperator relative to neutrosophic bivector spaces V(I) = $V_1$ (I) ∪ $V_2$ (I) as in case of bivectors space.

We can define the dimension of the neutrosophic bivector space $V_1$ (I) ∪ $V_2$ (I) is equal to the dimension of $V_1$ (I) as a neutrosophic vector space over K and the dimension of $V_2$(I) as the neutrosophic vector space over K.

Thus the dimension of the neutrosophic bivector space V(I) = $V_1$ (I) ∪ $V_2$ (I) over K is a pair (dim $V_1$ (I), dim $V_2$ (I)). If both dim $V_1$ (I) and dim $V_2$ (I) are finite we call V(I) a finite dimensional neutrosophic bivector space, other wise infinite. This even one of the neutrosophic vector spaces $V_1$(I) or $V_2$(I) is infinite dimensional, we call V(I) = $V_1$(I) ∪ $V_2$(I) to be infinite dimensional neutrosophic bivector space. It is important to mention here that $V_1$(I) and $V_2$(I) can be of different dimensions.

**Chapter Four**

# SUGGESTED PROBLEMS

Now as such this book only views bimatrix as a theory. In fact several more results about bimatrices can also be developed.

Here we propose some problems which are simple exercises and this will make the reader more familiar with the concept of bimatrices so when the neutrosophic bimatrices and fuzzy neutrosophic bimatrices are introduced, the reader may follow it easily.

1. Find the product of $A_B$ and $C_B$ where

$$A_B = \begin{bmatrix} 3 & 1 & 1 & 1 \\ 0 & 2 & -1 & 0 \\ 3 & 1 & 0 & 2 \end{bmatrix} \cup \begin{bmatrix} 5 & 3 \\ 7 & 1 \\ 3 & 2 \end{bmatrix}$$

and

$$C_B = \begin{bmatrix} 3 & 2 \\ 4 & 2 \\ 1 & 3 \\ -1 & 5 \end{bmatrix} \cup \begin{bmatrix} 3 & 1 & 1 & 1 & 0 \\ -1 & 0 & 2 & 4 & 5 \end{bmatrix}.$$

Find $A_B$ $C_B$. Prove $C_B$ $A_B$ is not defined. Find $A_B$ + $C_B$ and prove $A_B + C_B = C_B + A_B$.

2. Let

$$A_B = \begin{bmatrix} 0 & 1 & 1 & 0 \\ -1 & 0 & 0 & 1 \\ 2 & 2 & 0 & 1 \\ 5 & 0 & 1 & -5 \end{bmatrix} \cup \begin{bmatrix} 0 & 1 & 2 \\ 3 & -1 & 3 \\ 7 & 5 & 1 \end{bmatrix}$$

be a mixed square bimatrix find $A_B^{-1}$.

Show $\left(A_B^{-1}\right)^{-1}$ = $A_B$.

3. Prove (A B)' = B'A' where

$$A_B = \begin{bmatrix} 3 & 1 & 0 & 0 \\ 0 & 1 & 2 & -1 \\ 4 & 2 & 0 & 1 \\ 9 & 2 & 0 & 0 \end{bmatrix} \cup \begin{bmatrix} 7 & 6 & 5 & 4 \\ 0 & 1 & 0 & 2 \\ 7 & 0 & 8 & 1 \\ 0 & 0 & 5 & 4 \end{bmatrix}$$

$$B_B = \begin{bmatrix} 0 & 1 & 2 & 3 \\ 1 & 2 & 3 & 0 \\ 1 & 2 & 0 & 3 \\ 0 & 1 & 2 & 1 \end{bmatrix} \cup \begin{bmatrix} 5 & 4 & 6 & 3 \\ 3 & 5 & 4 & 6 \\ 6 & 3 & 5 & 4 \\ 4 & 6 & 5 & 3 \end{bmatrix}.$$

4. Find $A'_B + B'_B$ where $A_B$ and $B_B$ are given as in problem (3). Is $(A_B + B_B)' = A'_B + B'_B$?

5. Let

$$A_B = \begin{bmatrix} 3 & -2 & 4 & 0 \\ 9 & 0 & 1 & 1 \\ 8 & 1 & 0 & 1 \\ 6 & 0 & 0 & 5 \end{bmatrix} \cup \begin{bmatrix} 9 & 0 & 1 & 1 \\ 12 & 1 & 0 & -1 \\ 6 & 0 & 1 & 1 \\ 8 & 0 & 1 & 0 \end{bmatrix}$$

using partitioning of bimatrices find the bideterminant value of $A_B$.

6. Using biLaplace expansion expand the bimatrix $A_B$ where

$$A_B = \begin{bmatrix} 8 & 0 & 0 & 0 & 6 \\ 9 & 5 & 1 & 0 & 2 \\ 12 & 0 & 1 & 0 & 1 \\ 0 & 1 & 2 & 1 & 0 \\ 0 & 8 & 0 & 5 & 0 \end{bmatrix} \cup \begin{bmatrix} 9 & 2 & 0 & 6 & 1 \\ -1 & 3 & 8 & -1 & 4 \\ 9 & 2 & -1 & 0 & 2 \\ 2 & 2 & 4 & 0 & 3 \\ 0 & 3 & 1 & 0 & 6 \end{bmatrix}.$$

7. Write the bimatrix $A_B$ given in problem (6) as the sum of a skew symmetric and symmetric bimatrices.

8. Let $V = V_1 \cup V_2$ be an (n, p) dimensional bivector space over the field F and let $W = W_1 \cup W_2$ be an (m, q) dimensional bivector space over F. For each pair of ordered bases $(B_1, B_2)$ for $V = V_1 \cup V_2$ and $(B_1^1, B_2^1)$ of $W = W_1 \cup W_2$ respectively, prove the collection of all the functions which assigns to a linear bitransformation $T = T_1 \cup T_2$, a bimatrix relative is an isomorphism between the bivector space $L_F$ $(V = V_1 \cup V_2, W = W_1 \cup W_2)$ and the space of $(m \times n, q \times p)$ bimatrices over the field F.

9. Let $V = V_1 \cup V_2$ be a finite dimensional vector space over the field F and let

$$B = \left\{B_1 = \left\{\alpha_1^1, \ldots, \alpha_n^1\right\},\ B_2 = \left\{\alpha_1^2, \ldots, \alpha_n^2\right\}\right\}$$

be an ordered basis of $V = V_1 \cup V_2$. Derive a change of basis rule for V.

10. Let V and W be bivector spaces over the field F and let T be a linear bitransformation from V into W and suppose V is finite dimensional, then prove rank (T) + nullity (T) = dim V = (dim $V_1$, dim $V_2$).

11. Prove if V and W are bivector spaces over the same field F, then the set of all linear bitransformations from V into W together with addition and scalar multiplication as defined in this book is a bivector space over the field F.

    (Note we shall denote the set of all linear bitransformation of the bivector spaces V and W over F by $L^B_F$ (V, W)).

12. Let V be a n dimensional bivector space over the field F and let W be a m dimensional bivector space over F. Let $L^B_F$ (V, W) be a finite dimensional bivector space, what is the dimension of $L^B_F$ (V, W)?

13. Prove if V is a bivector space over F, then the set of all linear bioperators on V denoted by $L^B_F$ (V, V) is a bivector space over F.

14. Prove if $T = T_1 \cup T_2$ from $V = V_1 \cup V_2$ on to $W = W_1 \cup W_2$ is an invertible linear bitransformation and $U = U_1 \cup U_2$ from $W = W_1 \cup W_2$ to $Z = Z_1 \cup Z_2$ is an invertible linear Bitransformation, then

    $$(UT) = (U_1T_1) \cup (U_1 T_1)$$

    is invertible and

    $(UT)^{-1} = T^{-1} U^{-1}$ .

    (Note $(UT)^{-1} = [(U_1 \cup U_2) (T_1 \cup T_2)]^{-1}$
    $= (U_1T_1)^{-1} \cup (U_2 T_2)^{-1}$
    $= T_1^{-1}U_1^{-1} \cup T_2^{-1}U_2^{-1}$
    $= T^{-1} U^{-1}$ ).

15. If T is a linear bitransformation from $V = V_1 \cup V_2$ into $W = W_1 \cup W_2$, then T is non singular if and only if T carries each linearly independent subset of $V = V_1 \cup V_2$ on to a linearly independent subset of $W = W_1 \cup W_2$.

16. For a linear bioperator $T = T_1 \cup T_2$ prove a theorem analogous to primary decomposition theorem. Hence or otherwise prove that any linear bioperator T can be written as a sum of a diagonal bioperator and a nilpotent bioperator.

17. Prove if T is a linear bioperator on a finite dimensional bivector space V over the field F, then T = D + N, DN = ND where D is a diagonalizable bioperator on V and N is a nilpotent bioperator on V.

18. Obtain any interesting result on linear bioperator.

19. Define semi simple linear bioperator. Illustrate it with example.

20. Let T be a linear bioperator on a finite dimensional bivector space $V = V_1 \cup V_2 = Q^6 \cup Q^4$, $T_1$ (x, y, z, w, t, u) = (x – y, y – z, z – w, w – t, t – u, u – x) and $T_2$ (x, y, z, w) = (x + y, z + w, z, x). Is $T = T_1 \cup T_2$ diagonalizable?

21. Let T: $V = V_1 \cup V_2 \to V = V_1 \cup V_2$ be given by the bimatrix $A = A_1 \cup A_2$ relative to the standard basis i.e.,

$$A = \begin{bmatrix} 9 & 0 & 0 & -1 \\ 0 & 2 & 0 & 0 \\ 0 & 0 & 5 & 1 \\ -1 & 0 & 0 & 0 \end{bmatrix} \cup \begin{bmatrix} 1 & 0 & 0 \\ 2 & 1 & 0 \\ 3 & 0 & 2 \end{bmatrix}.$$

Find null bispace of T.

22. Is T in problem 21 diagonalizable? Find the characteristic bivalues of T.

23. Let T : V → V where V = $Q^3 \cup Q^5$ be a linear bioperator. If the related or associated bimatrix of T with the standard basis be

$$A = \begin{bmatrix} 3 & 0 & 0 \\ 0 & 2 & 0 \\ -1 & 1 & 1 \end{bmatrix} \cup \begin{bmatrix} -1 & 0 & 0 & 0 & 0 \\ 6 & 2 & 0 & 0 & 0 \\ 1 & 1 & 1 & 0 & 0 \\ 0 & 0 & 1 & 2 & 0 \\ 1 & 1 & 0 & -1 & -4 \end{bmatrix}.$$

Does there exist linear bioperator $E_i$ on V such that $T = \sum C_i E_i$ i.e., if $T = T_1 \cup T_2$ can $T_1$ and $T_2$ be written as

$$\sum_i C_i^1 E_i^1 \text{ and } \sum_i C_i^2 E_i^2 \ ?$$

Justify your claim.

24. Let V be a bivector space over reals, T a linear bioperator on V. Let the related bimatrix of $T = T_1 \cup T_2$ with respect to the standard basis is given by

$$A = A_1 \cup A_2 = \begin{bmatrix} 3 & 0 & 0 & 1 \\ 0 & 2 & 1 & 0 \\ 0 & 1 & 3 & 0 \\ 1 & 0 & 0 & -1 \end{bmatrix} \cup \begin{bmatrix} 6 & 0 & 2 \\ 0 & 3 & 1 \\ 2 & 1 & 5 \end{bmatrix}.$$

Does there exist a D-diagonalizable part and N nil potent part of T such that T = D + N and ND = ND?

25. Define Jordan biform for bimatrices, illustrate with examples.

26. Find the Jordan biform of the strong bimatrix $A = A_1 \cup A_2$, where

$$A = \begin{bmatrix} 0 & 1 & 0 & 0 \\ 0 & 0 & 2 & 0 \\ 0 & 0 & 0 & 3 \\ 0 & 0 & 0 & 0 \end{bmatrix} \cup \begin{bmatrix} 2 & 0 & 0 & 0 & 0 & 0 \\ 1 & 2 & 0 & 0 & 0 & 0 \\ -1 & 0 & 2 & 0 & 0 & 0 \\ 0 & 1 & 0 & 2 & 0 & 0 \\ 1 & 1 & 1 & 1 & 2 & 0 \\ 0 & 0 & 0 & 0 & 1 & -1 \end{bmatrix}.$$

(Here $A_1$ is got as the differentiation operator on the space of polynomials of degree less than or equal to 3 and is represented in the 'natural' ordered basis by the matrix $A_1$ and $A_2$ is the complex matrix, with entries from the complex field.)

27. If $A = A_1 \cup A_2$ is a bimatrix with entries from the complex field. $A_1$ is a $5 \times 5$ complex matrix with characteristic polynomial $f_1 = (x-2)^3 (x+7)^2$ and minimal polynomial $p = (x-2)^2 (x+7)$ and $A_2$ is a $6\times 6$ complex matrix with characteristic polynomial $f_2 = (x+2)^4 (x-1)^2$. Find the Jordan biform of $A = A_1 \cup A_2$. How many possible Jordan biforms exist for A?

28. Find the biminimal polynomials of the following bimatrices with real entries

a. $$A = \begin{bmatrix} 0 & -1 & 1- \\ 1 & 0 & 0 \\ -1 & 0 & 0 \end{bmatrix} \cup \begin{bmatrix} c & 0 & -1 \\ 0 & c & 1 \\ -1 & 1 & c \end{bmatrix},$$

b. $$A_1 = \begin{bmatrix} 3 & -4 & -4 \\ -1 & 3 & 2 \\ 2 & -4 & -3 \end{bmatrix} \cup \begin{bmatrix} \text{Cos}\,\theta & \text{Sin}\theta \\ -\text{Sin}\theta & \text{Cos}\theta \end{bmatrix}.$$

29. Is the bimatrix

$$A = \begin{bmatrix} 1 & 3 & 3 \\ 3 & 1 & 3 \\ -3 & -3 & -5 \end{bmatrix} \cup \begin{bmatrix} 2 & 0 & 0 & 0 \\ 1 & 2 & 0 & 0 \\ 0 & 1 & 2 & 0 \\ 0 & 0 & 1 & 2 \end{bmatrix}$$

invertible?
Find $B = B_1 \cup B_2$ such that B is similar to A.

30. Prove primary decomposition theorem in case of finite dimensional neutrosophic bivector spaces.

31. Obtain an example of a bidiagonalizable linear bioperator for a neutrosophic bivector space.

32. Given

$$A_N = \begin{bmatrix} 0 & I & 0 \\ 2 & -2 & 2 \\ 2 & -3 & 2 \end{bmatrix} \cup \begin{bmatrix} I & 2 \\ 0 & 2 \end{bmatrix}$$

is a neutrosophic bimatrix. Is A similar to a neutrosophic triangular bimatrix?

33. Given

$$A_N = \begin{bmatrix} 3 & I & 0 \\ 1 & I & 2 \\ 0 & 0 & 1 \end{bmatrix} \cup \begin{bmatrix} 0 & I \\ 2 & 1 \end{bmatrix}$$

find the characteristic bivalues and bivectors associated with $A_N$. Are the bivalues neutrosophic? Are the bivectors neutrosophic bivectors?

34. Prove some interesting results on linear bioperators in case of neutrosophic bivector spaces.

35. Find the minimal bipolynomial associated with the neutrosophic bimatrix

$$A_N = \begin{bmatrix} I & I & 0 & 0 \\ -1 & -1 & 0 & 0 \\ -2 & -2 & 2 & 1 \\ 1 & 1 & -1 & 0 \end{bmatrix} \cup \begin{bmatrix} 3 & 2 & 0 \\ 1 & 0 & 1 \\ 0 & 2 & -1 \end{bmatrix}.$$

Is the minimal bipolynomial a neutrosophic bipolynomial? (Note: A polynomial is a neutrosophic polynomial if the coefficient field over which the polynomial is defined is a neutrosophic field).

36. Can you find $a_1$, $b_1$ and $a_2$, $b_2$ in the neutrosophic bimatrix

$$A_N = \begin{bmatrix} 0 & 0 & 0 & 0 \\ a_1 & 0 & 0 & 0 \\ 0 & b_1 & 0 & I \\ 0 & 0 & 0 & 1 \end{bmatrix} \cup \begin{bmatrix} 0 & 0 & 0 \\ a_2 & 0 & I \\ 0 & b_2 & I \end{bmatrix}$$

using the sole condition $A_N$ is a diagonalizable neutrosophic bimatrix.

37. Let

$$A_N = \begin{bmatrix} 0 & I & 1 \\ -2 & 4I & 0 \\ 1 & 1 & 2 \end{bmatrix} \cup \begin{bmatrix} I & 0 & 0 & 0 \\ 2 & -1 & 0 & 2 \\ 1 & I & 0 & 1 \\ 1 & 1 & 0 & I \end{bmatrix}$$

be a neutrosophic bimatrix over Q (I). Does $A_N$ have characteristic bivalues? If $A_N$ has characteristic bivalues find the characteristic bivectors associated with it.

38. Find the associated neutrosophic bimatrix for the linear bioperator $T = T_1 \cup T_2$, on the neutrosophic bivector space V over Q(I) given by

$$V = Q(I) \times Q(I) \cup Q(I) \times Q(I) \times Q(I) \times Q(I)$$

defined by $T = (T_1 \cup T_2)$ $[(x, y), (x_1\ x_2\ x_3\ x_4)] = [(x + y, Iy), (x_1 + x_2, x_2 - x_3\ Ix_4, 2x_2)]$.

39. Given $V = V_1 \cup V_2$ is a neutrosophic bivector space over Q(I) where $V_1 = Q(I) \times Q(I) \times Q(I) \times Q(I)$ and $V_2 = M_{3 \times 4}$ (I) is a $3 \times 4$ neutrosophic matrix with entries from Q(I) and $W = W_1 \cup W_2$ where $W_1 = Q(I) \times Q(I) \times Q(I)$ and $W_2 = Q(I) \times Q(I) \times Q(I) \times Q(I)$ is the neutrosophic vector space over Q(I) Find two linear bitransformation $T: V \to W$ and $U : V \to W$ Prove $(T + U) : V \to W$ is a linear bitransformation from V to W. When will TU be defined?

40. Let $V (I) = V_1 \cup V_2$ be a neutrosophic bivector space defined over Q(I); where $V_1 = Q(I) \times Q(I) \times Q(I)$ and $V_2 = M_{4 \times 4}$ (I) i.e., set of all $4 \times 4$ matrices with entries from Q(I). Find a basis $B = (B_1, B_2)$ for V(I).

41. Define the following on a neutrosophic bivector space V(I).
    a. A bidiagonalizable linear bioperator T on a neutrosophic bivector space.
    b. Define a linear bioperator T on V(I) such that T o T = T.
    c. Define a linear bioperator on V(I) so that T is invertible.

# BIBLIOGRAPHY

Here we suggest the following book which my help the reader for further development of bimatrices.

# INDEX

**B**

**O**

**P**

**R**

**S**

**T**

**U**

**W**

**Z**

# About the Authors

**Dr.W.B.Vasantha Kandasamy** is an Associate Professor in the Department of Mathematics, Indian Institute of Technology Madras, Chennai, where she lives with her husband Dr.K.Kandasamy and daughters Meena and Kama. Her current interests include Smarandache algebraic structures, fuzzy theory, coding/ communication theory. In the past decade she has guided eight Ph.D. scholars in the different fields of non-associative algebras, algebraic coding theory, transportation theory, fuzzy groups, and applications of fuzzy theory of the problems faced in chemical industries and cement industries. Currently, six Ph.D. scholars are working under her guidance. She has to her credit 285 research papers of which 209 are individually authored. Apart from this, she and her students have presented around 329 papers in national and international conferences. She teaches both undergraduate and post-graduate students and has guided over 45 M.Sc. and M.Tech. projects. She has worked in collaboration projects with the Indian Space Research Organization and with the Tamil Nadu State AIDS Control Society. She has authored a Book Series, consisting of ten research books on the topic of Smarandache Algebraic Structures which were published by the American Research Press.

She can be contacted at vasantha@itm.ac.in
You can visit her work on the web at: http://mat.iitm.ac.in/~wbv

---

**Dr.Florentin Smarandache** is an Associate Professor of Mathematics at the University of New Mexico, Gallup Campus, USA. He published over 60 books and 80 papers and notes in mathematics, philosophy, literature, rebus. In mathematics his research papers are in number theory, non-Euclidean geometry, synthetic geometry, algebraic structures, statistics, and multiple valued logic (fuzzy logic and fuzzy set, neutrosophic logic and neutrosophic set, neutrosophic probability). He contributed with proposed problems and solutions to the Student Mathematical Competitions. His latest interest is in information fusion were he works with Dr.Jean Dezert from ONERA (French National Establishment for Aerospace Research in Paris) in creasing a new theory of plausible and paradoxical reasoning (DSmT).

He can be contacted at smarand@unm.edu

---

**K. Ilanthenral** is the editor of The Maths Tiger, Quarterly Journal of Maths. She can be contacted at ilanthenral@gmail.com